\documentclass{elsarticle}

\usepackage[left=1in,top=1in,right=1in,bottom=1in]{geometry}

\usepackage{amsmath}
\usepackage{amssymb}
\usepackage{placeins}
\usepackage{subfig}

\begin{document}

\begin{frontmatter}

\title{Stabilized Isogeometric Collocation Methods For Scalar Transport and Incompressible Fluid Flow}

\author[1]{Ryan M. Aronson}
\ead{rmaronso@stanford.edu}

\author[2]{Corey Wetterer-Nelson}
\ead{c.wetterer-nelson@kitware.com}

\author[3]{John A. Evans\corref{cor1}}
\ead{john.a.evans@colorado.edu}

\cortext[cor1]{Corresponding author}

\affiliation[1]{organization={Stanford University},
city={Stanford, CA},
postcode={94305},
country={USA}}

\affiliation[2]{organization={Kitware Inc.},
city={Clifton Park, NY},
postcode={12065},
country={USA}}

\affiliation[3]{organization={University of Colorado Boulder},
city={Boulder, CO},
postcode={80309},
country={USA}}



\begin{abstract}
    In this work we adapt classical residual-based stabilization techniques to the spline collocation setting. Inspired by the Streamline-Upwind-Petrov-Galerkin and Pressure-Stabilizing-Petrov-Galerkin methods, our stabilized collocation schemes address spurious oscillations that can arise from advection and pressure instabilities. Numerical examples for the advection-diffusion equation, Stokes equations, and incompressible Navier-Stokes equations show the effectiveness of the proposed stabilized schemes while maintaining the high-order convergence rates and accuracy of standard isogeometric collocation on smooth problems.
\end{abstract}


\end{frontmatter}

\section{Introduction}

\textit{Tom Hughes has been nothing short of an ideal mentor, role model, and friend to me (John Evans) since I started my doctoral studies under his supervision in August 2006.  I certainly would not be where I am today without him.  As such, it is my pleasure to dedicate this paper to the celebration of Tom’s 80th birthday.  This paper synthesizes Tom’s two greatest contributions to computational mechanics, isogeometric analysis (specifically isogeometric collocation) and residual-based stabilized methods.  Happy birthday Tom!}

Isogeometric Analysis (IGA) is an alternative to classical finite element analysis in which $C^0$-continuous piecewise polynomial basis functions are replaced by spline basis functions \cite{hughes_IGA, cottrell_IGAbook}. IGA has a number of advantages over classical finite element analysis, including the possibility for exact geometry representation independent of the mesh resolution for geometries represented by B-splines and Non-Uniform Rational B-splines (NURBS) as well as improved approximation properties \cite{evans_nwidths, sande2020explicit}. Moreover, the increased global smoothness of spline basis functions enables the use of collocation as an alternative to Galerkin's method. This approach was first studied in the context of IGA in \cite{auricchio_isogeometric_coll, reali_intro} and has since been successfully applied to many problems in solid mechanics \cite{auricchio_coll_elastostatics, evans_explicit, kruse2015largedef}. Many of the potential cost advantages of collocation are explained in  \cite{schillinger_cost_comp} and are largely due to the fact that collocation requires no spatial numerical integration steps. Instead, the strong form of the partial differential equation is evaluated at a set of collocation points to solve for the solution represented as a linear combination of spline basis functions.

In the fluid mechanics community, B-spline discretization methods were explored even before the invention of IGA. For example, \cite{kravchenko1999b} utilized a B-spline Galerkin method for the simulation of incompressible turbulent flows, while \cite{botella_collocation} used a mixed spline collocation technique to solve incompressible flow problems with the goal of reducing the cost of the evaluation of the nonlinear terms when compared to Galerkin methods. More recently, B-spline collocation has been used in conjunction with the Fourier-Galerkin method for direct numerical simulations of turbulent channel flows \cite{lee2015DNS}. In the IGA context, collocation solutions of scalar transport problems were considered in \cite{schillinger_cost_comp}. The authors of the current paper also developed divergence-conforming B-spline collocation schemes which exactly satisfy the incompressibility constraint at a pointwise level \cite{aronson2023divergence}.

A topic that has not been studied in detail is the stabilization of spline collocation methods for transport problems. Spurious oscillations in collocation solutions of flow problems can arise from both advection instabilities as well as pressure instabilities. Much like central difference schemes and unstabilized Galerkin finite element methods, spurious oscillations in the solution field (in the case of scalar transport) or the velocity field (in the case of incompressible flow) can appear upstream of sharp flow features, which is what we refer to as an advection instability. Pressure instabilities, on the other hand, manifest as spurious oscillations in the pressure field resulting from a poor choice of velocity and pressure discretization spaces. In the worst case scenario, the pressure solution is not unique due to the presence of spurious modes, and the velocity field can exhibit locking if over-constrained. Some preliminary work was done in the area of advection stabilization for IGA collocation in \cite{schillinger_cost_comp}, where, inspired by \cite{shapiro1981analysis, funaro1999spline}, the authors moved the collocation points in the upstream direction in order to lessen the magnitude of spurious velocity oscillations in highly advective linear transport problems. However, it is not clear how far the collocation points should be moved in this case, and moving points too far could result in collocation points that lie outside the physical domain. This scheme is also difficult to extend to nonlinear problems, where the local advection velocity can change between nonlinear iterations or time steps, as collocation points would need to be recomputed at every iteration. This in turn means that the point values of basis functions and their derivatives must be recomputed at every iteration, and thus the discrete, matrix form of the system is altered. It is possible that the set of basis functions which are nonzero at a collocation point may change between iterations, meaning the sparsity pattern of the system tangent matrix could change between iterations as well. Moreover, the idea of selecting upwind collocation points does not yield an obvious scheme to resolve pressure instabilities, meaning that carefully selected spaces must always be used for mixed problems. This precludes the use of equal-order interpolations, for example, which are convenient for implementation.

Meanwhile, many stabilized Galerkin finite element methods have been developed to alleviate advection and pressure instabilities and have been employed with great success in incompressible flow problems. The most classical stabilized methods are residual-based methods \cite{tezduyar1991stabilized, Brooks_supg, hughesPSPG, hughes1989gls}, but there are many other stabilization approaches, including edge or interior penalty stabilization \cite{burman2004edge, burman2007continuous} and local projection stabilization \cite{guermond1999stabilization, braack2006local}. Moreover, many approaches have been unified under the umbrella of variational multiscale (VMS) methods \cite{hughes1995multiscale, hughes1998variational}.

In this work we develop stabilized methods within the collocation framework inspired by residual-based stabilized finite element methods. We show by adapting the classical Streamline-Upwind-Petrov-Galerkin (SUPG) stabilization first introduced in \cite{Brooks_supg} to collocation we can effectively suppress oscillations in advection dominated flows. We also show that the classical Pressure-Stabilizing-Petrov-Galerkin (PSPG) stabilization \cite{hughesPSPG} can be written in a collocation framework and is effective at removing spurious pressure oscillations for equal-order velocity and pressure discretizations of incompressible Stokes and Navier-Stokes flows. Moreover, the use of the same approximation spaces for velocity and pressure means that the same set of collocation points can be used for both the momentum and continuity equations, further improving the simplicity and efficiency of the resulting method. When collocating the Navier-Stokes equations we also adapt Grad-Div stabilization \cite{olshanskii2009grad} to the collocation setting, which improves pointwise mass conservation. 

The rest of the manuscript is ordered as follows. We start with a description of collocation schemes, including appropriate choices of spline basis functions and associated collocation points. Then, in Section 3, we discuss the advection-diffusion equation and use it as motivation in developing advection-stabilized collocation schemes. Section 4 proceeds to the topic of incompressible Stokes flow, which is used to show the effectiveness of the pressure-stabilized collocation scheme. In this section, we also investigate the effect the differential order of the governing equations and we show that the stabilized scheme also achieves accelerated convergence rates when the rotational form of the Stokes equations are used, as seen in \cite{aronson2023divergence}. Finally, Section 5 develops stabilized collocation schemes for the incompressible Navier-Stokes equations which combine the advection and pressure stabilization techniques developed prior. A suite of benchmark problems shows the effectiveness of the stabilized methods in suppressing oscillations while retaining high-order convergence rates, and once again collocation schemes based on the rotational form of the Navier-Stokes equations achieve faster convergence rates relative to schemes based on the standard form of the equations. 

\section{Collocation Methods}

We start by describing collocation methods for solving partial differential equations (PDEs) in general, before describing the specifics of isogeometric collocation methods including appropriate choices of basis functions and collocation points.

\subsection{Collocation of the Strong Form}

The first description of collocation stems directly from the strong form of a partial differential equation. Consider a domain $\Omega$ with boundary $\partial \Omega$ over which an unknown field variable $\phi$ is defined. This variable satisfies a partial differential equation and boundary conditions

\bigskip

$$
\left\{ \hspace{5pt}
\parbox{5in}{
\noindent Given $\mathcal{L}$, $f$, $\mathcal{G}$, and $g$, find $\phi$ such that:
\begin{equation}
    \mathcal{L\{\phi\}} = f \textup{ in } \Omega,
\end{equation}
\begin{equation}
    \mathcal{G\{\phi\}} = g \textup{ on } \partial \Omega,
\end{equation}
}
\right.
$$

\bigskip

\noindent where $\mathcal{L}$ is the differential operator, $f$ is a forcing term, $\mathcal{G}$ is the boundary operator, and $g$ is the provided boundary data. 

Collocation methods approximate the solution of this problem with a linear combination of $n$ basis functions $\phi^h(\mathbf{x}) = \sum_{i = 1}^n N_i(\mathbf{x}) c_i$, where $\{N_i\}_{i=1}^{n}$ are the basis functions and $\{c_i\}_{i=1}^{n}$ are the unknown coefficients. In order to solve for the unknown coefficients we select $n$ collocation points $\{\mathbf{x}\}_{i=1}^n$ in the domain at which we will interpolate either the PDE (if the point is on the interior of the domain) or the boundary condition data (if the point is on the boundary):

\bigskip

$$
\left\{ \hspace{5pt}
\parbox{5in}{
\noindent Given $\mathcal{L}$, $f$, $\mathcal{G}$, and $g$, find $\phi^h(\mathbf{x}) = \sum_{i = 1}^n N_i(\mathbf{x}) c_i$ such that:
\begin{equation}
    \mathcal{L}\{\phi^h\}(\mathbf{x}_i) = f(\mathbf{x}_i) \quad \forall \mathbf{x}_i \in \Omega,
\end{equation}
\begin{equation}
    \mathcal{G}\{\phi^h\}(\mathbf{x}_i) = g(\mathbf{x}_i) \quad \forall \mathbf{x}_i \in \partial \Omega.
\end{equation}
}
\right.
$$

\bigskip

An alternative and equivalent way to construct collocation schemes is inspired by the weak formulation of the PDE. Using the set of basis functions defined above, we denote the trial function space as $S^h = \textup{span}\{N_i\}_{i=1}^{n}$. The space $S^h_g$ is defined as the set of all functions in $S^h$ which appropriately satisfy any provided Dirichlet boundary conditions. For a Galerkin method we also define $S^h_0$ to be the finite-dimensional test function space of all functions in $S^h$ which satisfy homogeneous boundary conditions along the specified Dirichlet boundaries. The weighted residual form of the PDE is then:

\bigskip

$$
\left\{ \hspace{5pt}
\parbox{5in}{
\noindent Given $\mathcal{L}$ and $f$, find $\phi^h \in S^h_g$ such that:
\begin{equation}
    \int_\Omega w^h (\mathcal{L}\{\phi^h\} - f) d \Omega = 0 \quad \forall w^h \in S^h_0.
\end{equation}
}
\right.
$$

\bigskip

If we instead select the test functions $w^h$ to be Dirac-delta functions located at the collocation points $\{\mathbf{x}\}_{i=1}^n$ selected above, we arrive at the same collocation scheme as generated by the strong form of the PDE. In this way, collocation can also be thought of as a Petrov-Galerkin scheme. This equivalence is explored in more detail in, for example,  \cite{gomez2016variational, deLorenzis2015isogeometric}. We will use this weak form motivation several times when developing stabilized collocation schemes. 

\subsection{B-splines and Isogeometric Analysis}

A full definition of a collocation scheme requires the appropriate definition of the basis functions $\{N_i\}_{i=1}^{n}$ and the collocation points $\{\mathbf{x}\}_{i=1}^n$. The basis functions selected must possess sufficient smoothness such that the residual function $\mathcal{L}\{\phi^h\} - f$ is continuous. For second-order operators $\mathcal{L}$, the basis must be globally $C^2$-continuous. One set of basis functions that can be constructed to satisfy this requirement are B-spline basis functions, which we utilize for our collocation schemes. 

B-spline basis functions in the one-dimensional setting are defined by a polynomial order $k$ and a knot vector $\Xi = \{\xi_1, ... \xi_{n+p+1}\}$, which defines the parametric locations of possible reduced continuity of the global basis functions, similar to element boundaries in standard finite element methods. The basis functions are defined through the Cox-de Boor recursion: The $k = 0$ basis functions are built as

\begin{equation}
    N_{i,0}(\xi) = \begin{cases} 
                1 & \xi_i \leq \xi \leq \xi_{i+1} \\
                0 & \text{otherwise},
              \end{cases}
\end{equation}

\noindent and higher-order bases are defined through

\begin{equation}
    N_{i,k}(\xi) = \frac{\xi - \xi_i}{\xi_{i+k} - \xi_i}N_{i,k-1}(\xi) + \frac{\xi_{i+k+1} - \xi}{\xi_{i+k+1} - \xi_{i+1}}N_{i+1,k-1}(\xi).
\end{equation}

At any knot in the knot vector $\Xi$, the basis functions are $C^{k-\ell}$-continuous, where $\ell$ is the multiplicity of the knot. This is the key property that enables collocation: higher order spline basis functions also have increased levels of continuity and this level of continuity is controllable. In this work we will always utilize maximal continuity splines, meaning that interior knots are never repeated in the knot vector $\Xi$. 

B-spline basis functions in higher spatial dimensions are generated by simply taking the tensor product of one dimensional B-spline bases in each direction, possibly using different polynomial degrees and knot vectors \cite{cottrell_IGAbook}. For simplicity, however, we will only consider spline bases with the same polynomial degrees and knot vectors in each direction in this work, meaning that only one degree and one knot vector define a spline space $S^h$ in any dimension. 

\subsection{Collocation at Greville Abscissae}

With basis functions defined, all that remains is the selection of appropriate collocation points. In this paper we will use the Greville points, defined in 1D as 

\begin{equation}
    \hat{\xi_i} = \frac{\xi_{i+1} + ... + \xi_{i+k}}{k},
\end{equation}

\noindent and in higher dimensions as the tensor product of the Greville points in each direction \cite{auricchio_isogeometric_coll, johnson_Greville_coll}. Note that by construction the number of Greville points for a set of spline basis functions is the same as the number of basis functions.

While the Greville points are quite simple to work with and give quality results in practice \cite{auricchio_isogeometric_coll, johnson_Greville_coll}, they have not been proved to yield stable interpolations in general and there are examples of unstable interpolations using non-uniform knot vectors and high polynomial degrees greater than or equal to 20 \cite{jia1988spline}. There are other possible sets of points at which the PDE could be interpolated. The Demko abscissae, for example, are provably stable for any knot vector and polynomial degree \cite{demko_abs}, but the points must be computed iteratively. The Cauchy-Galerkin points are also an option \cite{montardini_optimal, anitescu2015isogeometric}, though these points depend on the PDE under consideration and generally are not known \textit{a priori}.

\section{Advection-Diffusion and Scalar Transport}

Having described the basics of collocation methods, including basis function and collocation point definitions, we can now move on to consider collocation methods for flow and transport. We start with the simple case of scalar advection, which will allow us to describe a method for advection stabilization. 

\subsection{The Strong Form of Advection-Diffusion}

The strong form of the steady Advection-Diffusion equation subjected to Dirichlet boundary conditions is given by

\bigskip

$$
\left\{ \hspace{5pt}
\parbox{5in}{
\noindent Given $\kappa \in \mathbb{R}^+$, $\textbf{f} : \Omega \rightarrow \mathbb{R}^n$, $\textbf{u} : \Omega \rightarrow \mathbb{R}^n$, and $\textbf{g} : \partial \Omega \rightarrow \mathbb{R}^n$, find $\phi : \Omega \rightarrow \mathbb{R}$ such that:
\begin{align}
    \textbf{u} \cdot \nabla \phi - \kappa \Delta \phi = f \quad &\textup{in} \quad \Omega, \\
    \textbf{u} = \textbf{g} \quad &\textup{on} \quad \partial \Omega,
\end{align}
}
\right.
$$

\bigskip

\noindent where $\phi$ is the unknown being transported, $\textbf{u}$ is the advection velocity, $\kappa$ is the diffusion coefficient, and $\textbf{g}$ is the provided Dirichlet boundary data. The Péclet number is defined as $Pe = \frac{||u|| L}{\kappa}$, where $L$ is a length scale defined for the problem under consideration. This dimensionless quantity measures the relative strength of advection and diffusion in the problem, with $Pe > 1$ representing an advection-dominated flow and $Pe < 1$ representing a diffusion-dominated problem. It is well known that strongly advection-dominated flows can develop sharp features which can cause oscillations in solutions obtained with coarse, central-like discretizations (meaning those without any sense of upwinding). 

\subsection{SUPG Stabilization for the Advection-Diffusion Equation}

In Galerkin finite element methods these types of oscillations can be smoothed using the classical SUPG stabilization \cite{Brooks_supg}. To inspire our stabilized collocation scheme we start by stating an isogeometric Galerkin formulation including SUPG stabilization. We define spline spaces $S^h$, $S^h_g$, and $S^h_0$ as before, with $S^h_g$ satisfying the provided Dirichlet boundary data and $S^h_0$ satisfying homogeneous boundary conditions along the Dirichlet boundary. We assume that the basis functions are globally $C^1$ or smoother. Then the stabilized weak form is given by

\bigskip

$$
\left\{ \hspace{5pt}
\parbox{5in}{
\noindent Given $\kappa \in \mathbb{R}^+$, $\textbf{f} : \Omega \rightarrow \mathbb{R}^n$, $\textbf{u} : \Omega \rightarrow \mathbb{R}^n$, and $\textbf{g} : \partial \Omega \rightarrow \mathbb{R}^n$, find $\phi^h \in S^h_g$ such that, $\forall w^h \in S^h_0$ :

\begin{equation}
\begin{split}
    \int_{\Omega} (- \nabla w ^h\cdot \textbf{u} \phi^h + \kappa \nabla w^h \cdot \nabla \phi - w^hf) d\Omega\\ +  \int_{\Omega} \tau \textbf{u} \cdot \nabla w^h (\textbf{u} \cdot \nabla \phi^h - \kappa \Delta \phi^h - f) d\Omega = 0 ,
\label{eq:adv_diff_weak}
\end{split}
\end{equation}

}
\right.
$$

\bigskip

\noindent where $\tau$ is a stabilization constant to be defined.

To arrive at a form more amenable to collocation, we take Equation \eqref{eq:adv_diff_weak} and reverse integrate by parts once to reduce the number of derivatives on $w^h$:

\begin{equation}
    \int_{\Omega} w^h( \textbf{u} \cdot \nabla \phi^h - \kappa \Delta \phi^h - f - \nabla \cdot (\tau \textbf{u}  (\textbf{u} \cdot \nabla \phi^h - \kappa \Delta \phi^h - f))) d\Omega = 0.
    \label{eq:adv_diff_wr}
\end{equation}

Inspired by this formulation, we define a stabilized collocation scheme by defining the test functions $w^h$ in Equation \eqref{eq:adv_diff_wr} to be Dirac-delta functions located at the collocation points, as described previously. We let the discrete solution $\phi^h$ be a member of the spline space $S^h$ with dimension $n$ and Greville abscissae $\textbf{x}_i$ for $i = 1, ..., n$. Then the stabilized collocation scheme is defined by

\bigskip

$$
\left\{ \hspace{5pt}
\parbox{5in}{
\noindent Find $\phi^h \in S^h$ such that

\begin{equation}
\begin{split}
    &( \textbf{u} \cdot \nabla \phi^h - \kappa \Delta \phi^h - f \\ -  &\nabla \cdot (\tau \textbf{u}  (\textbf{u} \cdot \nabla \phi^h - \kappa \Delta \phi^h - f)) ) (\textbf{x}_i) = 0 \quad \forall \textbf{x}_i \in \Omega,
\end{split}
\end{equation}
\begin{equation}
    \phi^h(\textbf{x}_i) = g(\textbf{x}_i) \quad \forall \textbf{x}_i \in \partial \Omega. \label{eq:adv_diff_BC}
\end{equation}
}
\right.
$$

\bigskip

All that remains is to define the stabilization parameter $\tau$ at every collocation point. For this we utilize a standard definition \cite{tezduyar2003stabilization, tezduyar2000finite}:

\begin{equation}
    \tau = \frac{1}{\sqrt{(\frac{2 ||\mathbf{u}||}{h})^2 + (\frac{C_1 \kappa}{h^2})^2}},
\end{equation}

\noindent where, following \cite{bazilevs2007weak} $C_1$ is set to be 4 and $h$ is the mesh size defined by the average distance from each collocation point to each of its neighboring collocation points. We then use a spline interpolation in the same space $S^h$ as the discrete solution to interpolate the $\tau$ values, thus defining the derivatives of the stabilization parameter.

\subsubsection{A Note on Quadratic and Cubic Bases}

In the following scalar transport results we will report results using $k=2$ and $k=3$ bases, which deserves more detailed explanation. If no SUPG stabilization terms are included, $k=3$ collocation schemes are well-defined, as the basis is globally $C^2$ and the governing PDE is second-order. Using $k=2$, however, gives a basis which is only $C^1$ globally, as there are jumps in the second derivative of the global basis function at the knot locations. In the case of even polynomial degrees, the selected collocation points will never fall on interior knot locations, and thus the collocation scheme remains well-defined. 

When the stabilization terms are included, the order of the governing PDE has increased to third-order. As explained above, the $k=2$ collocation points will always lie on element interiors. Thus we simply set the third derivative terms in this case equal to zero. For $k=3$, there will be collocation points located at the same locations in space as interior knot boundaries, where there exists a jump in the third derivative. To handle this case, we use the average of the third derivative values on each side of the knot locations as the estimate for the third derivative of the basis at the collocation point. 

\subsubsection{A Note on Enforcement of Dirichlet Boundary Conditions}

Throughout this work we will use equations in the form of Equation \eqref{eq:adv_diff_BC} to denote our enforcement of Dirichlet boundary conditions. This corresponds to a simple collocation of the provided data along the domain boundary. However, we will also investigate many problems which are subjected to piecewise constant boundary conditions. In these cases we will enforce these boundary conditions by instead setting the value of the corresponding control coefficients to the boundary data and removing each collocation point on the boundary. This results in a spline approximation of the boundary condition which is regularized, and avoids oscillations near the locations of discontinuities in the boundary data which would arise if the boundary condition was simply collocated. 

\subsubsection{Relation to Upwind Collocation}

We can also relate the SUPG collocation scheme to schemes in which collocation points are moved in the upwind direction, under some assumptions. Let us denote the PDE residual as $R^h = \textbf{u} \cdot \nabla \phi^h - \kappa \Delta \phi^h - f$ so that the discrete, stabilized PDE can be written as 

\begin{equation}
    (R^h - \nabla \cdot (\tau \mathbf{u} R^h))(\mathbf{x}_i) = 0.
\end{equation}

If we assume that the product $\tau |\mathbf{u}|$ is constant, then we can write this as

\begin{equation}
    \left(R^h - \tau |\mathbf{u}| \nabla \cdot \left( \frac{\mathbf{u}}{|\mathbf{u}|} R^h \right) \right)(\mathbf{x}_i) = 0.
\end{equation}

If $\frac{\mathbf{u}}{|\mathbf{u}|}$ is also a constant (meaning that both $\tau$ and $\mathbf{u}$ are constants when combined with the previous assumption), the second term is a directional derivative. Thus this expression can be interpreted as the first two terms of the Taylor expansion of $R^h (\mathbf{x}_i - \tau \mathbf{u})$. But evaluating $R^h (\mathbf{x}_i - \tau \mathbf{u})$ is exactly a collocation scheme with the collocation points moved in the upwind direction by a distance of $|\tau \mathbf{u}|$. Therefore the SUPG collocation scheme is equivalent to this upwinded collocation scheme up to the higher-order terms of the Taylor series. 

\subsection{Advection-Diffusion Results}

We now test the stabilized collocation method using standard benchmark problems in one and two spatial dimensions. In each dimension, we consider a boundary layer problem which demonstrates that the stabilization effectively suppresses oscillations in advection-dominated problems on coarse meshes. We also manufacture solutions to illustrate that the stabilization does not introduce excess error on well-resolved, smooth problems. 

\subsubsection{1D Boundary Layer Flow}

A classical test case for advection-stabilized methods is the 1D boundary layer flow. We solve the steady advection diffusion equation on the interval $[0, 1]$ with boundary conditions $g(0) = 0$ and $g(1) = 1$. It is well known that the exact solution to this problem is given by an exponential function which becomes sharper and more concentrated at the right boundary as the Péclet number increases. The large gradients near the right boundary commonly reveal the lack of advection stability in discretizations.

\begin{figure}
\centering
\subfloat[Without SUPG stabilization]{\label{sfig:BL_unstab}\includegraphics[width=.45\textwidth]{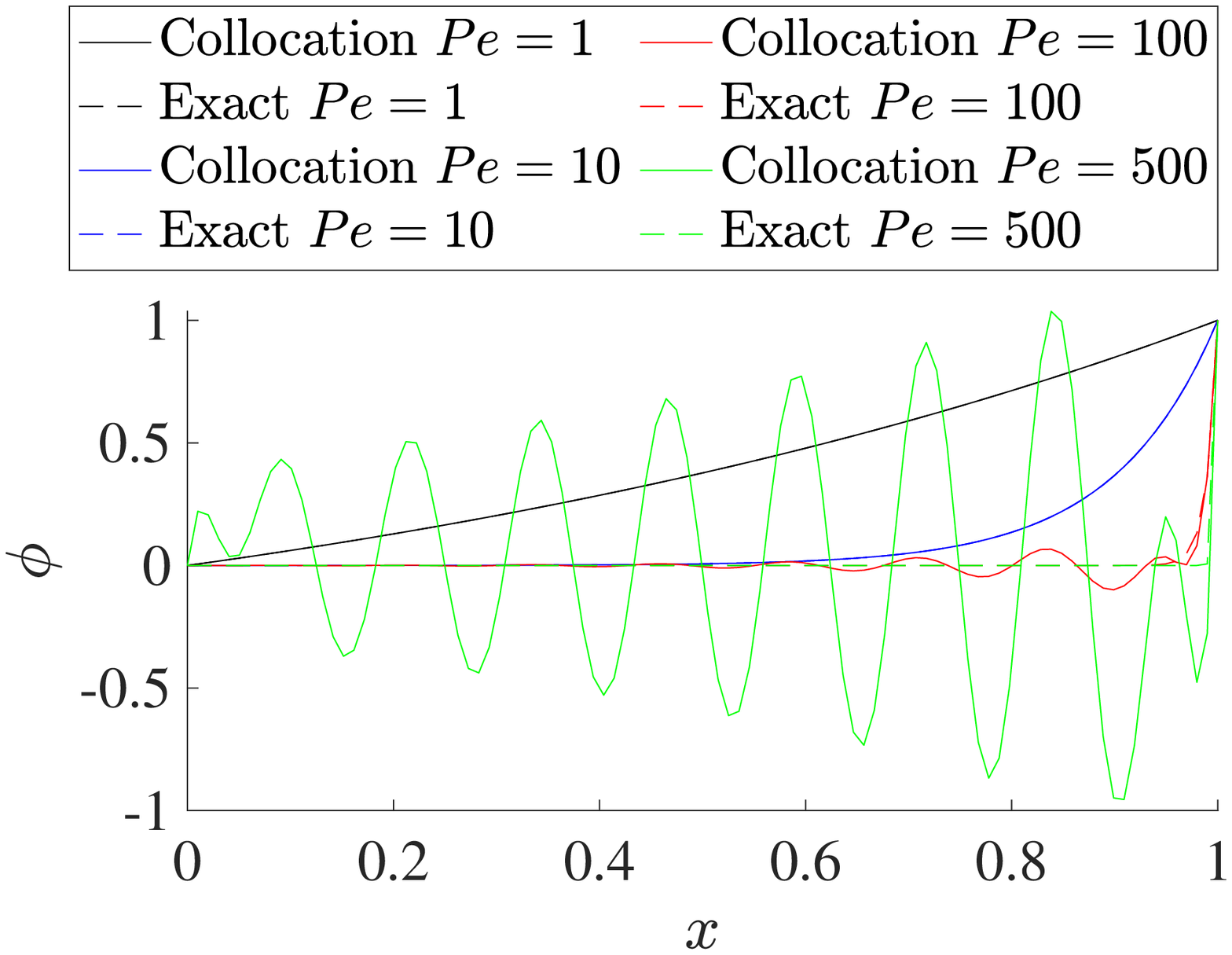}}\hfill
\subfloat[With SUPG stabilization]{\label{sfig:BL_stab}\includegraphics[width=.45\textwidth]{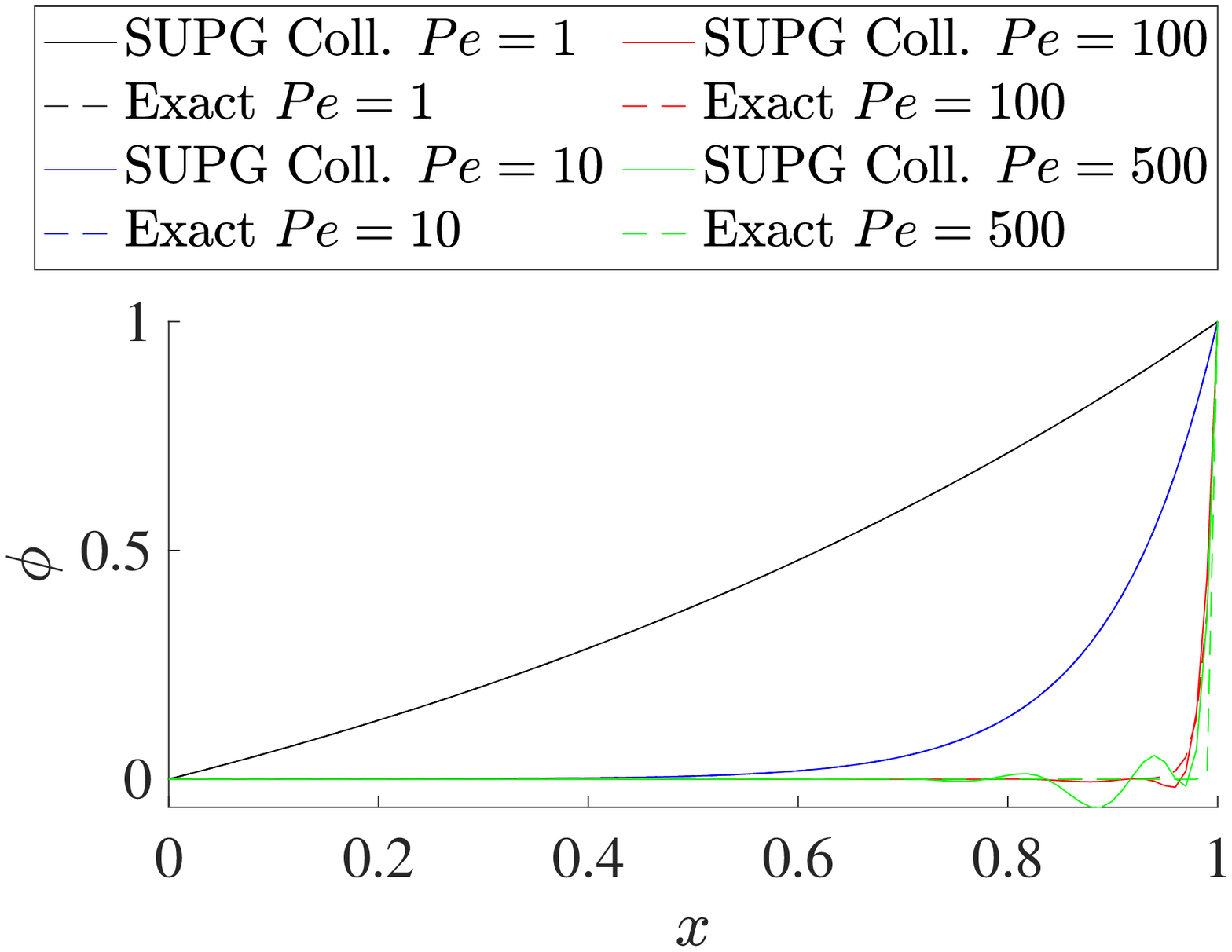}} \\
\caption{Collocation solution to boundary layer problem with and without SUPG stabilization, 16 elements and $k$ = 4}
\label{fig:BL}
\end{figure}

If we attempt to use a naïve collocation method with an underresolved mesh to solve this problem, we end up with results such as those shown in Figure \ref{sfig:BL_unstab}, which illustrates results obtained for a variety of Péclet numbers using a 16 element, $k=4$ mesh. As the problem becomes strongly advection-dominated the numerical solution becomes polluted by large oscillations, as expected. Figure \ref{sfig:BL_stab} shows the results generated by the stabilized collocation scheme on the same mesh. Clearly the inclusion of stabilization has greatly improved the results in the advection dominated regime, as the oscillations have been greatly reduced and localized near the boundary layer. Simultaneously, the approximations of results in lower Péclet number regimes remain reasonable.

\begin{figure}[t]
\centering
\includegraphics[width=.45\textwidth]{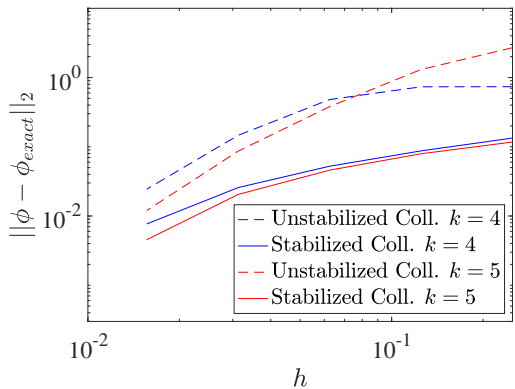}
\caption{$L^2$ errors in $Pe = 500$ solutions with and without SUPG stabilization}
\label{fig:BL_errors}
\end{figure}

For a more quantitative measure of the improved accuracy SUPG collocation provides in the advection dominated setting, we can consider the numerical solution error as a function of mesh resolution. Figure \ref{fig:BL_errors} shows the $L^2$ errors of both the unstabilized and stabilized numerical solution as a function of the Bézier mesh size $h$ for $k = 4$ and $k = 5$ bases with $Pe = 500$. Clearly in severely underresolved simulations the stabilized collocation scheme performs much better, returning almost an order of magnitude smaller error. As the meshes are refined, errors in both schemes seem to approach an asymptotic convergence rate, which we will explore next.

\subsubsection{1D Manufactured Solution}

\begin{figure}
\centering
\subfloat[Unstabilized $L^2$ error]{\label{sfig:unstab_l2}\includegraphics[width=.45\textwidth]{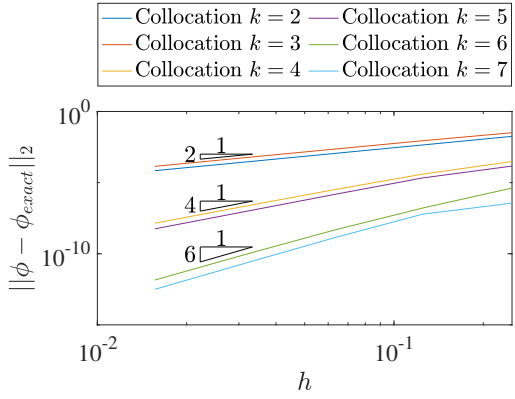}}\hfill
\subfloat[Stabilized $L^2$ error]{\label{sfig:stab_l2}\includegraphics[width=.45\textwidth]{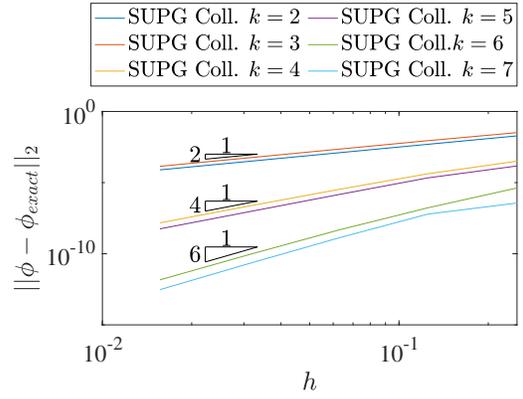}} \\
\subfloat[Unstabilized $H^1$ error]{\label{sfig:unstab_h1}\includegraphics[width=.45\textwidth]{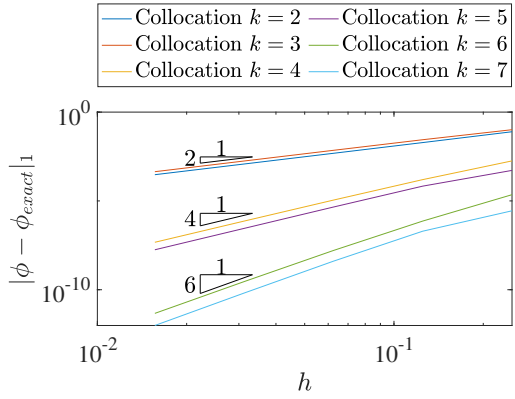}}\hfill
\subfloat[Stabilized $H^1$ error]{\label{sfig:stab_h1}\includegraphics[width=.45\textwidth]{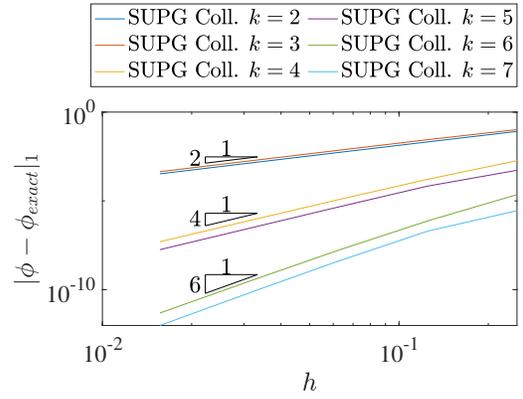}} \\
\caption{Errors in 1D advection-diffusion manufactured solution with and without SUPG stabilization}
\label{fig:ad_1d_conv}
\end{figure}

We use the method of manufactured solutions to determine the asymptotic order of accuracy of the collocation methods, beginning in 1D. The exact solution is chosen to be $\phi_{exact} = \sin(\pi x)$ over the unit interval from 0 to 1 with $Pe = 1$. Figure \ref{fig:ad_1d_conv} shows the convergence of both the $L^2$ norm and $H^1$ semi-norm of the error for both stabilized and unstabilized methods. In both cases we recover a rate of $k$ for even polynomial degrees and a rate of $k-1$ for odd polynomial degrees. Also note that the addition of the stabilizing terms has not made a noticeable impact on the magnitude of the errors in this well-resolved case. 

\subsubsection{Advection Skew to Mesh}

\begin{figure}[t]
\centering
\includegraphics[width=2.5in]{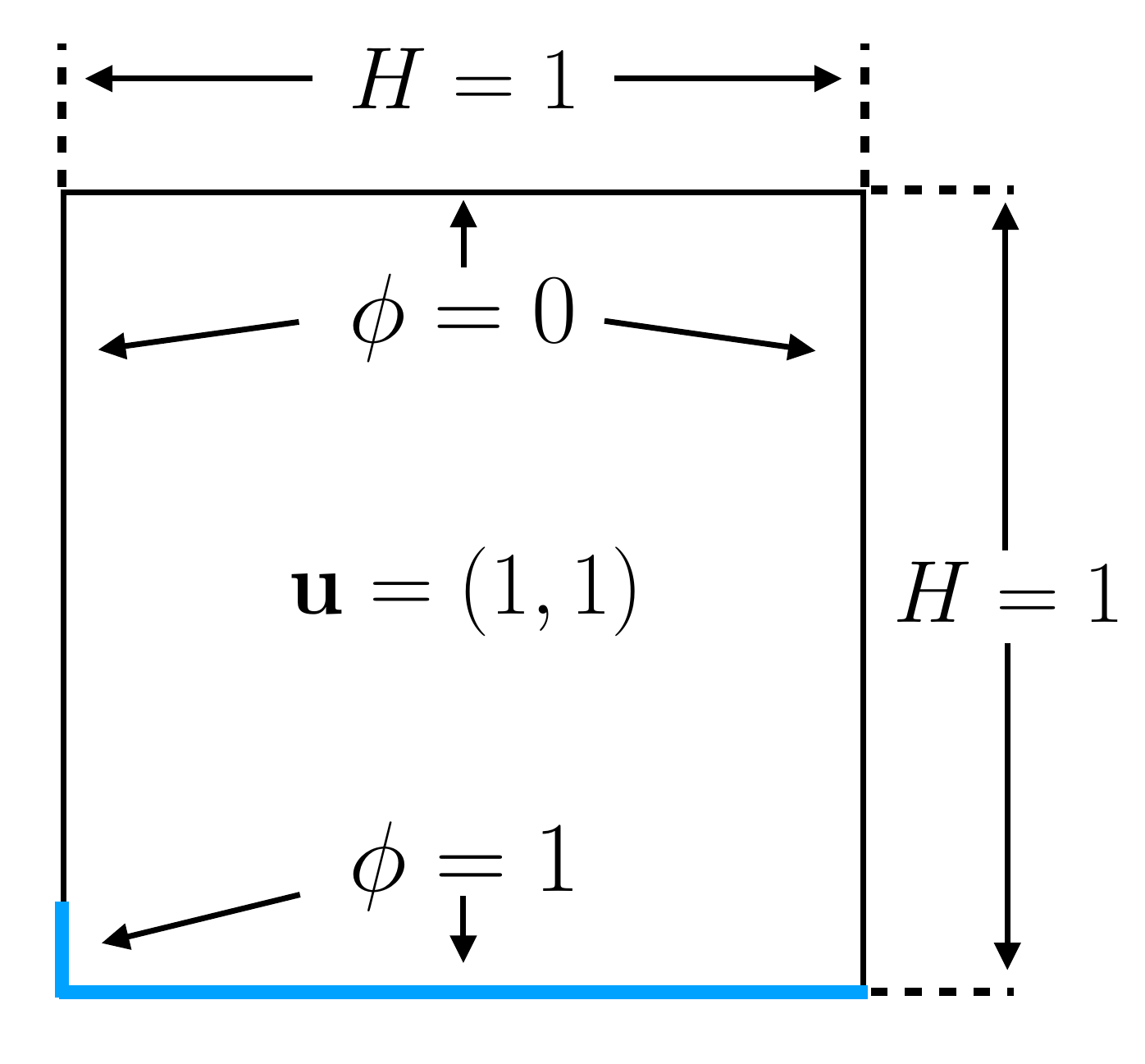}
\caption{Schematic of advection skew to mesh in 2D}
\label{fig:skew_setup}
\end{figure}

\begin{figure}
\centering
\subfloat[Without SUPG stabilization]{\label{sfig:skew_unstab}\includegraphics[width=.45\textwidth]{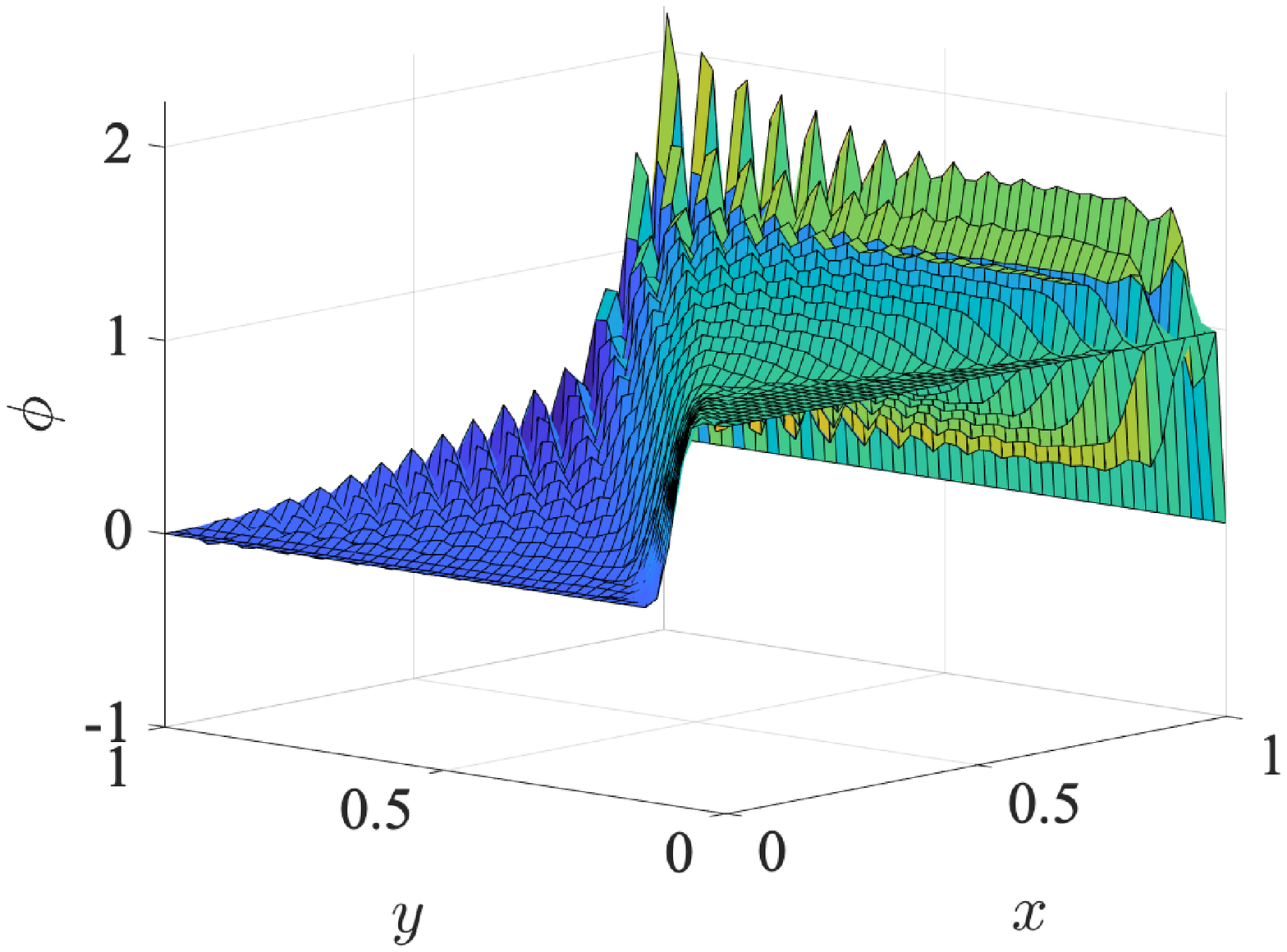}}\hfill
\subfloat[With SUPG stabilization]{\label{sfig:skew_stab}\includegraphics[width=.45\textwidth]{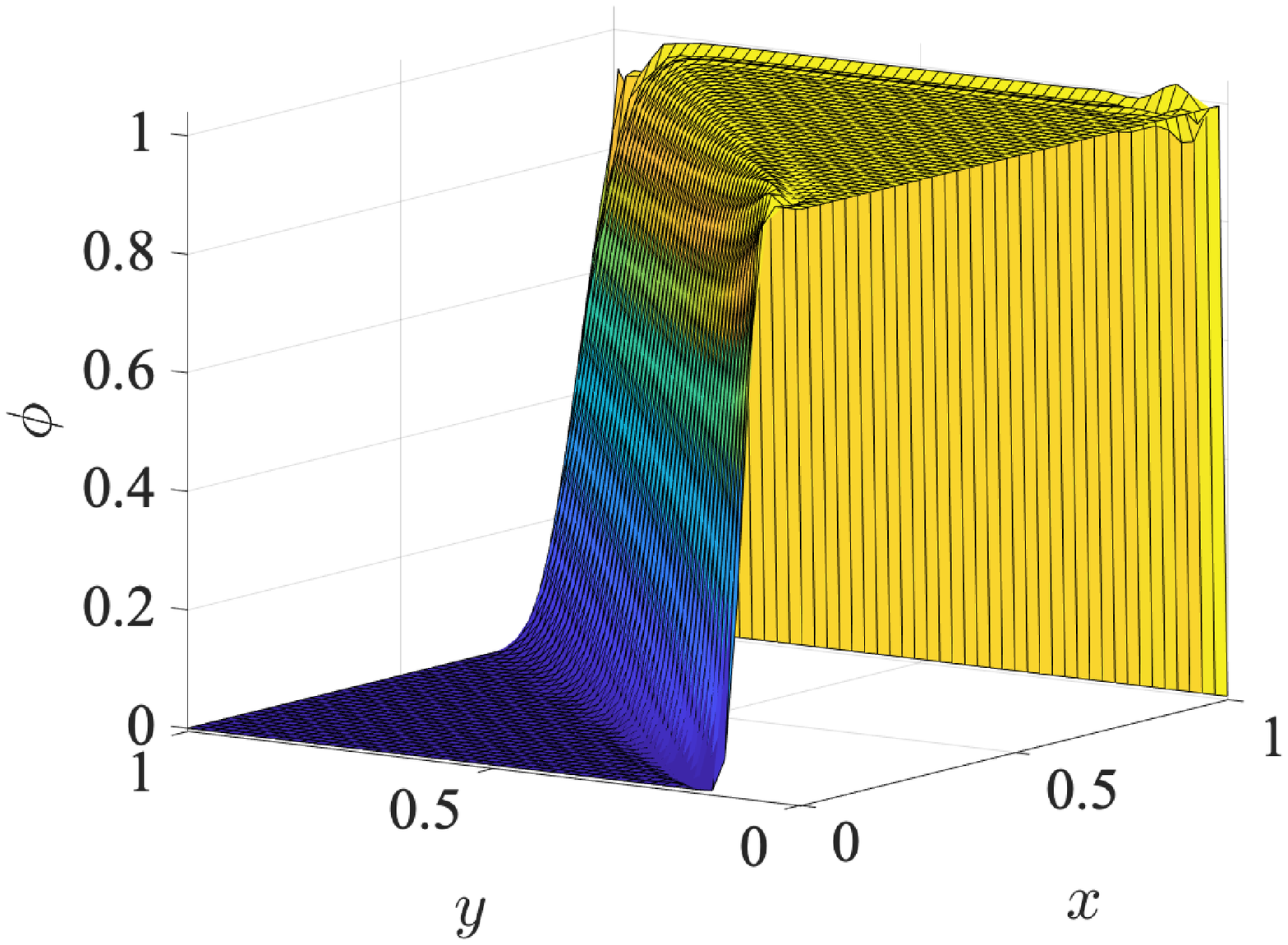}} \\
\caption{Collocation solution to skew advection problem with and without SUPG stabilization, $Pe=1000$, 32$^2$ elements and $k$ = 4}
\label{fig:skew}
\end{figure}

Another classical test case for stabilized methods is advection skew to the mesh. This is in some ways a generalization of the 1D boundary layer problem into 2D, and the setup is shown in Figure \ref{fig:skew_setup}. In particular, we consider the unit square domain ($H = 1$) and prescribe a Dirichlet boundary condition of 1 along the bottom of the square and along the first tenth of the left vertical side. Everywhere else along the boundary the solution is set to zero. We also set $Pe = 1000$, which will result in very sharp boundary layers forming in the transition regions of the solution. Figure \ref{sfig:skew_unstab} shows the solution when collocation is applied without any stabilization on a mesh of $32^2$ elements with $k = 4$. Once again the solution is corrupted by spurious oscillations. Figure \ref{sfig:skew_stab} shows that when SUPG stabilization is included the solution is much better quality.

\subsubsection{2D Manufactured Solution}

\begin{figure}
\centering
\subfloat[Unstabilized $L^2$ error]{\label{sfig:unstab_l2_2D}\includegraphics[width=.45\textwidth]{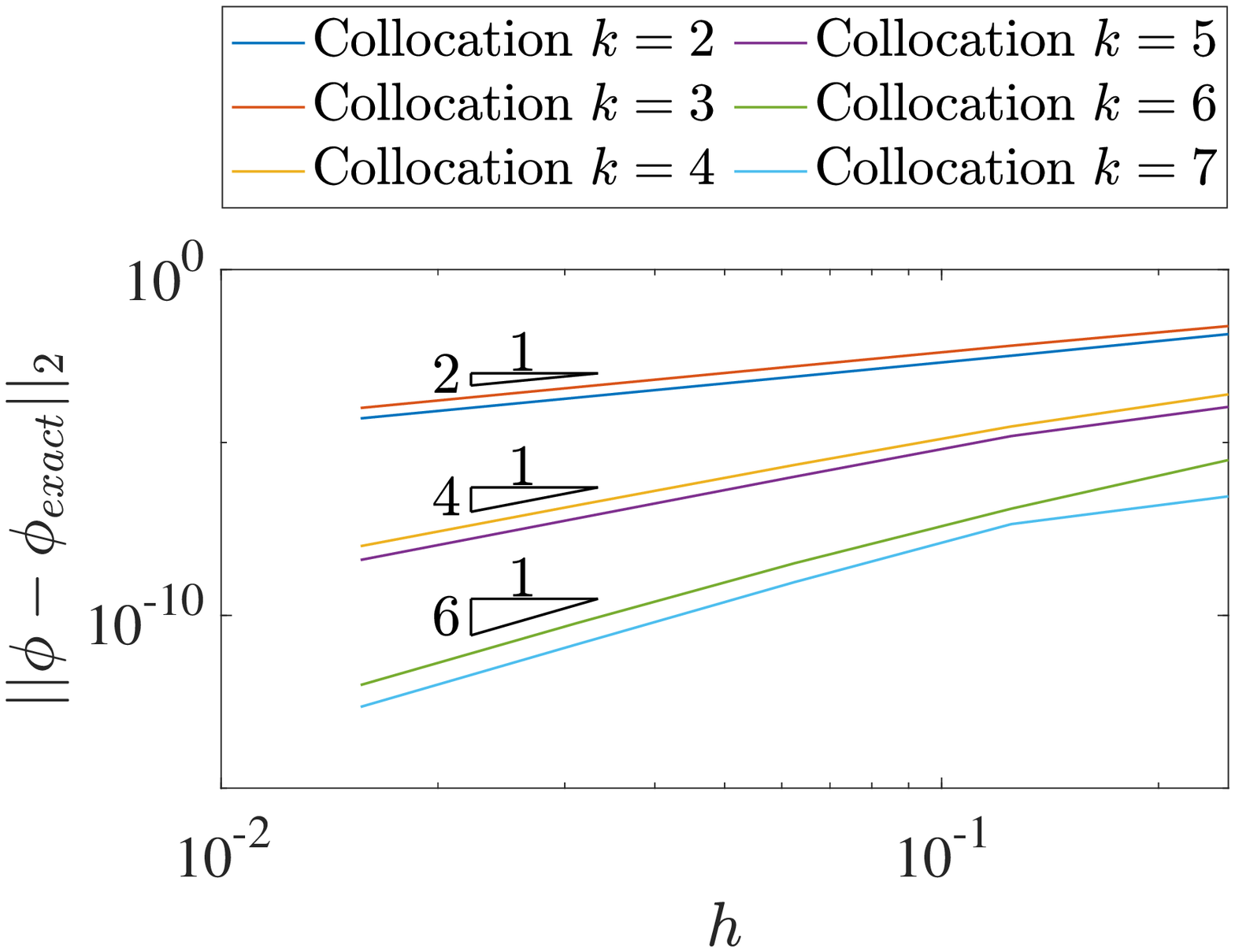}}\hfill
\subfloat[Stabilized $L^2$ error]{\label{sfig:stab_l2_2D}\includegraphics[width=.45\textwidth]{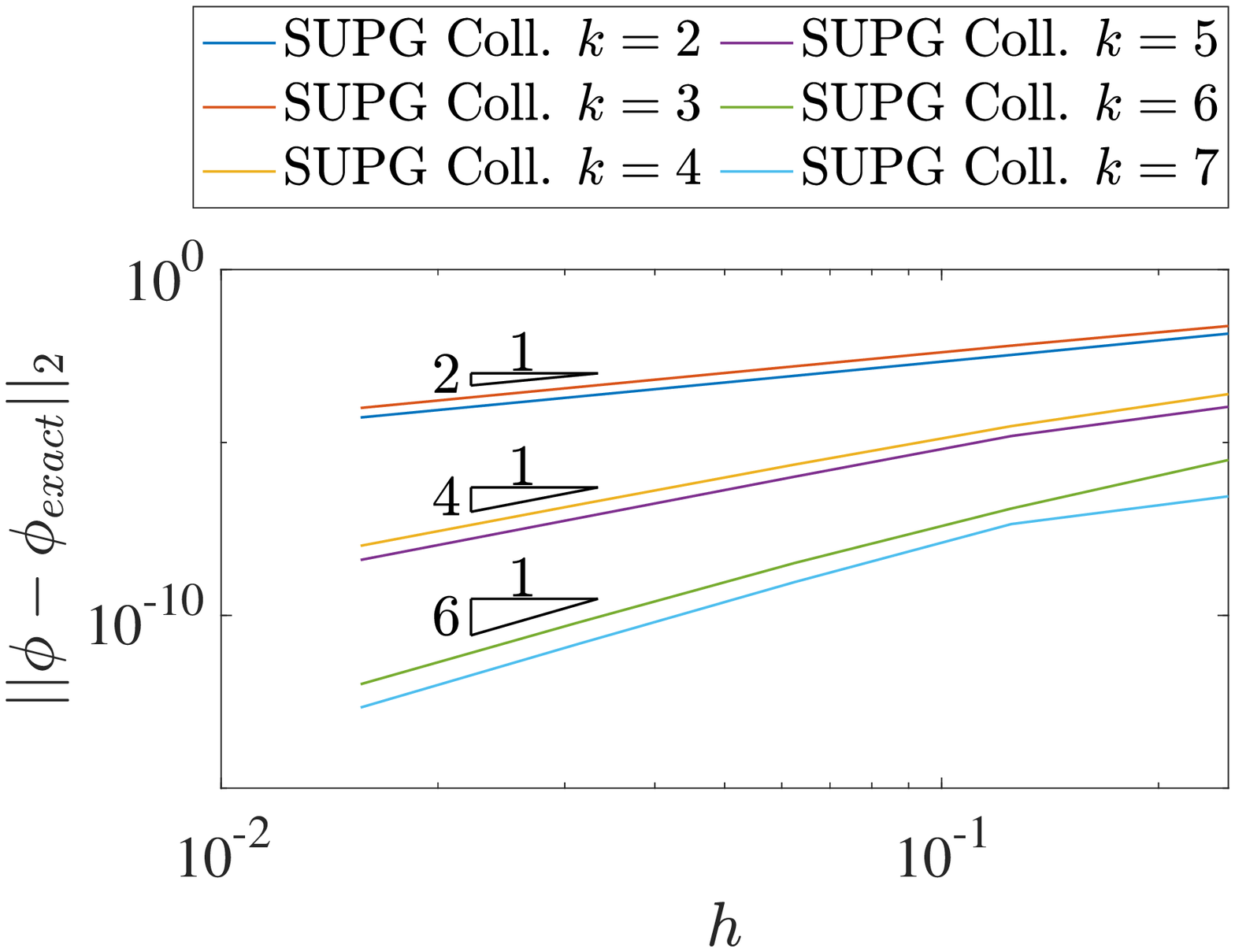}} \\
\subfloat[Unstabilized $H^1$ error]{\label{sfig:unstab_h1_2D}\includegraphics[width=.45\textwidth]{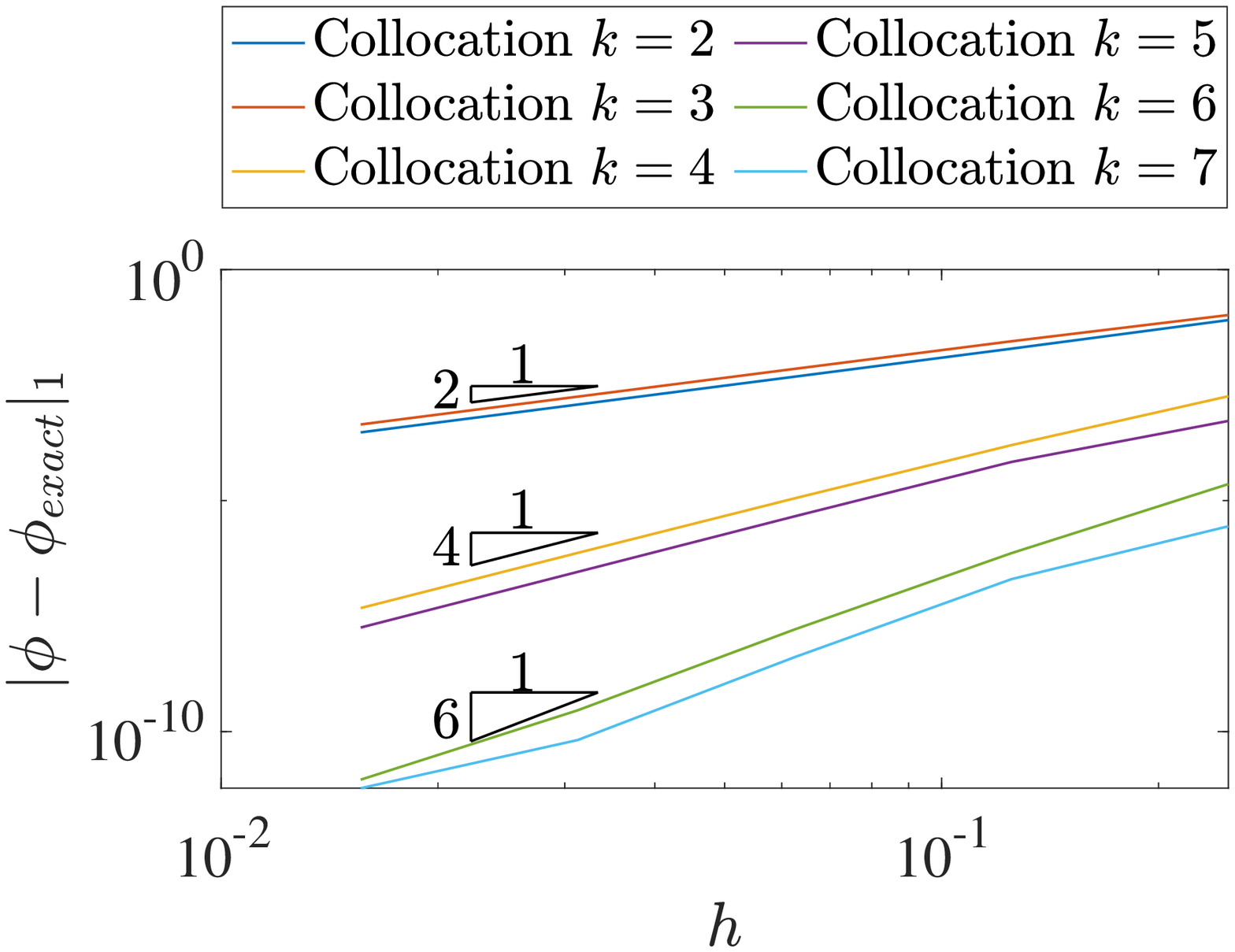}}\hfill
\subfloat[Stabilized $H^1$ error]{\label{sfig:stab_h1_2D}\includegraphics[width=.45\textwidth]{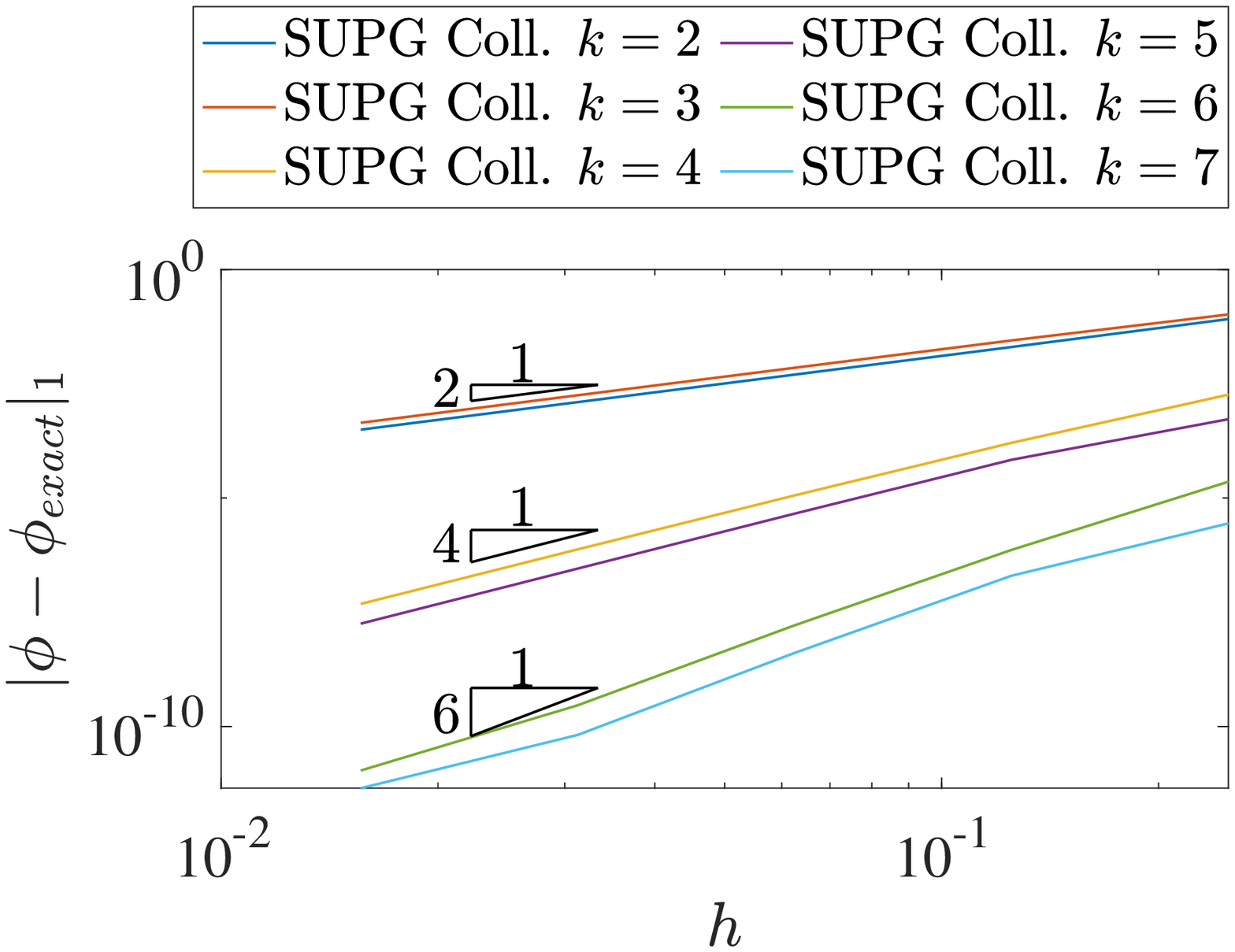}} \\
\caption{Errors in 2D advection-diffusion manufactured solution with and without SUPG stabilization}
\label{fig:ad_2d_conv}
\end{figure}

Like in the 1D setting we also investigate the convergence rates of the scheme on a smooth problem and ensure that the error in the stabilized scheme remains comparable to the unstabilized case. In this case we use a manufactured solution with exact solution $\phi_{exact} = \sin{(\pi x)}\sin{(\pi y)}$ over the unit square with $Pe = 1$. Figure \ref{fig:ad_2d_conv} shows the convergence of the errors. Once again we recover rates of $k$ for even degree bases and $k-1$ for odd degree bases, and the magnitude of the errors is similar between stabilized and unstabilized schemes for this well-resolved case.

\section{Incompressible Stokes Flow}

At this point we transition from advection stabilization and consider instead pressure stabilization, using the incompressible Stokes equations as a model problem. It is well known that well-posedness of the discrete system requires careful selection of the velocity and pressure approximation spaces, and approaches often represent the discrete velocity with a higher polynomial degree than the discrete pressure. Equal-order spaces, though not naturally stable, can be more convenient in implementation, which thus led to the development of pressure stabilization methods for Galerkin discretizations. In this section we will show that one of these methods, namely the PSPG method, can also be effectively used in the collocation setting to stabilize equal-order discretizations. 

\subsection{The Strong Form of the Stokes Equations}

The standard, velocity-pressure form of the steady Stokes equations subjected to Dirichlet boundary conditions is given by:

\bigskip

$$
\left\{ \hspace{5pt}
\parbox{5in}{
\noindent Given $\mu \in \mathbb{R}^+$, $\textbf{f} : \Omega \rightarrow \mathbb{R}^n$, and $\textbf{g} : \partial \Omega \rightarrow \mathbb{R}^n$, find $\mathbf{u} : \Omega \rightarrow \mathbb{R}^n$ and $p : \Omega \rightarrow \mathbb{R}$ such that:
\begin{align}
    -\mu \Delta \mathbf{u} + \nabla p = \mathbf{f} \quad &\textup{in} \quad \Omega \label{eq:stokes_mom}, \\
    \nabla \cdot \mathbf{u} = 0 \quad &\textup{in} \quad \Omega, \label{eq:stokes_mass}  \\
    \textbf{u} = \textbf{g} \quad &\textup{on} \quad \partial \Omega,
\end{align}
}
\right.
$$

\bigskip

\noindent where $\mathbf{u}$ is the unknown velocity field, $p$ is the unknown pressure field, $\mu$ is the viscosity, $\mathbf{f}$ is a prescribed forcing function, and $\mathbf{g}$ is the prescribed Dirichlet boundary data. Equation \eqref{eq:stokes_mom} is a statement of the conservation of momentum while Equation \eqref{eq:stokes_mass} states conservation of mass for an incompressible fluid. 

In \cite{aronson2023divergence}, it was also shown that divergence-conforming collocation schemes based on the rotational form, or vorticity-velocity-pressure form, of the Stokes and Navier-Stokes equations achieved accelerated convergence rates when compared with the standard, velocity-pressure scheme. We wish to see if this also holds true with stabilized methods and thus we also review the rotational form of the (2D) Stokes equations

\bigskip

$$
\left\{ \hspace{5pt}
\parbox{5in}{
\noindent Given $\mu \in \mathbb{R}^+$, $\textbf{f} : \Omega \rightarrow \mathbb{R}^2$, and $\textbf{g} : \partial \Omega \rightarrow \mathbb{R}^2$, find $\mathbf{u} : \Omega \rightarrow \mathbb{R}^2$, $p : \Omega \rightarrow \mathbb{R}$, and $\omega : \Omega \rightarrow \mathbb{R}$ such that:
\begin{align}
    \mu \nabla^{\perp} \omega + \nabla p = \mathbf{f} \quad &\textup{in} \quad \Omega,\\
    \nabla \cdot \mathbf{u} = 0 \quad &\textup{in} \quad \Omega,\\
    \omega - \nabla \times \mathbf{u} = 0 \quad &\textup{in} \quad \Omega ,\label{eq:stokes_rot_const}\\
    \textbf{u} = \textbf{g} \quad &\textup{on} \quad \partial \Omega.
\end{align}
}
\right.
$$

\noindent Here we have introduced the scalar vorticity $\omega$ as an additional unknown and Equation \eqref{eq:stokes_rot_const} as the constitutive law relating velocity and vorticity. The rotor, or perpendicular gradient, operator $\nabla^\perp$ is defined for scalar functions $\omega$ via $\nabla^\perp \omega = \left( \frac{\partial \omega}{\partial y}, -\frac{\partial \omega}{\partial x} \right)$. 

\bigskip

\subsection{PSPG Stabilization for Stokes Flow}

In the Galerkin setting, one method to stabilize equal-order interpolations of velocity and pressure is to include PSPG stabilization in the weak form. We define the discrete spline space $\mathbf{S}^h$ as the vector-valued spline space such that each component is a member of $S^h$. We similarly define $\mathbf{S}^h_g$ as the members of $\mathbf{S}^h$ satisfying the Dirichlet boundary conditions and $\mathbf{S}^h_0$ as the functions in $\mathbf{S}^h$ equal to zero on the Dirichlet boundaries. Thus the discrete velocity is a member of $\mathbf{S}^h_g$ and the velocity test functions lie in $\mathbf{S}^h_0$. For a pure Dirichlet problem, we must also constrain the discrete pressure and corresponding test functions to the space $\Tilde{S}^h$, where $\Tilde{S}^h$ is the space of all functions in $S^h$ with zero mean. Once again we assume that the basis functions are $C^1$ or smoother. Then the stabilized weak form of the velocity-pressure form of the Stokes equations is given by

\bigskip

$$
\left\{ \hspace{5pt}
\parbox{5in}{
Given $\mu \in \mathbb{R}^+$, $\textbf{f} : \Omega \rightarrow \mathbb{R}^n$, and $\textbf{g} : \partial \Omega \rightarrow \mathbb{R}^n$, find $\mathbf{u}^h \in \mathbf{S}^h_g$ and $p^h \in \Tilde{S}^h$ such that :

\begin{align}
    \int_{\Omega} (\mu \nabla \mathbf{w}^h \cdot \nabla \mathbf{u}^h - p^h \nabla \cdot \mathbf{w}^h) d\Omega = \int_{\Omega} \mathbf{w}^h \cdot \mathbf{f} d \Omega \quad \forall \mathbf{w}^h \in \mathbf{S}^h_0, \label{eq:stokes_mom_weak}\\
    \int_{\Omega} q^h (\nabla \cdot \mathbf{u}^h) d \Omega + \int_{\Omega} \tau (\nabla q^h \cdot \mathbf{R}^h) d \Omega = 0 \quad \forall q^h \in \Tilde{S}^h ,\label{eq:stokes_mass_weak} 
\end{align}
}
\right.
$$

\bigskip

\noindent where $\tau$ is a stabilization constant and we have defined the discrete momentum PDE residual $\mathbf{R}^h$ as:
\begin{equation}
    \mathbf{R}^h = -\mu \Delta \mathbf{u}^h + \nabla p^h - \mathbf{f}.
\end{equation}

Much like in the scalar transport setting, we develop inspiration for our stabilized collocation schemes by reversing the integration by parts within the weak forms of both the mass and momentum equations used in the Galerkin formulation. This operation is simple for the momentum equations, which become

\begin{equation}
    \int_\Omega \mathbf{w}^h \cdot (-\mu \Delta \mathbf{u}^h + \nabla p^h - \mathbf{f}) d\Omega = 0.
\end{equation}

\noindent When the continuity equation is integrated by parts we must note that the pressure test function is not required to satisfy any homogeneous Dirichlet boundary conditions. Thus the resulting boundary term does not disappear as in the momentum equations, and we end up with

\begin{equation}
    \int_\Omega q^h(\nabla \cdot \mathbf{u}^h - \nabla \cdot (\tau \mathbf{R}^h)) d\Omega + \int_{\partial \Omega} q^h \tau \mathbf{R}^h \cdot \mathbf{n} d \Gamma = 0.
\end{equation}

To turn this statement into a collocation scheme, we first need to manipulate this boundary term. This term resembles a Neumann boundary condition on the pressure field, and so we adopt an approach similar to the Enhanced Collocation approach developed in \cite{deLorenzis2015isogeometric}. In particular, after assuming a Dirac-delta test function $q^h$ centered on a boundary collocation point $\mathbf{x}_i$, the continuity equation on the boundary will be of the form

\begin{equation}
    (\nabla \cdot \mathbf{u}^h - \nabla \cdot (\tau \mathbf{R}^h) + \frac{C}{h_b}(\tau \mathbf{R}^h \cdot \mathbf{n}))(\mathbf{x}_i) = 0,
\end{equation}

\noindent where $C$ is a sufficiently large constant and $h_b$ is the mesh size in the boundary normal direction. For our tests, $C$ is set to be 1, though we will also include results to demonstrate the sensitivity of the results of a sample problem to this choice.

To fully define the collocation scheme, we again define the test functions $\mathbf{w}^h$ and $q^h$ to be Dirac-delta functions located at the collocation points associated with the spline space $S^h$. Note that the same set of collocation points are used for the momentum and mass balance equations, namely the Greville abscissae $\{\mathbf{x}\}_{i=1}^n$ of the spline space $S^h$. Then the stabilized collocation scheme is defined by

\bigskip

$$
\left\{ \hspace{5pt}
\parbox{5in}{
\noindent Find $\mathbf{u}^h \in \mathbf{S}^h$ and $p^h \in S^h$ such that, for $i = 1, ..., n$

\begin{equation}
    (-\mu \Delta \mathbf{u}^h + \nabla p^h - \mathbf{f})(\mathbf{x}_i) = \mathbf{0} \quad \forall \mathbf{x}_i \in \Omega,
\end{equation}
\begin{equation}
    (\nabla \cdot \mathbf{u}^h - \nabla \cdot (\tau \mathbf{R}^h))(\mathbf{x}_i) = 0 \quad \forall \mathbf{x}_i \in \Omega,
\end{equation}
\begin{equation}
    (\nabla \cdot \mathbf{u}^h - \nabla \cdot (\tau \mathbf{R}^h) + \frac{C}{h_b}(\tau \mathbf{R}^h \cdot \mathbf{n}))(\mathbf{x}_i) = 0 \quad \forall \mathbf{x}_i \in \partial \Omega,
\end{equation}
\begin{equation}
    \mathbf{u}^h(\textbf{x}_i) = \mathbf{g}(\mathbf{x}_i) \quad \forall \textbf{x}_i \in \partial \Omega.
\end{equation}

}
\right.
$$

\bigskip

\noindent Note that we have written two expressions for the continuity equation, with one being valid at collocation points within the domain and another being valid on the boundary collocation points. Finally, the stabilization constant is defined as 

\begin{equation}
    \tau = \frac{h^2}{C_2 \mu},
\end{equation}

\noindent where $C_2$ is set to be 4, $h$ is the average distance between a collocation point and its neighbors, and we fit a spline from the solution space $S^h$ to the $\tau$ values to define the necessary derivatives. 

This concludes the construction of the PSPG stabilized collocation scheme based on the velocity-pressure form of the Stokes equations. We follow the same steps to derive a collocation scheme based on the rotational form of the equations. In particular, we include the PSPG stabilization term in the weak form of these equations in the same manner as above, with the only change being how the momentum residual is defined to be consistent with the new form of the equations

\begin{equation}
    \mathbf{R}^h = \nabla^{\perp} \omega^h + \nabla p^h - \mathbf{f}.
\end{equation}

By reverse integrating by parts as before, we can arrive at a new set of collocation equations for the rotational form of the Stokes equations given by

\bigskip

$$
\left\{ \hspace{5pt}
\parbox{5in}{
\noindent Find $\mathbf{u}^h \in \mathbf{S}^h$, $p^h \in \Tilde{S}^h$, and $\omega^h \in S^h$ such that, for $i = 1, ..., n$
\begin{equation}
    (\nabla^{\perp} \omega^h + \nabla p^h - \mathbf{f})(\mathbf{x}_i) = \mathbf{0} \quad \forall \mathbf{x}_i \in \Omega,
\end{equation}
\begin{equation}
    (\nabla \cdot \mathbf{u}^h - \nabla \cdot (\tau \mathbf{R}^h))(\mathbf{x}_i) = 0 \quad \forall \mathbf{x}_i \in \Omega
\end{equation}
\begin{equation}
    (\omega^h - \nabla \times \mathbf{u}^h)(\mathbf{x}_i) = 0 \quad \forall \mathbf{x}_i \in \Omega,
\end{equation}
\begin{equation}
    (\nabla \cdot \mathbf{u}^h - \nabla \cdot (\tau \mathbf{R}^h) + \frac{C}{h_b}(\tau \mathbf{R}^h \cdot \mathbf{n}))(\mathbf{x}_i) = 0 \quad \forall \mathbf{x}_i \in \partial \Omega,
\end{equation}
\begin{equation}
    \mathbf{u}^h(\textbf{x}_i) = \mathbf{g}(\mathbf{x}_i) \quad \forall \textbf{x}_i \in \partial \Omega.
\end{equation}
}
\right.
$$

\bigskip

\noindent Note that we use the same spline space to discretize the vorticity as velocity and pressure, and thus use the same collocation points for evaluation of the constitutive law.

\subsection{Stokes Equations Results}

With the stabilized collocation schemes fully defined, we now test them on some classic benchmark problems for Stokes flow. We start with a manufactured solution to evaluate and compare the numerical errors and asymptotic convergence rates of the two schemes. We then consider a lid-driven cavity to further demonstrate that the stabilized schemes effectively remove any spurious pressure oscillations. 

\subsubsection{Manufactured Vortex Solution}

We test the accuracy of the stabilized collocation schemes for the Stokes equations with a manufactured solution, defined on the unit square domain by the following velocity and pressure exact solutions

\begin{equation}
    \bar{\textbf{u}} = \left[
\begin{array}{c}
2e^x(-1+x)^2x^2(y^2-y)(-1+2y) \\
(-e^x(-1+x)x(-2+x(3+x))(-1+y)^2y^2)
\end{array}
\right],
\end{equation}

\begin{equation}
    \left.
\begin{array}{ccc}
\bar{p} & = & (-424+156e+(y^2-y)(-456+e^x(456+x^2(228-5(y^2-y))+ \\ &&
2x(-228+(y^2-y))+2x^3(-36+(y^2-y))+x^4(12+(y^2-y))))).
\end{array}
\right.
\end{equation}

\noindent We select $\mu = 1$ and start by considering the velocity-pressure formulation, meaning that the forcing term is defined as

\begin{equation}
    \textbf{f} = -\mu \Delta \bar{\textbf{u}} + \nabla \bar{p}.
\end{equation}

\begin{figure}
\centering
\subfloat[Velocity $L^2$ error]{\label{sfig:stokes_vel_l2}\includegraphics[width=.45\textwidth]{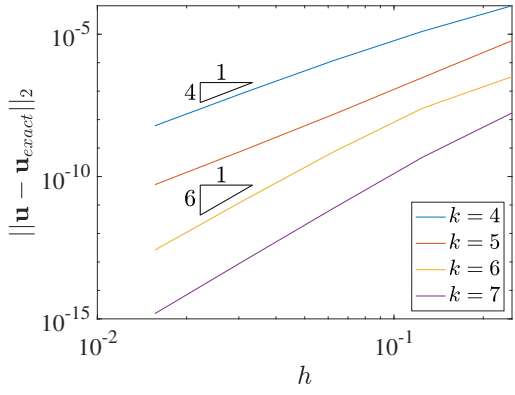}}\hfill
\subfloat[Pressure $L^2$ error]{\label{sfig:stokes_press_l2}\includegraphics[width=.45\textwidth]{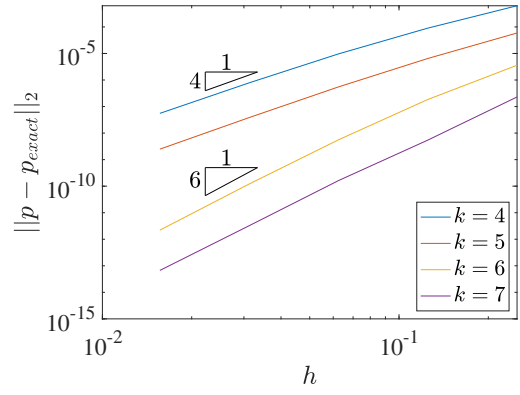}} \\
\subfloat[Velocity $H^1$ error]{\label{sfig:stokes_vel_h1}\includegraphics[width=.45\textwidth]{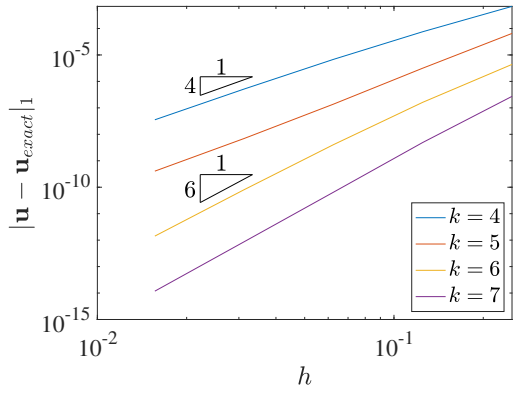}}\hfill
\subfloat[Pressure $H^1$ error]{\label{sfig:stokes_press_h1}\includegraphics[width=.45\textwidth]{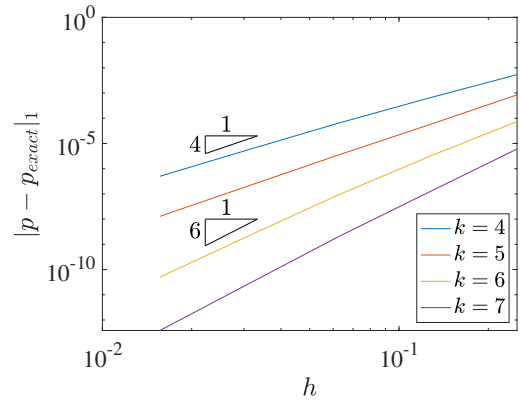}} \\
\caption{Errors in 2D Stokes manufactured solution using velocity-pressure form}
\label{fig:stokes_2d_conv}
\end{figure}

Figure \ref{fig:stokes_2d_conv} shows the $L^2$ norm and $H^1$ semi-norm errors for both velocity and pressure as functions of the polynomial degree and mesh size. For even polynomial degrees $k$ we recover convergence rates of $k$ in both norms for both velocity and pressure, and for odd polynomial degrees the solutions converge at a rate of $k-1$. This is consistent with the rates seen in previous studies of isogeometric collocation. 

\begin{figure}
\centering
\subfloat[Velocity $L^2$ error]{\label{sfig:stokes_rot_vel_l2}\includegraphics[width=.45\textwidth]{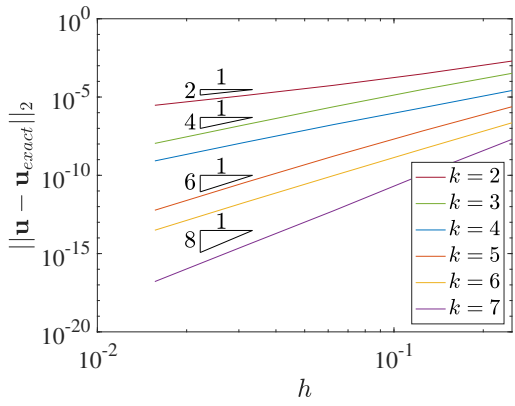}}\hfill
\subfloat[Pressure $L^2$ error]{\label{sfig:stokes_rot_press_l2}\includegraphics[width=.45\textwidth]{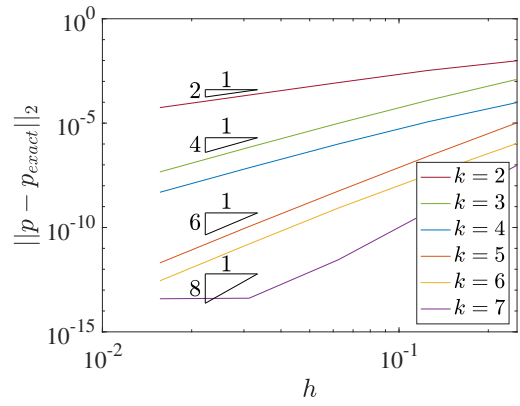}} \\
\subfloat[Velocity $H^1$ error]{\label{sfig:stokes_rot_vel_h1}\includegraphics[width=.45\textwidth]{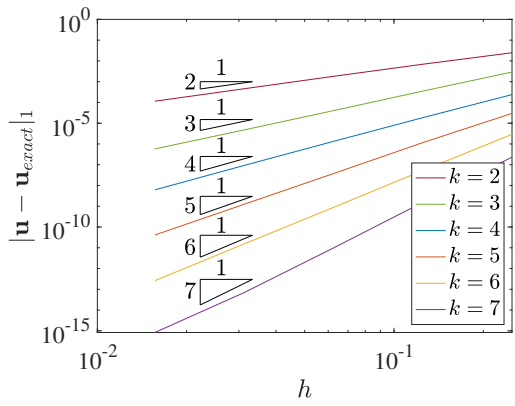}}\hfill
\subfloat[Pressure $H^1$ error]{\label{sfig:stokes_rot_press_h1}\includegraphics[width=.45\textwidth]{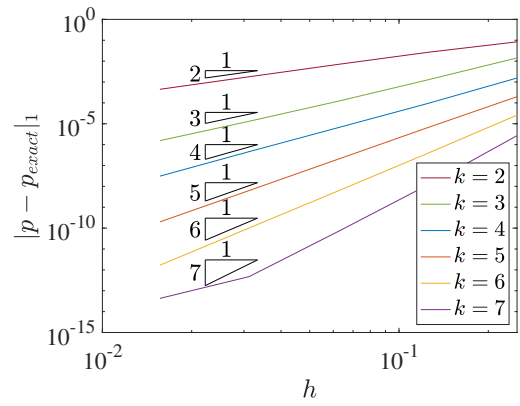}} \\
\caption{Errors in 2D Stokes manufactured solution using rotational form}
\label{fig:stokes_rot_2d_conv}
\end{figure}

Using the same manufactured solution, Figure \ref{fig:stokes_rot_2d_conv} shows the convergence of $L^2$ norm and $H^1$ semi-norm of the error for the scheme based on the rotational formulation of the Stokes equations. In this case the exact vorticity is given by

\begin{equation}
    \bar{\omega} = \frac{\partial \bar{u}_y}{\partial x} - \frac{\partial \bar{u}_x}{\partial y},
\end{equation}

\noindent and the forcing becomes 

\begin{equation}
    \mathbf{f} = \mu \nabla^{\perp} \bar{\omega} + \nabla \bar{p}.
\end{equation}

\noindent Note that because of the reduced order of the governing equations the scheme is well-defined for polynomial degrees of 2 and 3, as the stabilized equations are second-order. When $k = 2$ the basis is globally $C^1$, but the collocation points will not lie on knot boundaries where the second derivatives are discontinuous. For $k = 3$ the basis is $C^2$, meaning this stabilized scheme is well-defined everywhere. Now the errors for odd polynomial degrees converge like $k+1$ in $L^2$ and $k$ in $H^1$, which match the optimal rates seen in Galerkin methods much like the three-field schemes in \cite{aronson2023divergence}. Meanwhile, the convergence rates for even polynomial degrees are the same as those obtained with the two-field formulation.

\begin{figure}
\centering
\subfloat[Velocity $L^2$ error]{\label{sfig:stokes_vel_l2_edge}\includegraphics[width=.45\textwidth]{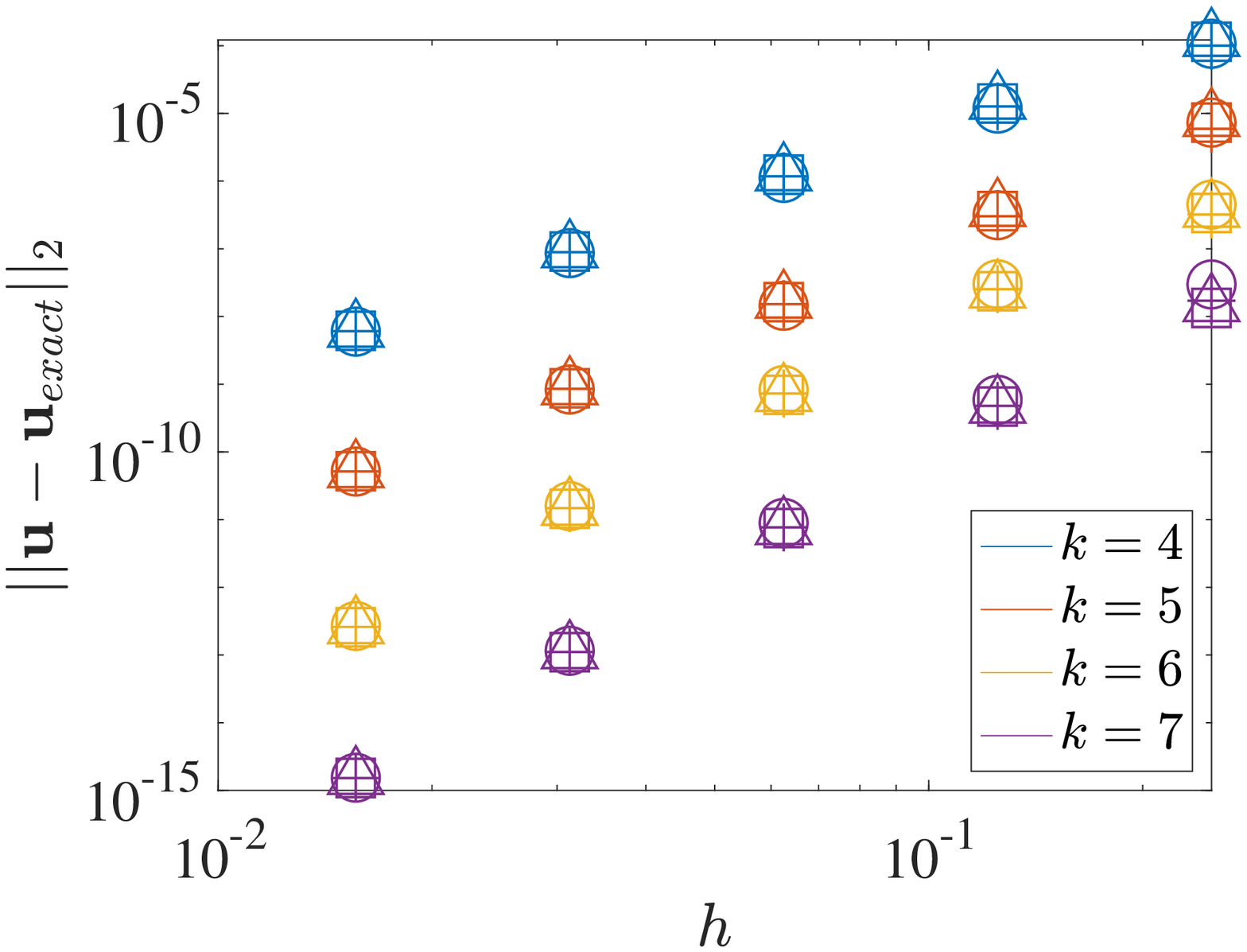}}\hfill
\subfloat[Pressure $L^2$ error]{\label{sfig:stokes_press_l2_edge}\includegraphics[width=.45\textwidth]{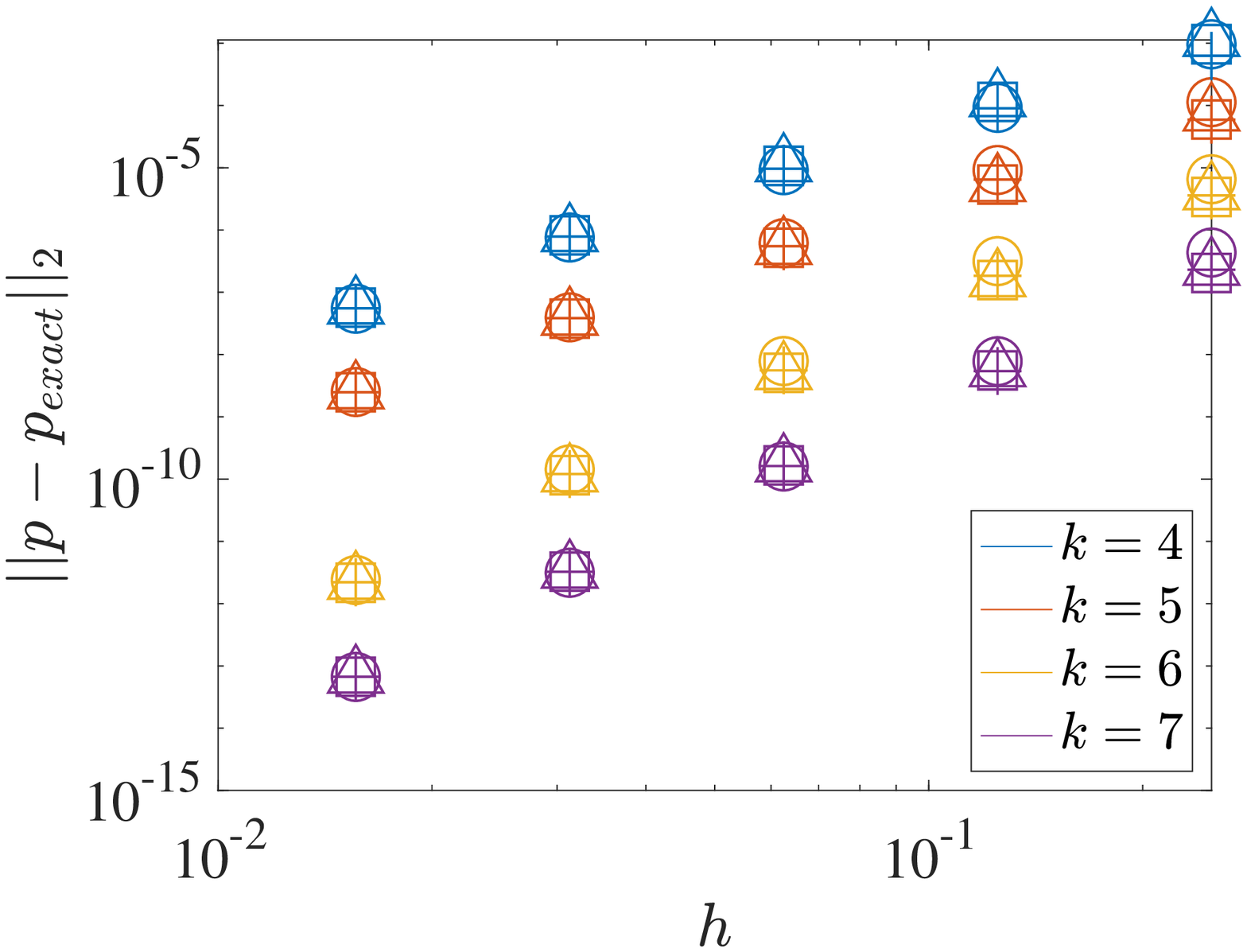}} \\
\subfloat[Velocity $H^1$ error]{\label{sfig:stokes_vel_h1_edge}\includegraphics[width=.45\textwidth]{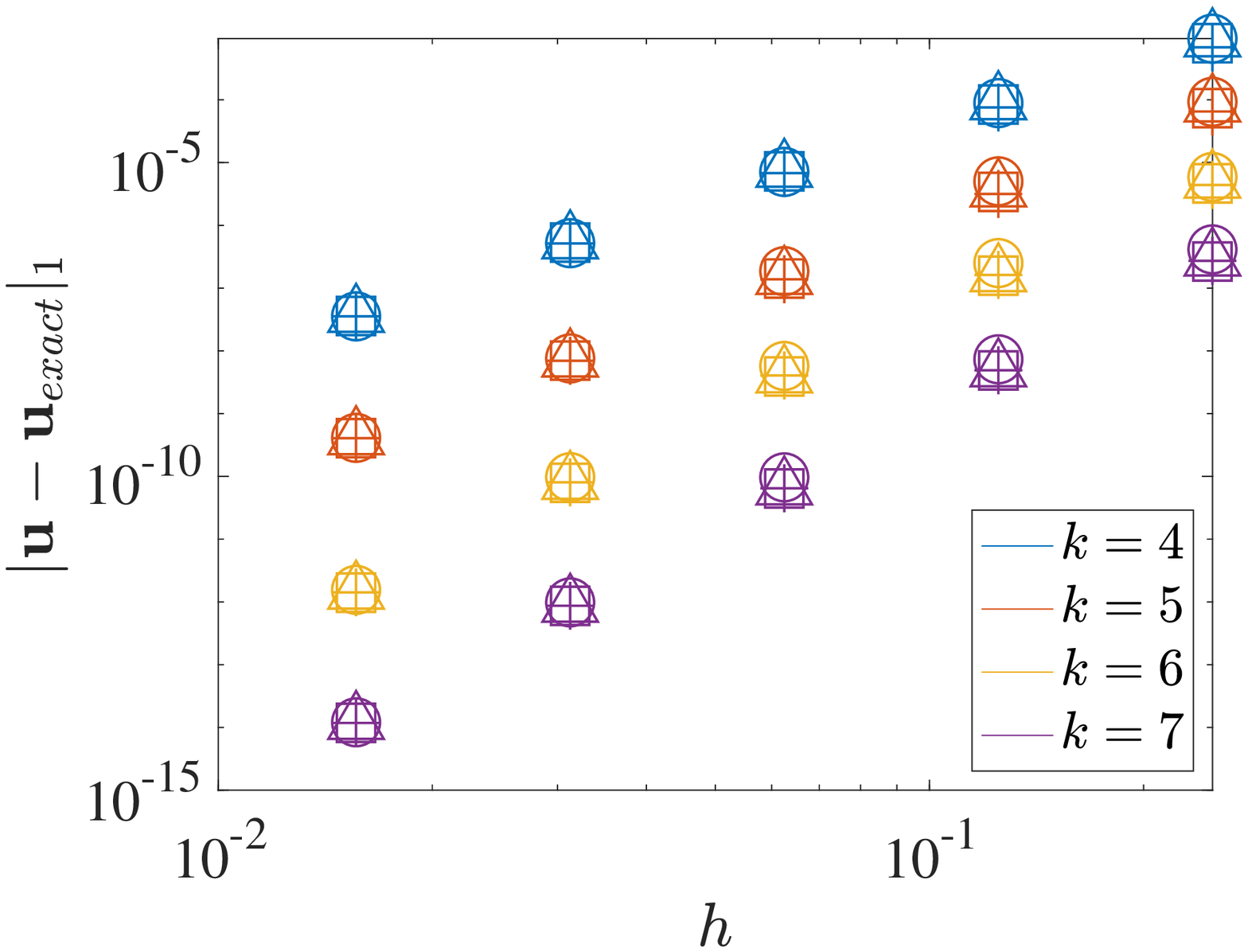}}\hfill
\subfloat[Pressure $H^1$ error]{\label{sfig:stokes_press_h1_edge}\includegraphics[width=.45\textwidth]{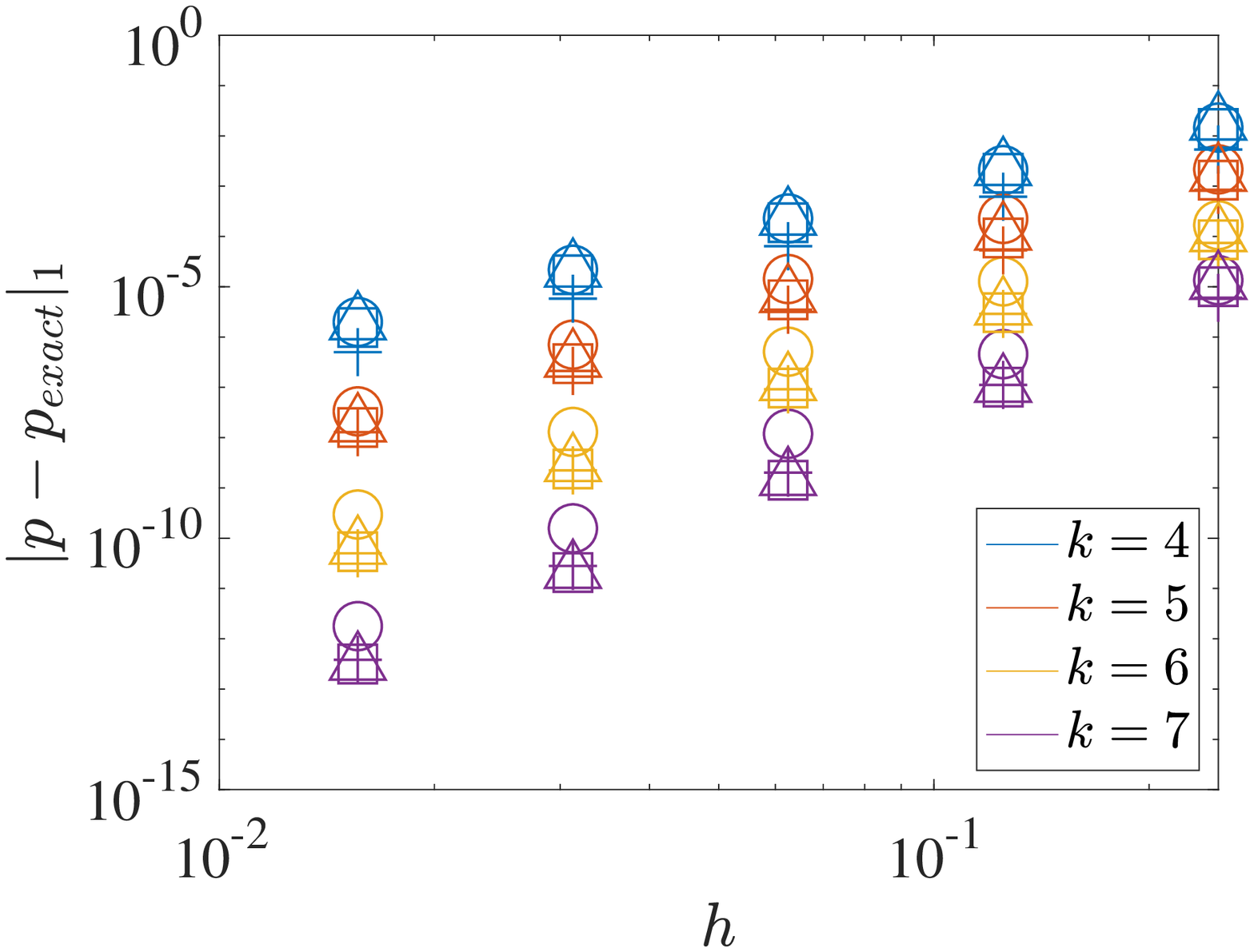}} \\
\caption{Effect of edge constant $C$ on errors in 2D Stokes manufactured solution using velocity-pressure form. $C = 0$ represented by circles, $C = 1$ represented by plus, $C = 10$ represented by squares, and $C = 100$ represented by triangles.}
\label{fig:stokes_2d_conv_edge}
\end{figure}

We also use this manufactured solution to examine the sensitivity of the collocation schemes to the constant $C$ used in the boundary term within the PSPG stabilization. Figure \ref{fig:stokes_2d_conv_edge} shows how the errors from the collocation scheme based on the two-field form of the equations varies as a function of $C$. For the most part we see that the errors are fairly insensitive to the value of the constant. The main noticeable difference is in the $H^1$ semi-norm of the pressure error, where increasing $C$ up to 10 reduces the error, and this effect is more noticeable for larger polynomial degrees. 

\begin{figure}
\centering
\subfloat[Velocity $L^2$ error]{\label{sfig:stokes_rot_vel_l2_edge}\includegraphics[width=.45\textwidth]{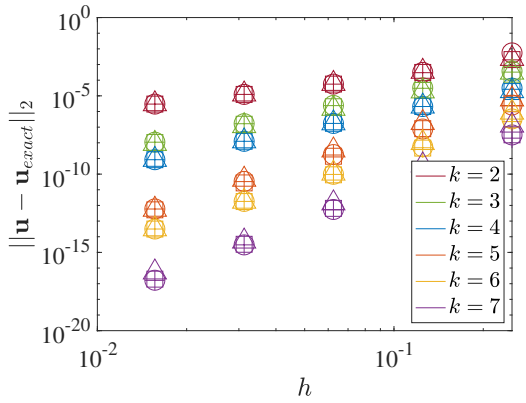}}\hfill
\subfloat[Pressure $L^2$ error]{\label{sfig:stokes_rot_press_l2_edge}\includegraphics[width=.45\textwidth]{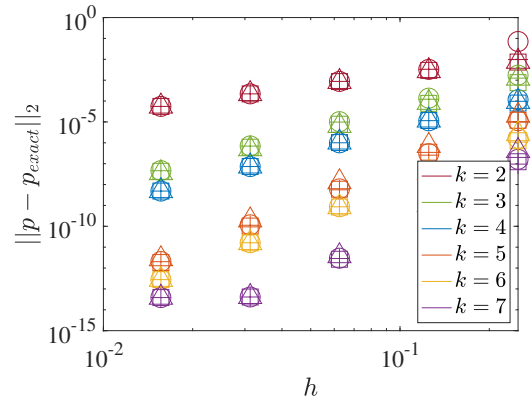}} \\
\subfloat[Velocity $H^1$ error]{\label{sfig:stokes_rot_vel_h1_edge}\includegraphics[width=.45\textwidth]{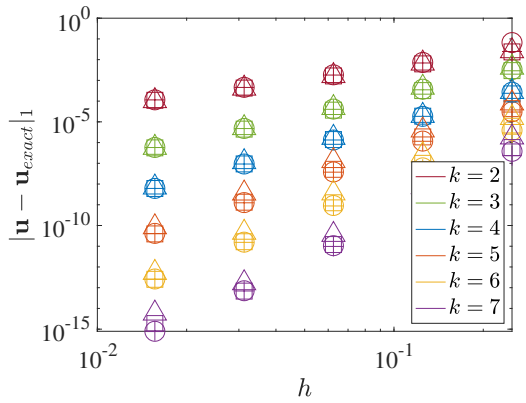}}\hfill
\subfloat[Pressure $H^1$ error]{\label{sfig:stokes_rot_press_h1_edge}\includegraphics[width=.45\textwidth]{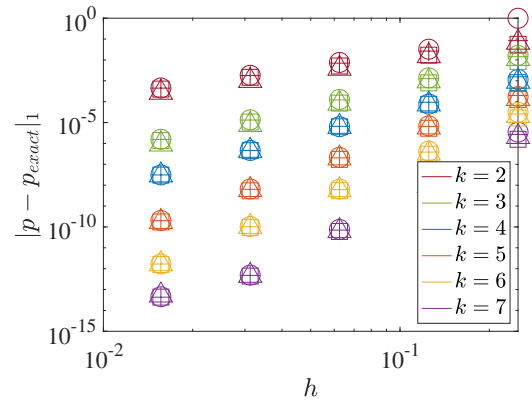}} \\
\caption{Effect of edge constant $C$ on errors in 2D Stokes manufactured solution using rotational form. $C = 0$ represented by circles, $C = 1$ represented by plus, $C = 10$ represented by squares, and $C = 100$ represented by triangles.}
\label{fig:stokes_rot_2d_conv_edge}
\end{figure}

Figure \ref{fig:stokes_rot_2d_conv_edge} details the results of the same test applied to the rotational form based collocation scheme. These errors are also fairly insensitive to $C$, including the $H^1$ pressure error. There are some small effects, such as the pressure being less accurate on the coarsest mesh with $k = 2$ and $C = 0$, and some small increases in the velocity errors for $k = 7$ and $C = 100$.

\subsubsection{2D Lid-Driven Cavity}

\begin{figure}
\centering
\subfloat[Velocity streamlines]{\label{sfig:stokes_cavity_sl}\includegraphics[width=.45\textwidth]{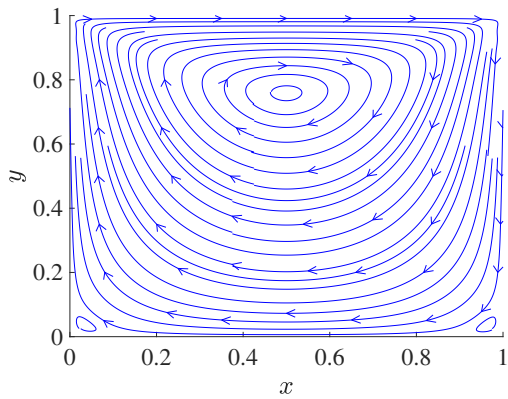}}\hfill
\subfloat[Pressure field]{\label{sfig:stokes_cavity_p}\includegraphics[width=.45\textwidth]{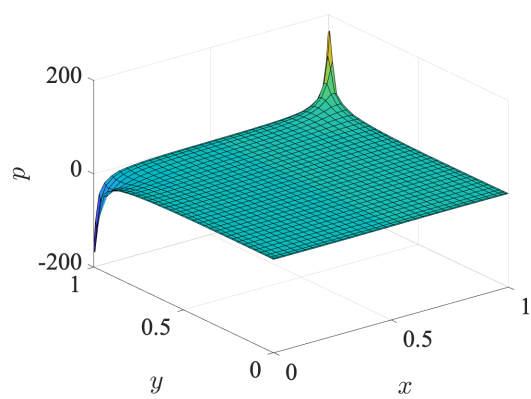}} \\
\caption{Lid driven cavity results with velocity-pressure formulation, 32$^2$ elements and $k$ = 4}
\label{fig:stokes_cavity}
\end{figure}

\begin{figure}
\centering
\subfloat[Velocity streamlines]{\label{sfig:stokes_rot_cavity_sl}\includegraphics[width=.45\textwidth]{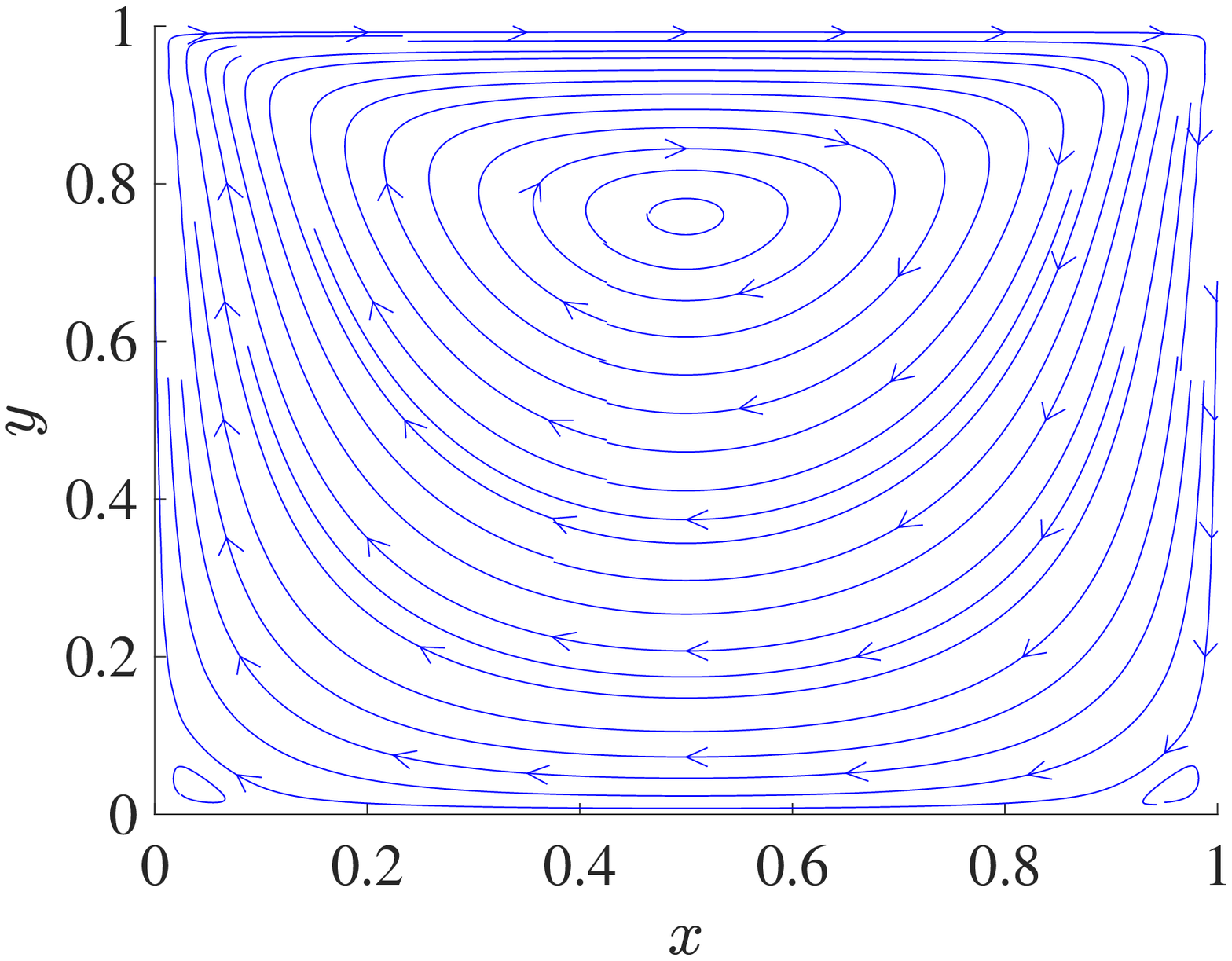}}\hfill
\subfloat[Pressure field]{\label{sfig:stokes_rot_cavity_p}\includegraphics[width=.45\textwidth]{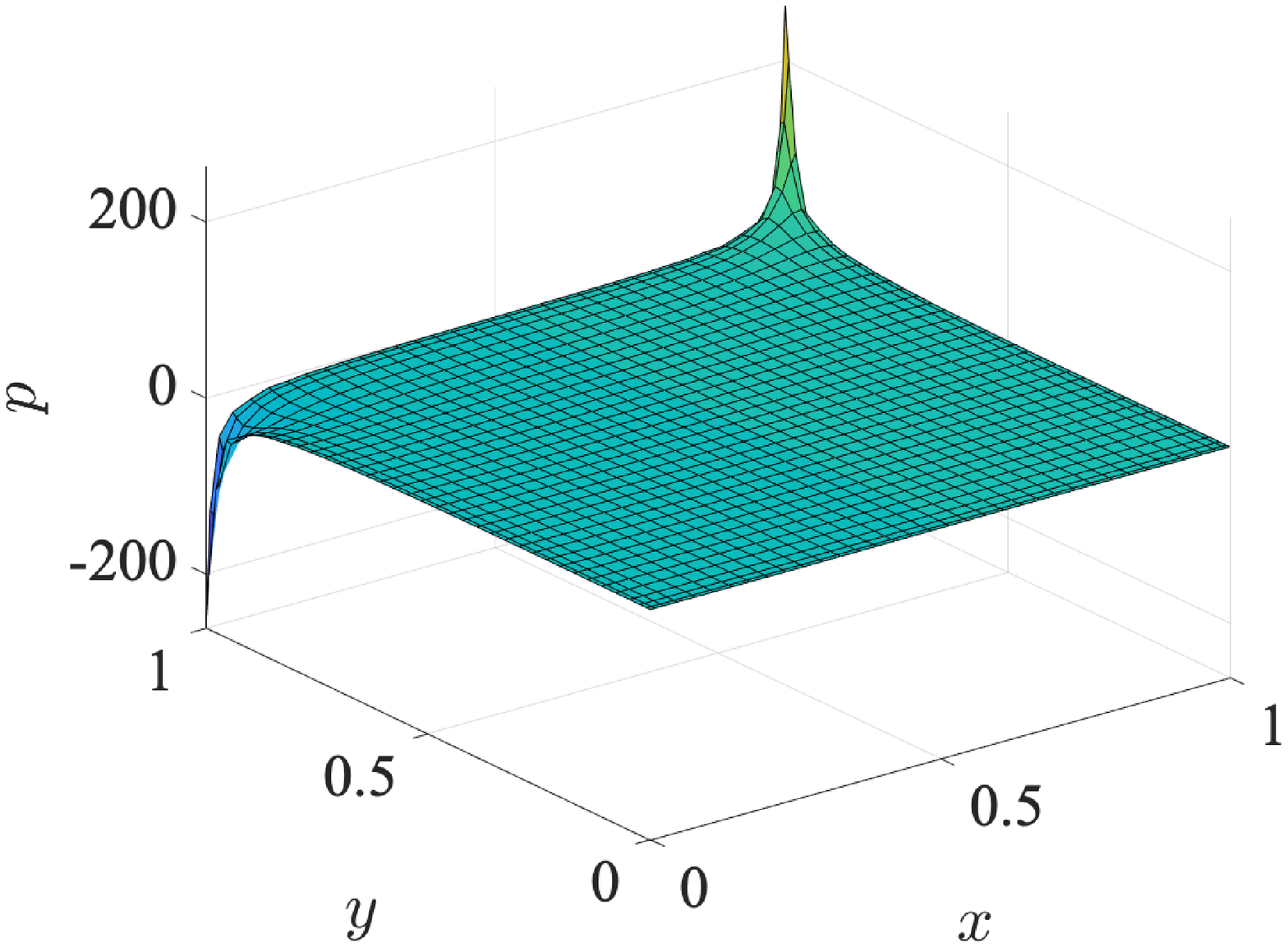}} \\
\caption{Lid driven cavity results with rotational formulation, 32$^2$ elements and $k$ = 4}
\label{fig:stokes_rot_cavity}
\end{figure}

We also test the stabilized collocation schemes for the Stokes equations by applying them to the 2D lid-driven cavity problem. In this case we consider the flow on the unit square with the top wall moving to the right at unit speed. Figure \ref{fig:stokes_cavity} shows the results with the velocity-pressure form of the Stokes equations and Figure \ref{fig:stokes_rot_cavity} shows the rotational result. Clearly both schemes are able to recover the symmetric velocity field that is expected in the Stokes case and the pressure field is not polluted with spurious oscillations.

\section{Incompressible Navier-Stokes Flow}

The final set of governing equations we consider in this paper are the nonlinear, steady, incompressible Navier-Stokes equations. In this case, numerical solutions can exhibit both advection and pressure instabilities and thus we combine the efforts of the previous sections to define collocation schemes which include both SUPG and PSPG stabilization. 

\subsection{The Strong Form of the Navier-Stokes Equations}

The standard, velocity-pressure form of the steady Navier-Stokes equations with Dirichlet boundary conditions is given by:

\bigskip

$$
\left\{ \hspace{5pt}
\parbox{5in}{
\noindent Given $\nu \in \mathbb{R}^+$, $\textbf{f} : \Omega \rightarrow \mathbb{R}^n$, and $\textbf{g} : \partial \Omega \rightarrow \mathbb{R}^n$, find $\mathbf{u} : \Omega \rightarrow \mathbb{R}^n$ and $p : \Omega \rightarrow \mathbb{R}$ such that:
\begin{align}
    -\nu \Delta \mathbf{u} + \mathbf{u} \cdot \nabla \mathbf{u} + \nabla p = \mathbf{f} \quad &\textup{in} \quad \Omega,\\
    \nabla \cdot \mathbf{u} = 0 \quad &\textup{in} \quad \Omega,\\
    \textbf{u} = \textbf{g} \quad &\textup{on} \quad \partial \Omega.
\end{align}
}
\right.
$$

\bigskip

\noindent Here we have defined the kinematic viscosity as $\nu$, while rest of the variables are defined as in the Stokes case. In the Navier-Stokes setting we also wish to explore the use of the rotational form, thus writing the governing equations of momentum and mass balance as a system of first-order differential equations. The rotational form of the Navier-Stokes equations in 2D is given as

\bigskip

$$
\left\{ \hspace{5pt}
\parbox{5in}{
\noindent Given $\nu \in \mathbb{R}^+$, $\textbf{f} : \Omega \rightarrow \mathbb{R}^2$, and $\textbf{g} : \partial \Omega \rightarrow \mathbb{R}^2$, find $\mathbf{u} : \Omega \rightarrow \mathbb{R}^2$, $P : \Omega \rightarrow \mathbb{R}$, and $\omega : \Omega \rightarrow \mathbb{R}$ such that:
\begin{align}
    \nu \nabla^{\perp} \omega + \omega \hat{\mathbf{k}} \times \mathbf{u} + \nabla P = \mathbf{f} \quad &\textup{in} \quad \Omega, \label{eq:NS_rot_mom}\\
    \nabla \cdot \mathbf{u} = 0 \quad &\textup{in} \quad \Omega,\\
    \omega - \nabla \times \mathbf{u} = 0 \quad &\textup{in} \quad \Omega,\\
    \textbf{u} = \textbf{g} \quad &\textup{on} \quad \partial \Omega.
\end{align}
}
\right.
$$

\bigskip

In these equations $P$ is the total pressure, defined in relation to the kinematic pressure via via $P = p + \frac{1}{2}\textbf{u}\cdot\textbf{u}$. This is different from the Stokes setting, where the pressure variable had the same meaning in the two- and three-field cases. We have also used the convention that the velocity $\mathbf{u}$ is defined in the $\hat{\mathbf{i}}$ and $\hat{\mathbf{j}}$ (or $x$ and $y$) directions, while the scalar vorticity value $\omega$ defines a vector vorticity purely in the $\hat{\mathbf{k}}$ direction, orthogonal to the velocity components. The $\hat{\mathbf{k}}$ component of Equation \eqref{eq:NS_rot_mom} is trivially satisfied for 2D flows. 

\subsection{Stabilization for Navier-Stokes Flow}

Similar to the previous sections, we start by stating a stabilized weak form upon which a Galerkin finite element method would be based. Utilizing the same trial and test space definitions as in the Stokes case, the stabilized weak form of the velocity-pressure Navier-Stokes equations is given by

\bigskip

$$
\left\{ \hspace{5pt}
\parbox{5in}{
Given $\nu \in \mathbb{R}^+$, $\textbf{f} : \Omega \rightarrow \mathbb{R}^n$, and $\textbf{g} : \partial \Omega \rightarrow \mathbb{R}^n$, find $\mathbf{u}^h \in \mathbf{S}^h_g$ and $p^h \in \Tilde{S}^h$ such that :
\begin{equation}
\begin{split}
    \int_{\Omega} (\nu \nabla \mathbf{w}^h \cdot \nabla \mathbf{u}^h + &\mathbf{w}^h \cdot (\mathbf{u}^h \cdot \nabla \mathbf{u}^h) - p^h \nabla \cdot \mathbf{w}^h) d\Omega \\ -  &\int_{\Omega} \mathbf{w}^h \cdot \mathbf{f} d \Omega\\ +  &\int_{\Omega} \tau_{SUPG}( \mathbf{u}^h \cdot \nabla \mathbf{w}^h) \cdot \mathbf{R}^h d \Omega\\ +  &\int_{\Omega} \tau_{GD} (\nabla \cdot \mathbf{w}^h) (\nabla \cdot \mathbf{u}^h) d \Omega = 0 \quad \forall \mathbf{w}^h \in \mathbf{S}^h_0,  \label{eq:ns_mom_weak}
\end{split}
\end{equation}
\begin{equation}
    \int_{\Omega} q^h (\nabla \cdot \mathbf{u}^h) d \Omega + \int_{\Omega} \tau_{PSPG} (\nabla q^h \cdot \mathbf{R}^h) d \Omega = 0 \quad \forall q^h \in \Tilde{S}^h.  \label{eq:ns_mass_weak}
\end{equation}
}
\right.
$$

\bigskip

The third integral in Equation \eqref{eq:ns_mom_weak} is the SUPG stabilization term and the following term defines the Grad-Div stabilization \cite{olshanskii2009grad}, which is sometimes included in the definition of SUPG. The PSPG term in Equation \eqref{eq:ns_mass_weak} takes the same form as in the Stokes equations, but with the nonlinear momentum equation defining the residual

\begin{equation}
    \mathbf{R}^h = -\nu \Delta \mathbf{u}^h + \mathbf{u}^h \cdot \nabla \mathbf{u}^h + \nabla p^h - \mathbf{f}.
\end{equation}

The integration by parts of the continuity equation is done in the same manner as in the Stokes equations. Integrating the stabilized momentum equation by parts and again noting that the test function satisfies homogeneous Dirichlet boundary conditions yields

\begin{equation}
\begin{split}
    \int_{\Omega} \mathbf{w}^h \cdot &(-\nu \Delta \mathbf{u}^h + \mathbf{u}^h \cdot \nabla \mathbf{u}^h + \nabla p^h - \mathbf{f})d \Omega \\ &- \int_{\Omega} \mathbf{w}^h \cdot ( \nabla \cdot (\tau_{SUPG} \mathbf{u}^h \otimes \mathbf{R}^h ) ) d \Omega \\ &- \int_{\Omega} \mathbf{w}^h \cdot \nabla (\tau_{GD} \nabla \cdot \mathbf{u}^h) d \Omega = 0.
\end{split}
\end{equation}

With the above manipulations, we are ready to define a collocation scheme for the incompressible Navier-Stokes equations. Once again, let $\mathbf{w}^h$ and $q^h$ be Dirac-delta functions located at the Greville points $\{\mathbf{x}\}_{i=1}^n$ of the spline space $S^h$. This yields the following fully discrete collocation scheme, where boundary integrals in the continuity equation have again been treated using the Enhanced Collocation treatment

\bigskip

$$
\left\{ \hspace{5pt}
\parbox{5in}{
\noindent Find $\mathbf{u}^h \in \mathbf{S}^h$ and $p^h \in \Tilde{S}^h$ such that, for $i = 1, ..., n$

\begin{equation}
\begin{split}
    (-\nu \Delta \mathbf{u}^h + &\mathbf{u}^h \cdot \nabla \mathbf{u}^h + \nabla p^h - \mathbf{f} \\ & - \nabla \cdot (\tau_{SUPG} \mathbf{u}^h \otimes \mathbf{R}^h ) \\ &- \nabla (\tau_{GD} \nabla \cdot \mathbf{u}^h))(\mathbf{x}_i) = \mathbf{0} \quad \forall \mathbf{x}_i \in \Omega,
\end{split}
\end{equation}
\begin{equation}
    (\nabla \cdot \mathbf{u}^h - \nabla \cdot (\tau_{PSPG} \mathbf{R}^h))(\mathbf{x}_i) = 0 \quad \forall \mathbf{x}_i \in \Omega,
\end{equation}
\begin{equation}
    (\nabla \cdot \mathbf{u}^h - \nabla \cdot (\tau_{PSPG} \mathbf{R}^h) + \frac{C}{h_b}(\tau_{PSPG} \mathbf{R}^h \cdot \mathbf{n}))(\mathbf{x}_i) = 0 \quad \forall \mathbf{x}_i \in \partial \Omega,
\end{equation}
\begin{equation}
    \mathbf{u}^h(\textbf{x}_i) = \mathbf{g}(\mathbf{x}_i) \quad \forall \textbf{x}_i \in \partial \Omega.
\end{equation}

}
\right.
$$

\bigskip

\noindent This defines a nonlinear system of equations which can be solved with a Newton-Raphson method. All that remains is the specification of the stabilization constants. For our purposes we select

\begin{equation}
    \tau_{SUPG} = \tau_{PSPG} = \frac{1}{\sqrt{(\frac{2 ||\mathbf{u}||}{h})^2 + (\frac{C_3 \nu}{h^2})^2}},
    \label{eq:tau_supg_pspg}
\end{equation}

\noindent and 

\begin{equation}
    \tau_{GD} = \frac{2 h^2}{\nu},
\end{equation}

\noindent with $C_3$ equal to 4 and spline interpolations of the $\tau$ values as used previously to compute derivative terms.

We now follow similar steps to construct the stabilized collocation scheme based on the rotational form of the Navier Stokes equations. We remark at this point that Grad-Div stabilization has been shown in \cite{layton2009accuracy} to be especially important in incompressible flow computations based on the rotational form of the Navier-Stokes equations. Using the same test and trial function spaces as above (with the minor modification that the pressure now denotes the total pressure rather than the kinematic pressure), a sample stabilized weak form is written as

\bigskip

$$
\left\{ \hspace{5pt}
\parbox{5in}{
Given $\nu \in \mathbb{R}^+$, $\textbf{f} : \Omega \rightarrow \mathbb{R}^n$, and $\textbf{g} : \partial \Omega \rightarrow \mathbb{R}^n$, find $\mathbf{u} \in \mathbf{S}^h_g$ and $P \in \Tilde{S}^h$ such that :
\begin{equation}
\begin{split}
    \int_{\Omega} (\nu \mathbf{w}^h \cdot \nabla^{\perp} \omega^h + &\mathbf{w}^h \cdot (\omega^h \hat{\mathbf{k}} \times \textbf{u}^h) - P^h \nabla \cdot \mathbf{w}^h) d\Omega \\ &- \int_{\Omega} \mathbf{w}^h \cdot \mathbf{f} d \Omega \\ &+ \int_{\Omega} \tau_{SUPG} (\omega^h \hat{\mathbf{k}} \times \mathbf{w}^h) \cdot \mathbf{R}^h d \Omega \\ &+ \int_{\Omega} \tau_{GD} (\nabla \cdot \mathbf{w}^h) (\nabla \cdot \mathbf{u}^h) d \Omega = 0 \quad \forall \mathbf{w}^h \in \mathbf{S}^h_0,
\end{split}
\end{equation}
\begin{equation}
    \int_{\Omega} q^h (\nabla \cdot \mathbf{u}^h) d \Omega + \int_{\Omega} \tau_{PSPG} (\nabla q^h \cdot \mathbf{R}^h) d \Omega = 0 \quad \forall q^h \in \Tilde{S}^h, 
\end{equation}
}
\right.
$$

\bigskip

\noindent where now the momentum residual is defined as 

\begin{equation}
    \mathbf{R}^h = \nu \nabla^{\perp} \omega^h + \omega^h \hat{\mathbf{k}} \times \mathbf{u}^h + \nabla P^h - \mathbf{f}.
\end{equation}

Using properties of the scalar triple product on the SUPG term and reverse integrating by parts transforms the momentum equation into

\begin{equation}
\begin{split}
    \int_{\Omega} \mathbf{w}^h \cdot &(\nu \nabla^{\perp} \omega^h + \omega^h \hat{\mathbf{k}} \times \textbf{u}^h + \nabla P^h - \mathbf{f}) d\Omega \\ &+ \int_{\Omega} \mathbf{w}^h \cdot (\tau_{SUPG} (\mathbf{R}^h \times \omega^h \hat{\mathbf{k}}))d \Omega \\ &- \int_{\Omega} \mathbf{w}^h \cdot \nabla (\tau_{GD} \nabla \cdot \mathbf{u}^h) d \Omega = 0.
\end{split}
\end{equation}

With the continuity equation treated as before, we define our last stabilized collocation scheme by

\bigskip

$$
\left\{ \hspace{5pt}
\parbox{5in}{
\noindent Find $\mathbf{u}^h \in \mathbf{S}^h$, $p^h \in \Tilde{S}^h$, and $\omega^h \in S^h$ such that, for $i = 1, ..., n$

\begin{equation}
\begin{split}
    (\nu \nabla^{\perp} \omega^h + &\omega^h \hat{\mathbf{k}} \times \textbf{u}^h + \nabla P^h - \mathbf{f} \\ &+ \tau_{SUPG} (\mathbf{R}^h \times \omega^h \hat{\mathbf{k}}) \\ &- \nabla (\tau_{GD} \nabla \cdot \mathbf{u}^h))(\mathbf{x}_i) = \mathbf{0} \quad \forall \mathbf{x}_i \in \Omega,
\end{split}
\end{equation}
\begin{equation}
    (\nabla \cdot \mathbf{u}^h - \nabla \cdot (\tau_{PSPG} \mathbf{R}^h))(\mathbf{x}_i) = 0 \quad \forall \mathbf{x}_i \in \Omega,
\end{equation}
\begin{equation}
    (\nabla \cdot \mathbf{u}^h - \nabla \cdot (\tau_{PSPG} \mathbf{R}^h) + \frac{C}{h_b}(\tau_{PSPG} \mathbf{R}^h \cdot \mathbf{n}))(\mathbf{x}_i) = 0 \quad \forall \mathbf{x}_i \in \partial \Omega,
\end{equation}
\begin{equation}
    \mathbf{u}^h(\textbf{x}_i) = \mathbf{g}(\mathbf{x}_i) \quad \forall \textbf{x}_i \in \partial \Omega.
\end{equation}

}
\right.
$$

\bigskip

\noindent This nonlinear system is also solved by a Newton-Raphson approach, and the stabilization constant $\tau_{SUPG} = \tau_{PSPG}$ is selected to be the same as in the velocity-pressure scheme but divided by 10. We have seen that in some high Reynolds number cases we need to weaken the SUPG/PSPG stabilization constant in order to achieve convergence of the Newton solver, suggesting that perhaps the scaling of $\tau$ with Reynolds number should be different for the rotational form. We remark that the stability of the SUPG method when applied to the rotational form of the Navier-Stokes equations is analytically studied in \cite{lube2002stable}, but it has rarely been used in numerical computations to our knowledge. 

\subsubsection{A Note on Neumann Boundary Conditions}

Throughout this work we have focused on pure Dirichlet problems for concreteness and ease of definition. However, in the following section, we also solve a problem including Neumann boundary conditions to ensure that the schemes remain viable. We assume the Neumann boundary condition will be of form 

\begin{equation}
    -\nu \nabla \mathbf{u} \cdot \mathbf{n} + p \mathbf{n} = \mathbf{h}
    \label{eq:neumann_bc}
\end{equation}

\noindent where $h$ is the provided boundary data. Then, for every collocation point along the Neumann boundary $\partial \Omega_N$ we simply collocate Equation \eqref{eq:neumann_bc}, leading to

\begin{equation}
    (-\nu \nabla \mathbf{u} \cdot \mathbf{n} + p \mathbf{n})(\mathbf{x_i}) = \mathbf{h}(\mathbf{x}_i) \quad \forall \textbf{x}_i \in \partial \Omega_N
\end{equation}

The rest of the scheme remains the same outside the Neumann boundary. There are other approaches to enforcing Neumann boundary conditions which are summarized in \cite{deLorenzis2015isogeometric}. In particular, the Enhanced Collocation approach, where a penalty term is added to the momentum equations to enforce the boundary conditions, is shown to be much more accurate than the collocated approach shown here in some cases. However, the simple collocation approach has produced sufficient quality results here and its simplicity allows us to focus on the stabilization terms. 

\subsection{Navier-Stokes Results}

In this last major section, we apply our collocation schemes to some classical Navier-Stokes test cases. The same manufactured solution as in the Stokes case is tested first, and we show that the inclusion of nonlinear convection and SUPG stabilization minimally impacts the errors. Then we move on to Kovasznay flow where we also obtain convergence rates at least as fast as expected. Finally, we simulate the lid-driven cavity at variety of Reynolds numbers and show that these results agree well with standard benchmark data.

\subsubsection{Manufactured Solution and Convergence Rates in 2D}

To ensure that we retain the same accuracy in as in the linear Stokes case, we start by considering the same manufactured solution as above with $\nu = 1$. Then the forcing term takes the form 

\begin{equation}
    \textbf{f} = -\nu \Delta \bar{\textbf{u}} + \bar{\mathbf{u}} \cdot \nabla \bar{\mathbf{u}} + \nabla \bar{p}, 
\end{equation}

\noindent in the velocity-pressure formulation, while for the rotational form we use 

\begin{equation}
    \textbf{f} = \nu \nabla^\perp \bar{\omega} + \overline{\omega} \times \bar{\textbf{u}} + \nabla \bar{P}.
\end{equation}

\begin{figure}
\centering
\subfloat[Velocity $L^2$ error]{\label{sfig:ns_vel_l2}\includegraphics[width=.45\textwidth]{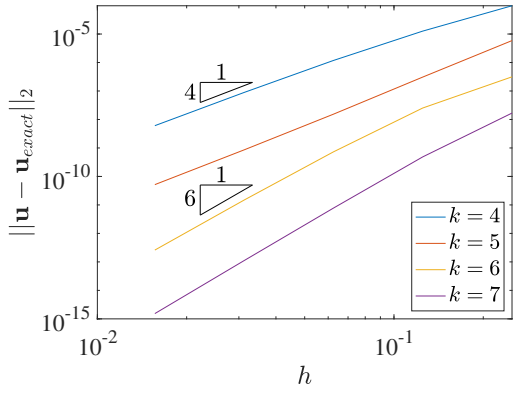}}\hfill
\subfloat[Pressure $L^2$ error]{\label{sfig:ns_press_l2}\includegraphics[width=.45\textwidth]{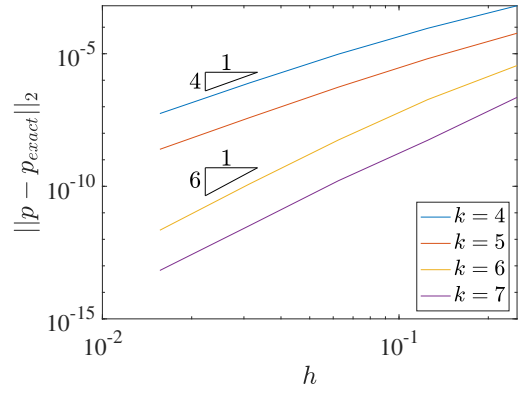}} \\
\subfloat[Velocity $H^1$ error]{\label{sfig:ns_vel_h1}\includegraphics[width=.45\textwidth]{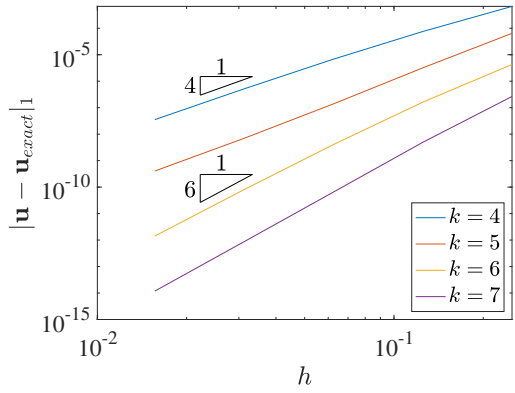}}\hfill
\subfloat[Pressure $H^1$ error]{\label{sfig:ns_press_h1}\includegraphics[width=.45\textwidth]{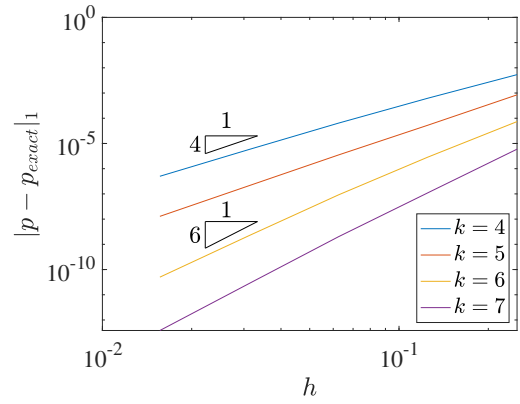}} \\
\caption{Errors in 2D Navier-Stokes manufactured solution using velocity-pressure form}
\label{fig:ns_2d_conv}
\end{figure}

\begin{figure}
\centering
\subfloat[Velocity $L^2$ error]{\label{sfig:ns_rot_vel_l2}\includegraphics[width=.45\textwidth]{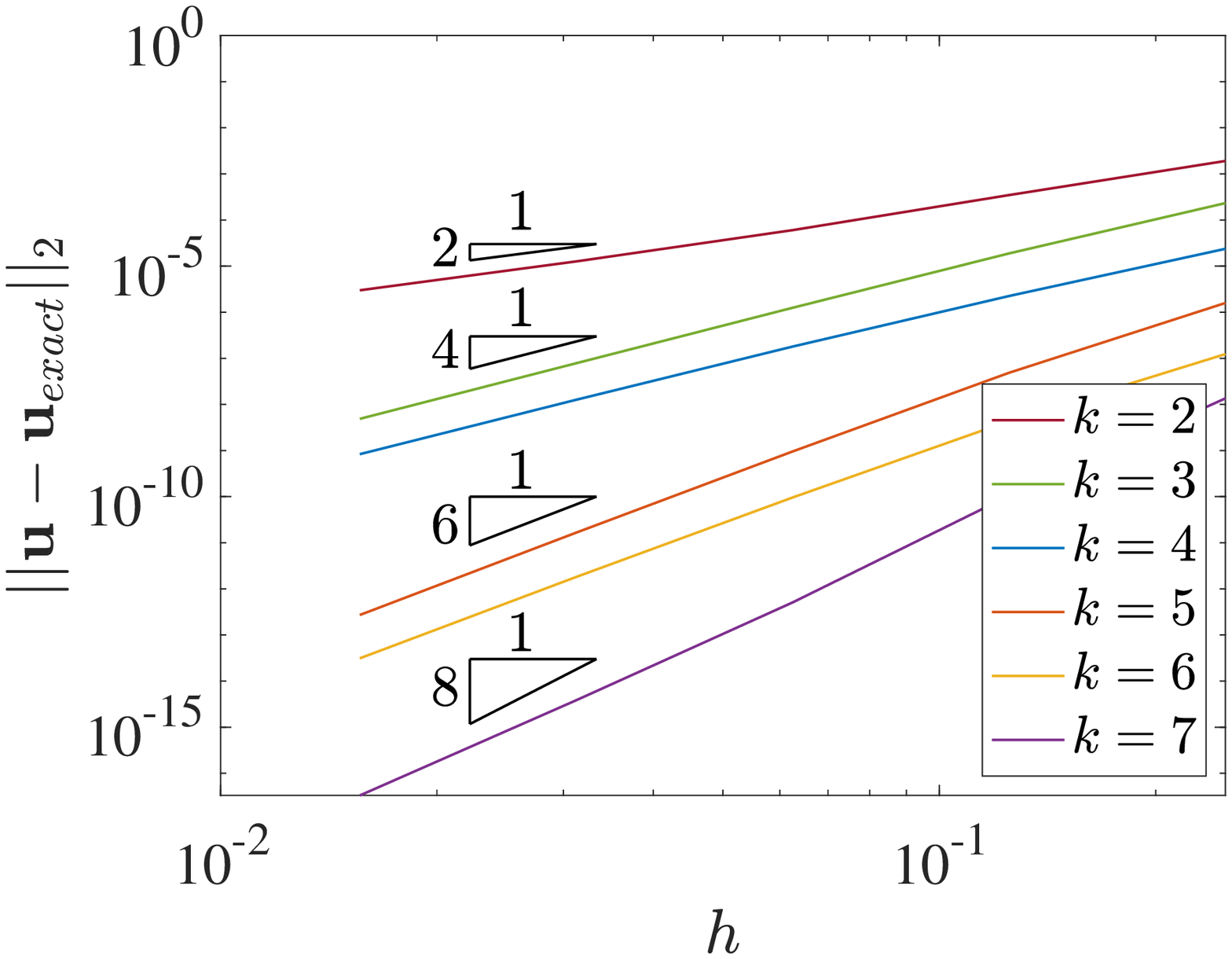}}\hfill
\subfloat[Pressure $L^2$ error]{\label{sfig:ns_rot_press_l2}\includegraphics[width=.45\textwidth]{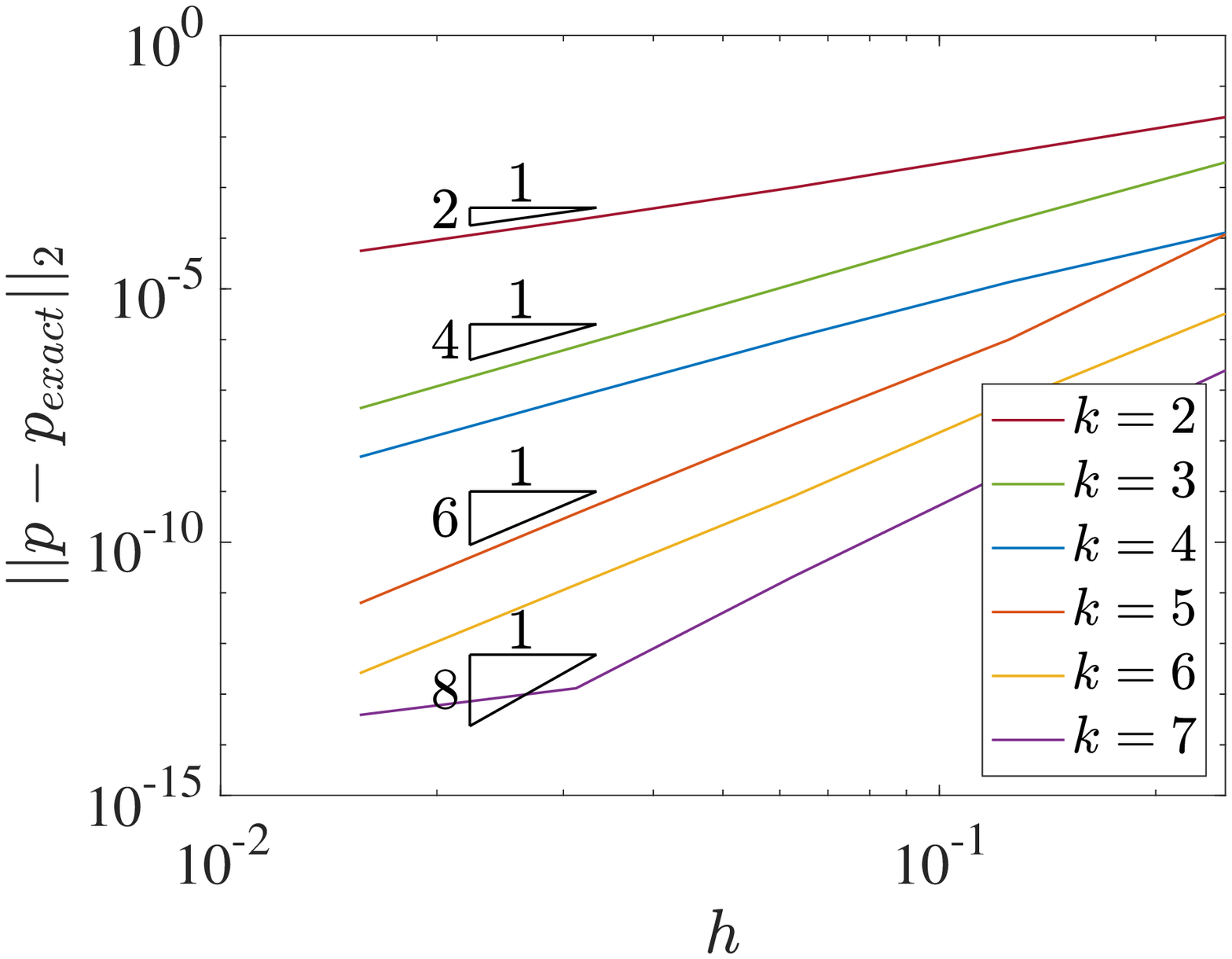}} \\
\subfloat[Velocity $H^1$ error]{\label{sfig:ns_rot_vel_h1}\includegraphics[width=.45\textwidth]{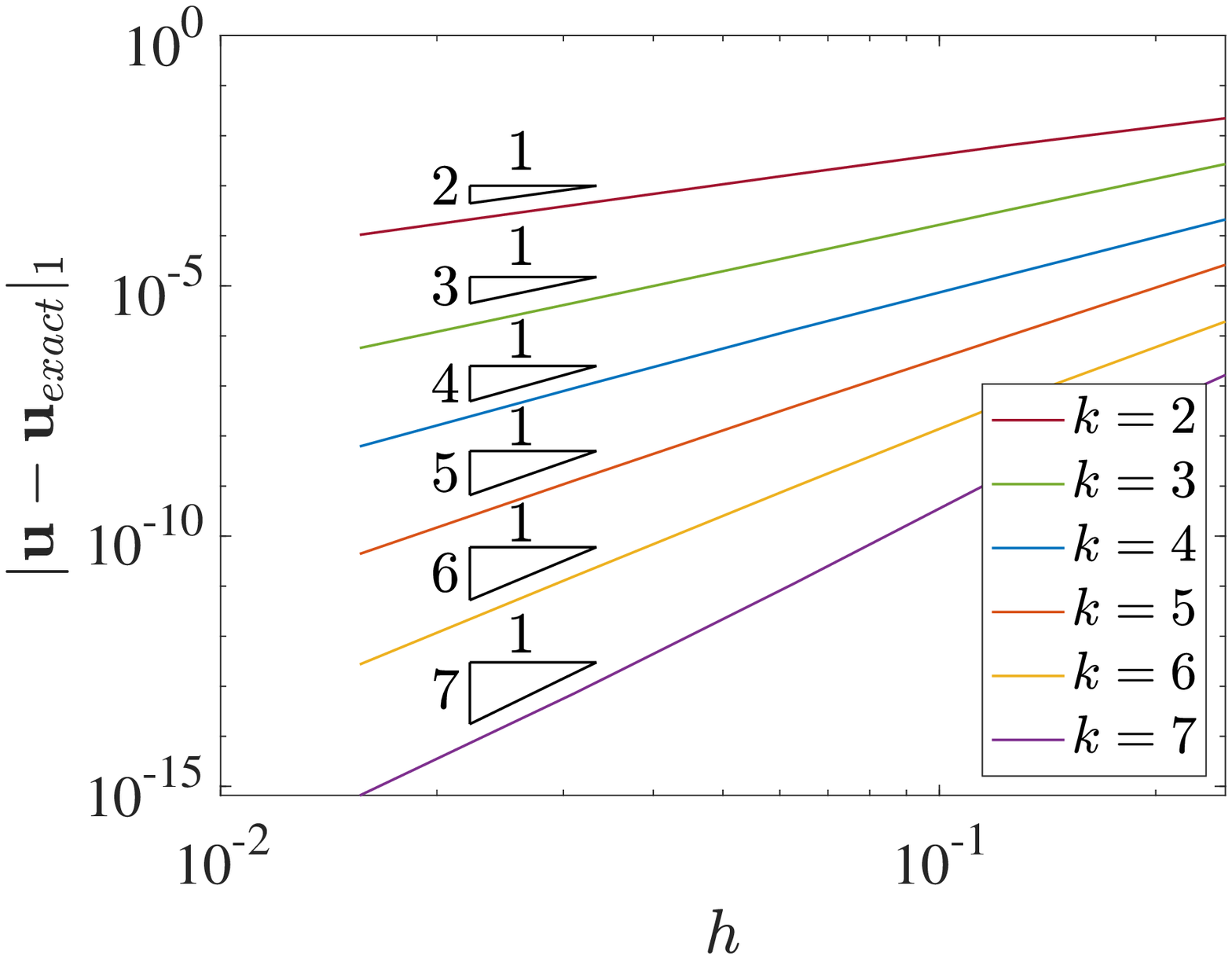}}\hfill
\subfloat[Pressure $H^1$ error]{\label{sfig:ns_rot_press_h1}\includegraphics[width=.45\textwidth]{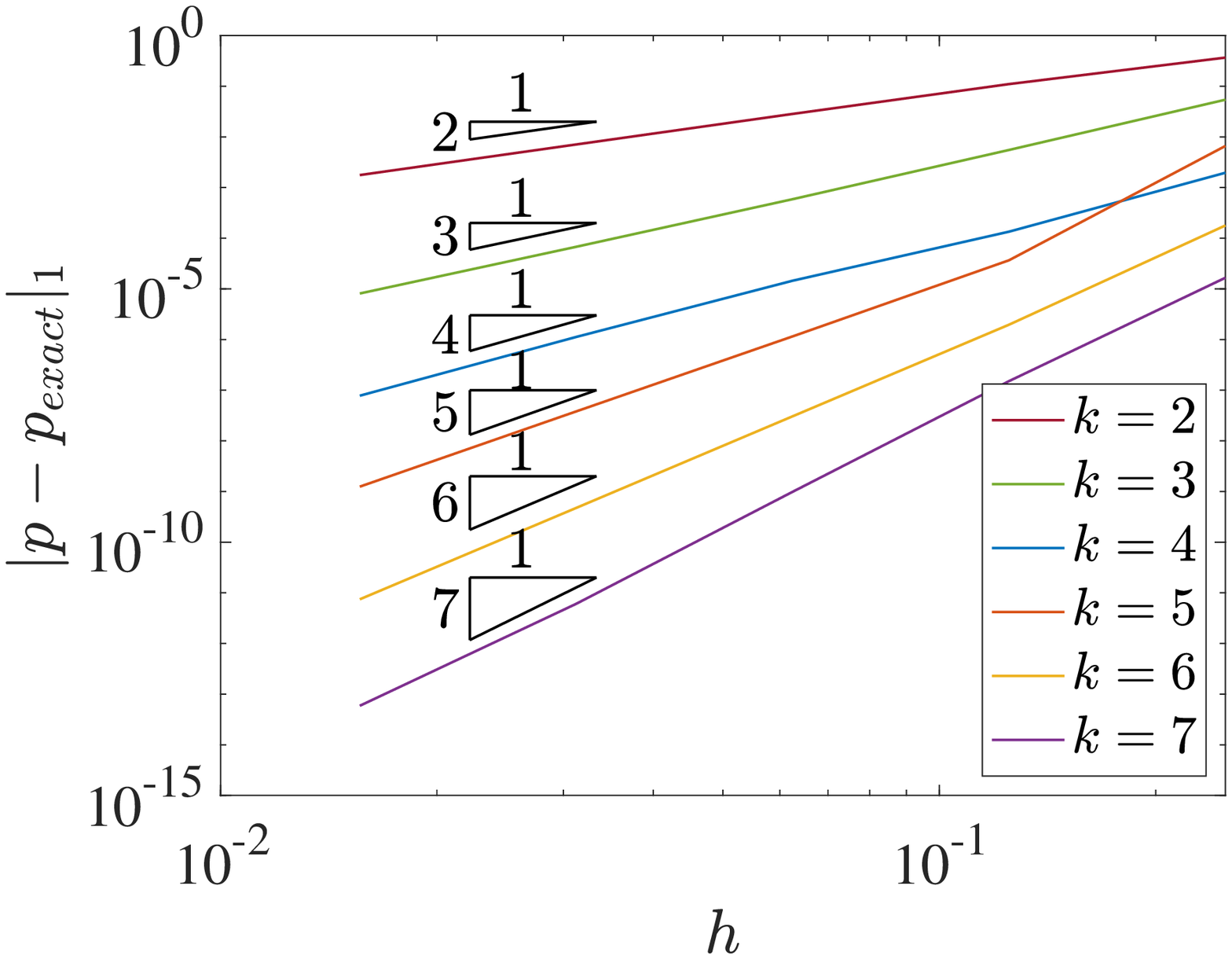}} \\
\caption{Errors in 2D Navier-Stokes manufactured solution using rotational form}
\label{fig:ns_rot_2d_conv}
\end{figure}

Figure \ref{fig:ns_2d_conv} shows the convergence of the numerical errors as a function of basis polynomial degree and mesh resolution for the velocity-pressure form. Figure \ref{fig:ns_rot_2d_conv} details the equivalent results obtained with the rotational form of the Navier-Stokes equations. Note that we are presenting the errors in the kinematic pressure $p$ in the rotational form results. The inclusion of nonlinear convection and SUPG stabilization has minimal effect on the convergence rates when compared with those achieved in the Stokes equation setting. The magnitudes of the errors in the velocity-pressure scheme are almost identical to those of the velocity-pressure Stokes equation schemes. We do see an increase in error magnitude in some cases in the Navier-Stokes results using the rotational form, however. We have found that this is likely due to our choice of stabilization parameter $\tau_{SUPG} = \tau_{PSPG}$. If this parameter is increased to the expression used for the velocity-pressure formulation then these errors agree better with those seen in the Stokes case. However, we have found that this choice stalls the convergence of the Newton solver for high Reynolds number versions of the lid-driven cavity test case shown below. Thus we choose to present the results obtained for all examples with a consistent choice of stabilization parameter.

\subsubsection{Kovasznay Flow}

Next we move on to consider Kovasznay flow, which represents the flow behind an infinite two dimensional grid. This problem also possesses an exact solution and is also often used to perform convergence studies for Navier-Stokes discretizations. The analytically derived exact solution \cite{kovasznay1948laminar} is given by 

\begin{equation}
    \bar{\textbf{u}} = \left[
    \begin{array}{c}
    1-e^{\lambda x} \cos(2\pi y) \\
    \frac{\lambda}{2\pi}e^{\lambda x} \sin(2\pi y)
    \end{array}
    \right],
\end{equation}

\begin{equation}
    \bar{p} = \frac{1 - e^{\lambda x}}{2},
\end{equation}

\noindent with the parameter $\lambda$ defined as

\begin{equation}
    \lambda = \frac{Re}{2} - \sqrt{\frac{Re^2}{4} + 4 \pi^2},
\end{equation}

\noindent and the Reynolds number $Re$ defined as the reciprocal of the viscosity $\nu$. We solve over the domain $[-0.5, 1] \times [-0.5,  0.5]$ with $Re = 40$. We use the exact solution as a Dirichlet boundary condition on the left, top, and  bottom boundaries, and define a Neumann boundary condition at the outflow on the right boundary. The forcing terms take the same forms as in the manufactured solution above. 

\begin{figure}
\centering
\subfloat[Velocity $L^2$ error]{\label{sfig:ns_kovasznay_vel_l2}\includegraphics[width=.45\textwidth]{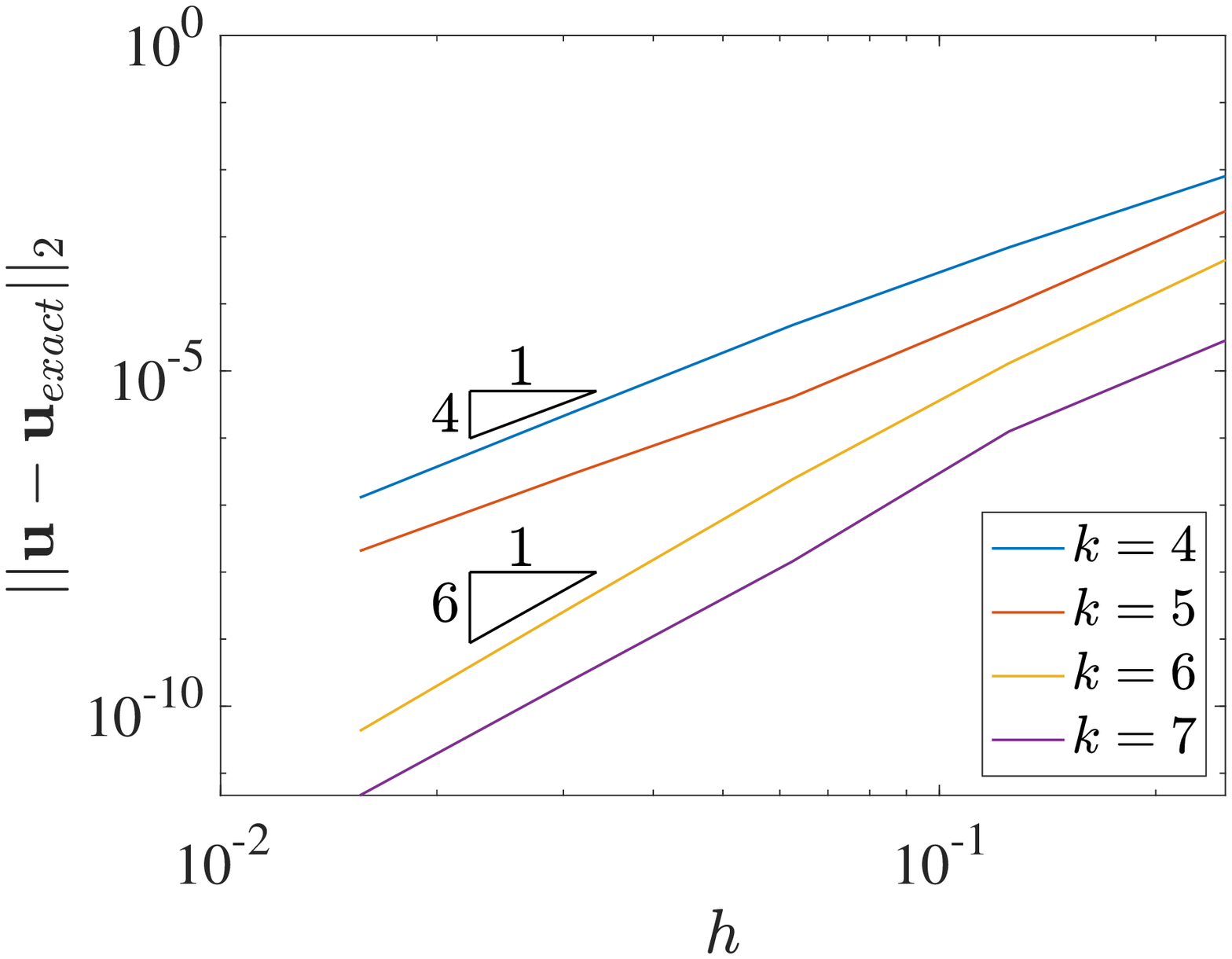}}\hfill
\subfloat[Pressure $L^2$ error]{\label{sfig:ns_kovasznay_press_l2}\includegraphics[width=.45\textwidth]{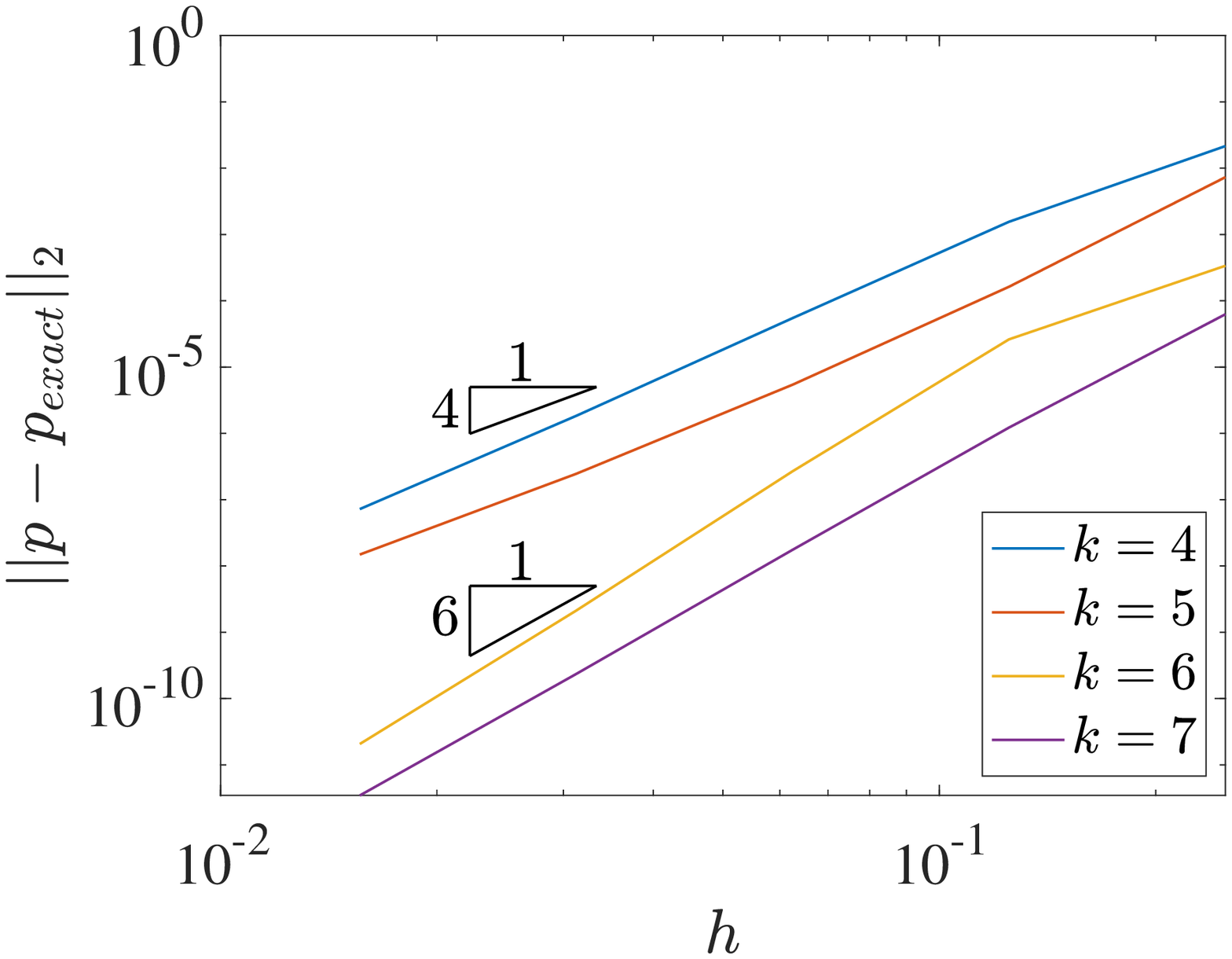}} \\
\subfloat[Velocity $H^1$ error]{\label{sfig:ns_kovasznay_vel_h1}\includegraphics[width=.45\textwidth]{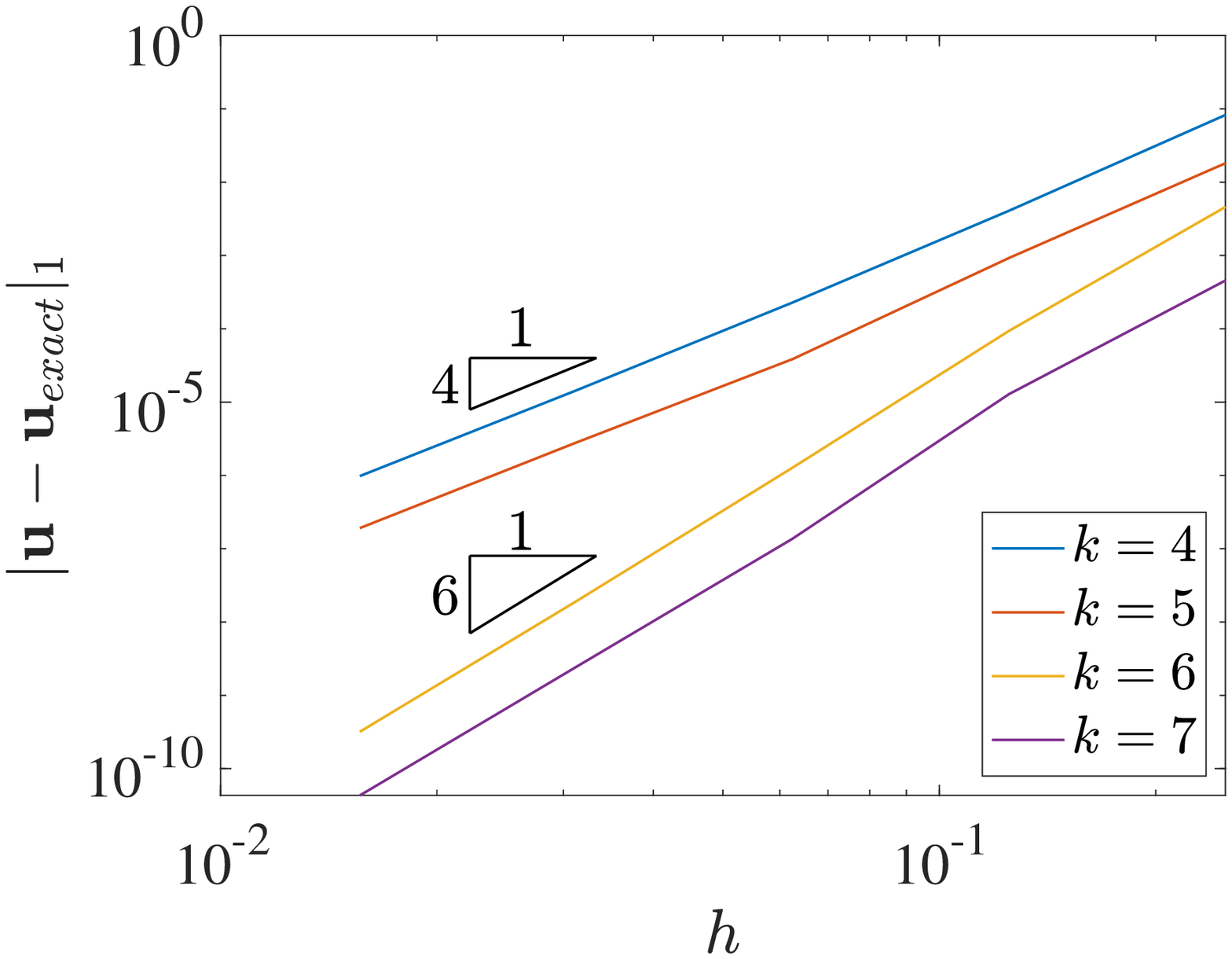}}\hfill
\subfloat[Pressure $H^1$ error]{\label{sfig:ns_kovasznay_press_h1}\includegraphics[width=.45\textwidth]{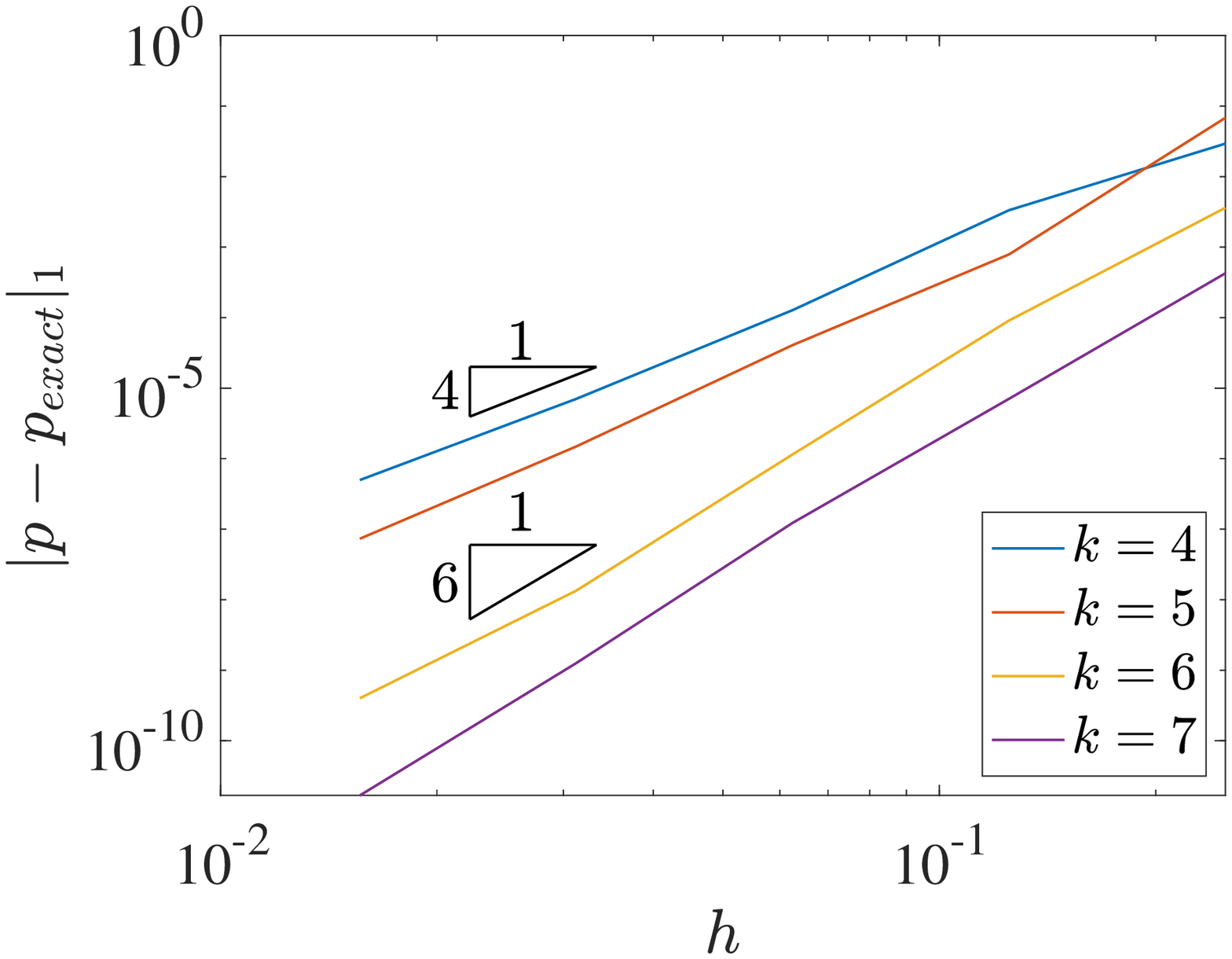}} \\
\caption{Errors in Kovasznay flow using velocity-pressure form}
\label{fig:ns_kovasznay_2d_conv}
\end{figure}

Figure \ref{fig:ns_kovasznay_2d_conv} shows the convergence of errors using the two-field formulation. The expected rates of $k$ and $k-1$ for even and odd $k$, respectively, are again recovered in this case. 

\begin{figure}
\centering
\subfloat[Velocity $L^2$ error]{\label{sfig:ns_rot_kovasznay_vel_l2}\includegraphics[width=.45\textwidth]{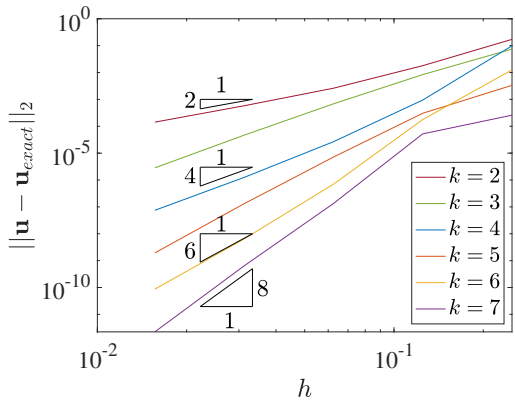}}\hfill
\subfloat[Pressure $L^2$ error]{\label{sfig:ns_rot_kovasznay_press_l2}\includegraphics[width=.45\textwidth]{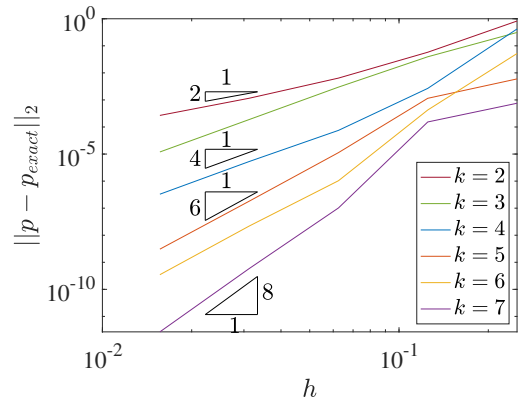}} \\
\subfloat[Velocity $H^1$ error]{\label{sfig:ns_rot_kovasznay_vel_h1}\includegraphics[width=.45\textwidth]{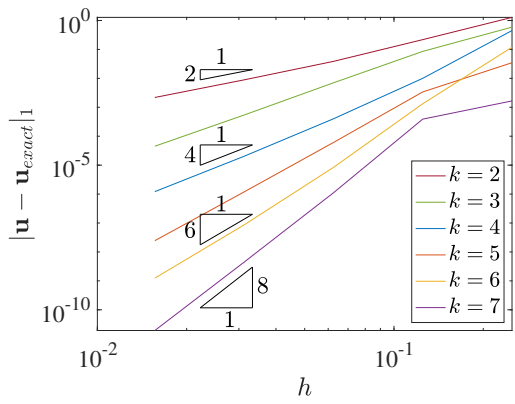}}\hfill
\subfloat[Pressure $H^1$ error]{\label{sfig:ns_rot_kovasznay_press_h1}\includegraphics[width=.45\textwidth]{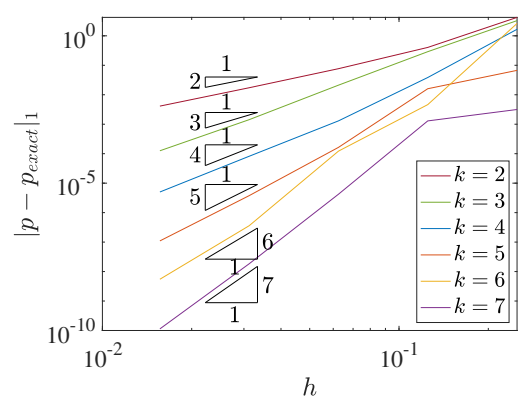}} \\
\caption{Errors in Kovasznay flow using rotational form}
\label{fig:ns_rot_kovasznay_2d_conv}
\end{figure}

Figure \ref{fig:ns_rot_kovasznay_2d_conv} shows the results obtained using the rotational form based method. We see that the method recovers the expected $L^2$ convergence rates of $k$ and $k+1$ for even and odd polynomial degrees. Interestingly, there seems to be some superconvergence occurring for odd polynomial degrees in the $H^1$ semi-norm of velocity, as they are converging almost like $k+1$ instead of the expected rate of $k$. We attribute this to pre-asymptotic behavior.

\subsubsection{2D Lid-Driven Cavity}

\begin{figure}
\centering
\subfloat[Re = 100 with 32$^2$ elements]{\label{sfig:Re100_32_sl}\includegraphics[width=.45\textwidth]{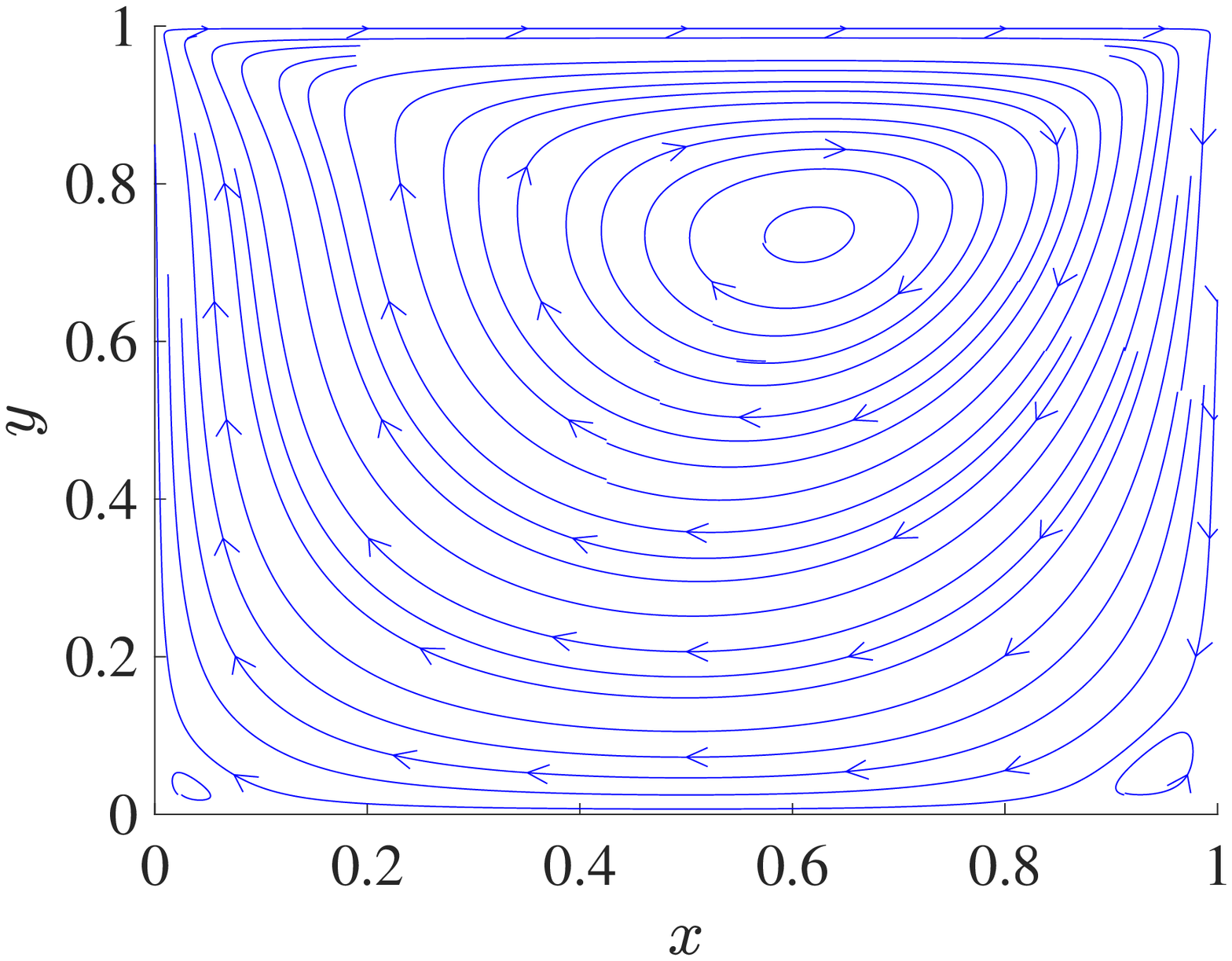}}\hfill
\subfloat[Re = 100 with 64$^2$ elements]{\label{sfig:Re100_64_sl}\includegraphics[width=.45\textwidth]{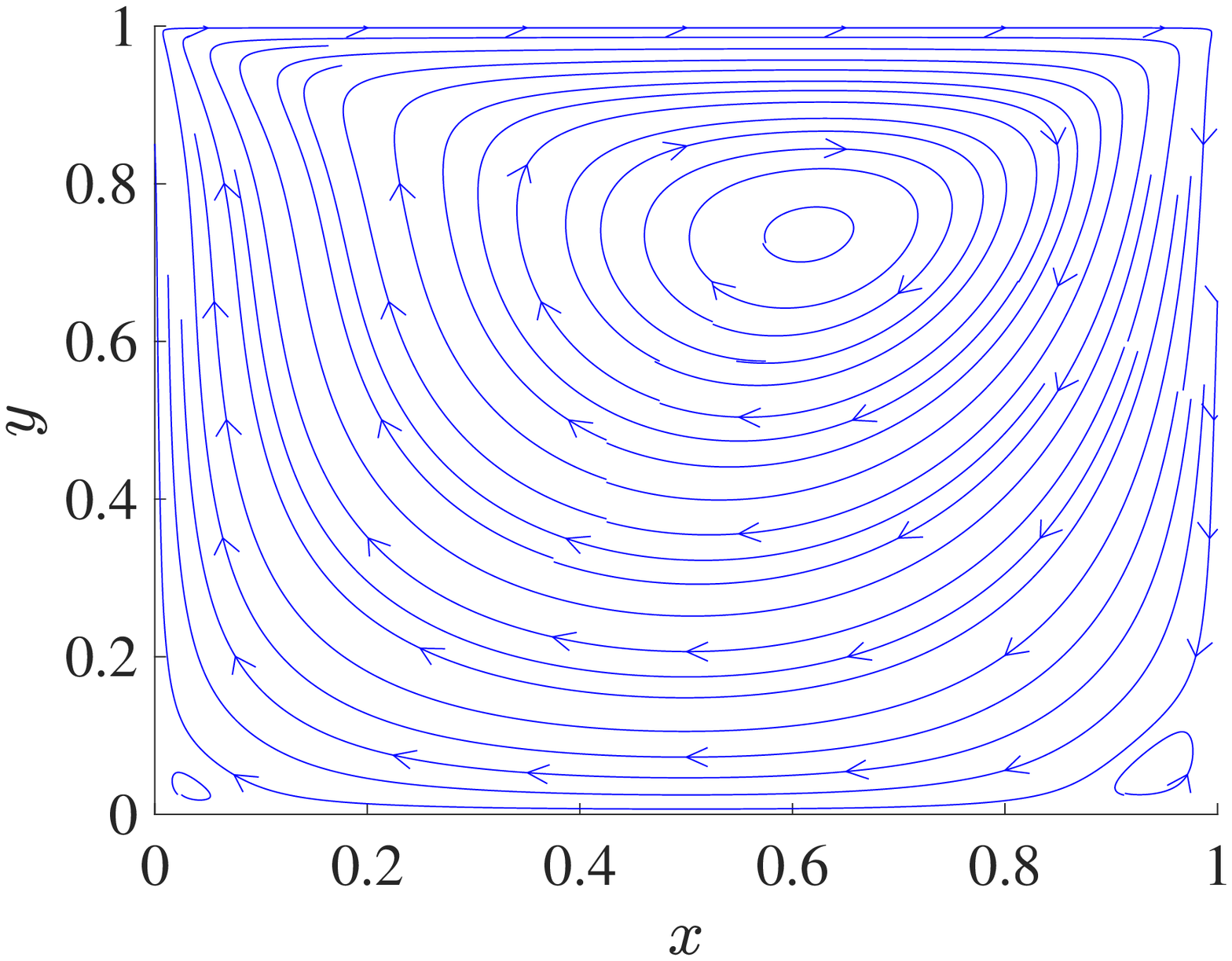}} \\
\subfloat[Re = 400 with 32$^2$ elements]{\label{sfig:Re400_32_sl}\includegraphics[width=.45\textwidth]{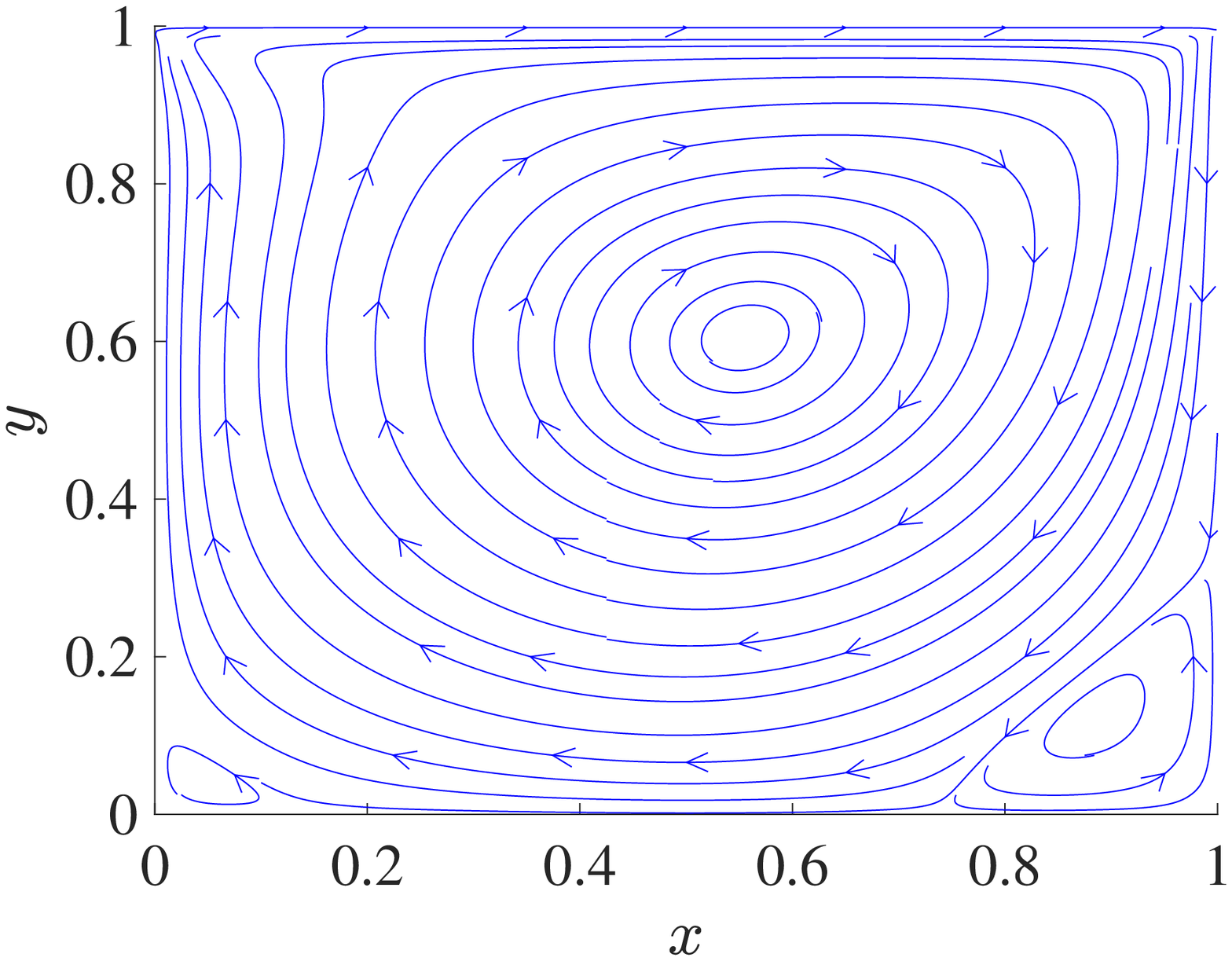}}\hfill
\subfloat[Re = 400 with 64$^2$ elements]{\label{sfig:Re400_64_sl}\includegraphics[width=.45\textwidth]{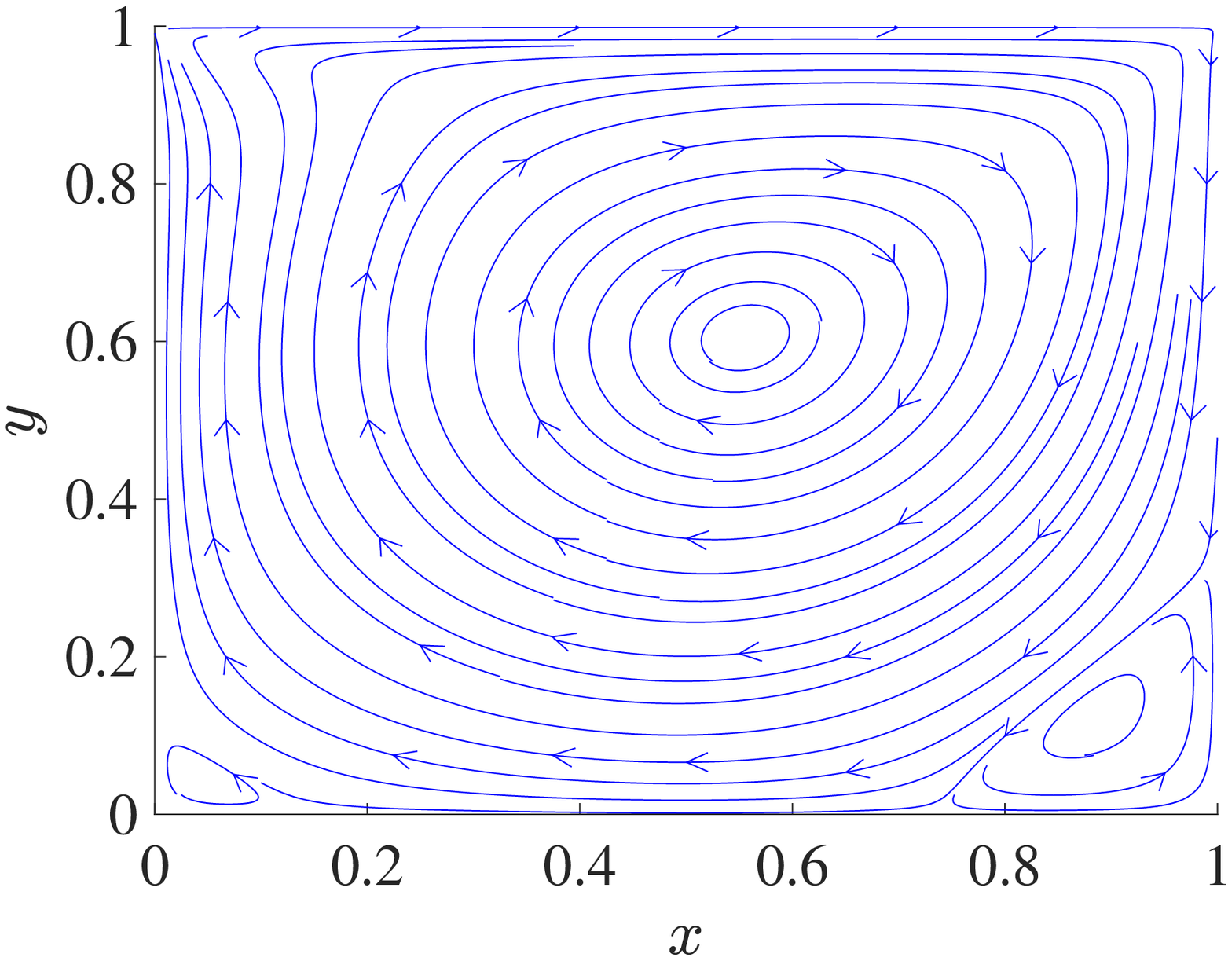}} \\
\subfloat[Re = 1000 with 32$^2$ elements]{\label{sfig:Re1000_32_sl}\includegraphics[width=.45\textwidth]{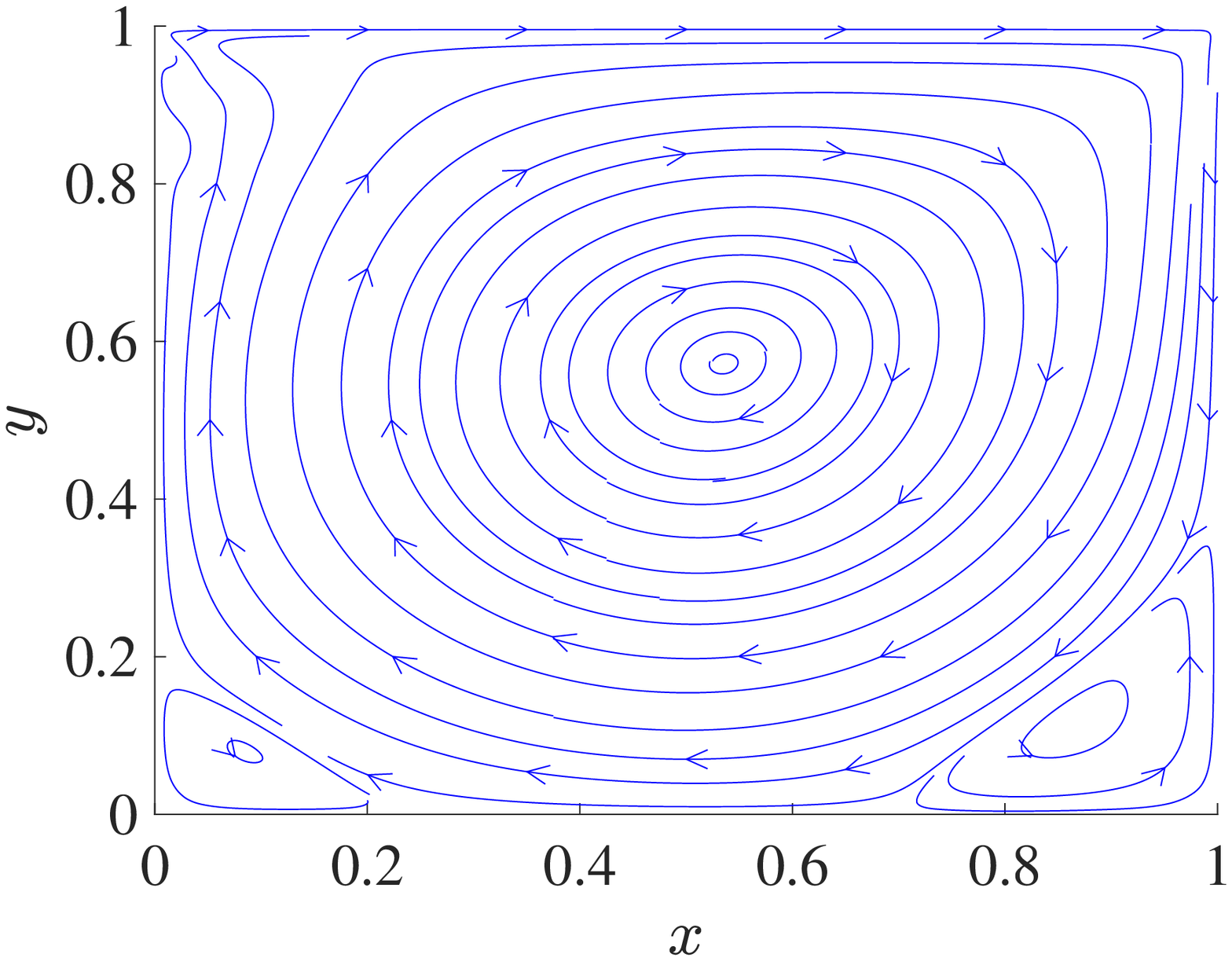}}\hfill
\subfloat[Re = 1000 with 64$^2$ elements]{\label{sfig:Re1000_64_sl}\includegraphics[width=.45\textwidth]{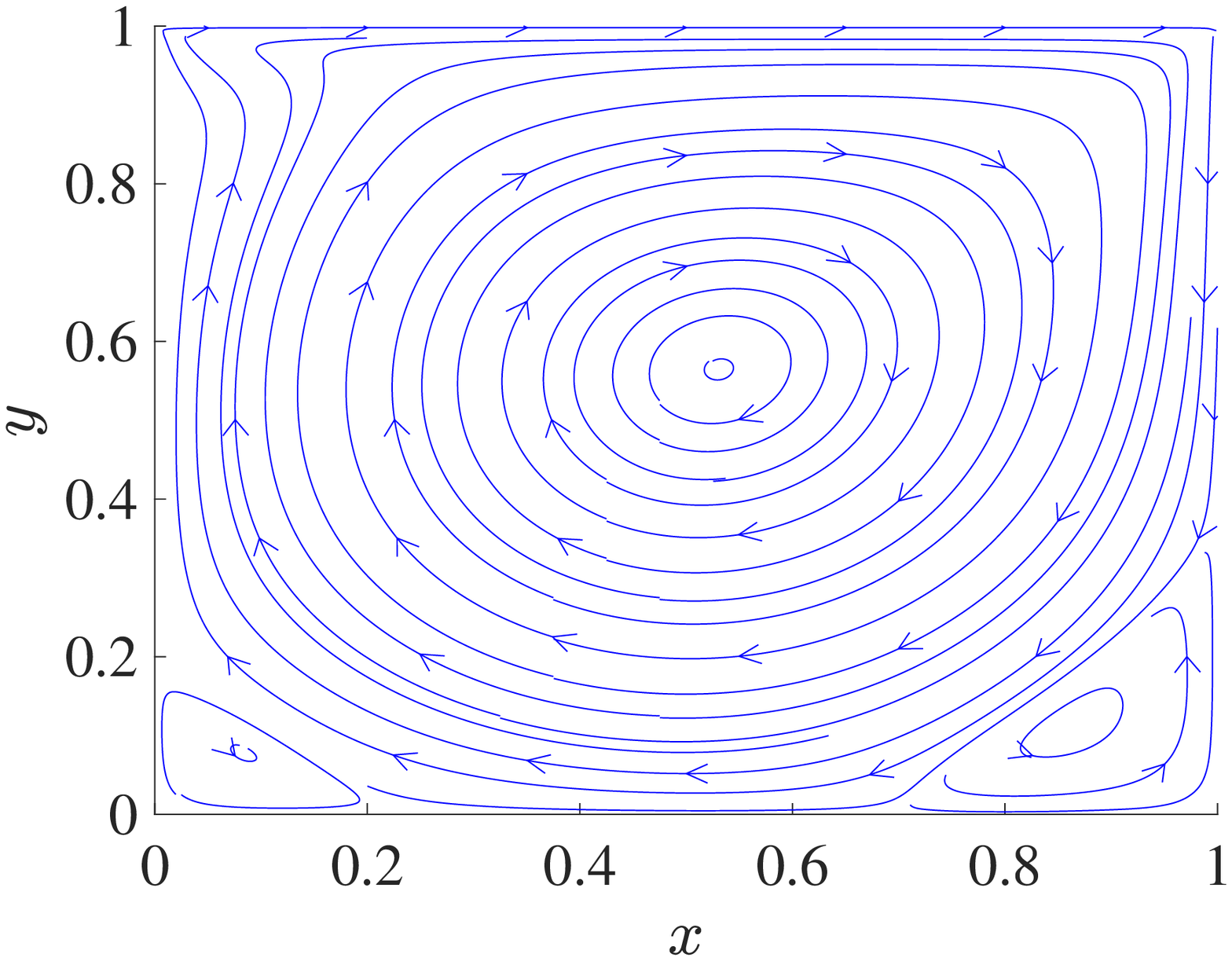}} \\
\caption{Streamlines of velocity for 2D lid-driven cavity flow using velocity-pressure form, $k$ = 5}
\label{fig:ns_2d_cavity_sl}
\end{figure}

We also apply our schemes to nonlinear lid-driven cavity flows at Reynolds numbers (defined by the reciprocal of the viscosity $\nu$ as the problem is posed on the unit square) of 100, 400, and 1000, and compare the results with those obtained in \cite{ghia_cavity}. Figure \ref{fig:ns_2d_cavity_sl} shows the flow streamlines obtained with a polynomial degree $k = 5$ and 32 (left) or 64 (right) elements in each direction of a stretched mesh with knots defined by 

 \begin{equation}
     \xi_i = \frac{1}{2} \left( 1 + \frac{\tanh(3ih - 2)}{\tanh(2)} \right) \quad \forall \xi_i \in \Xi.
 \end{equation}

\begin{figure}
\centering
\subfloat[Re = 100 with 32 elements]{\label{sfig:Re100_32p8_sl}\includegraphics[width=.45\textwidth]{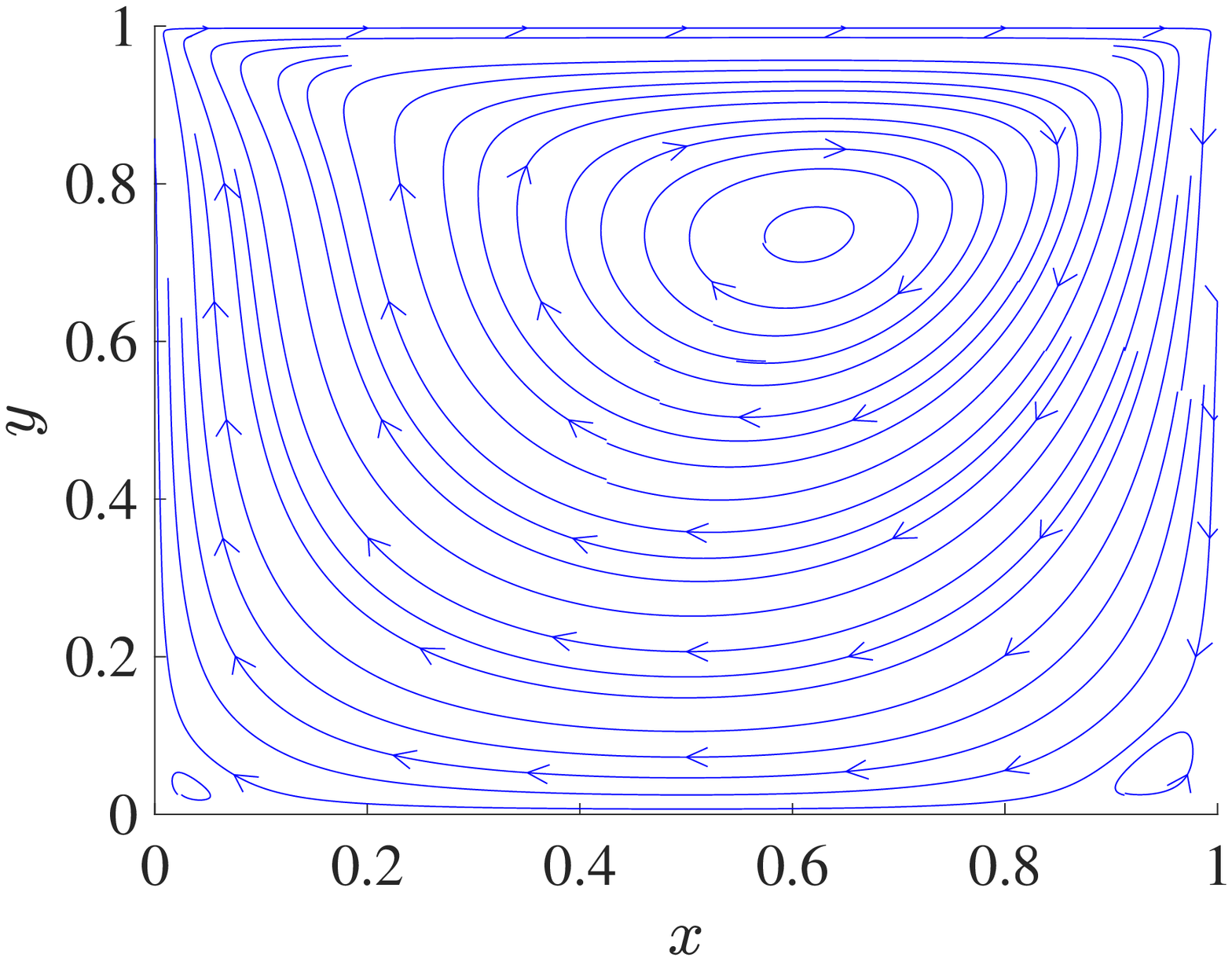}}\hfill
\subfloat[Re = 100 with 64$^2$ elements]{\label{sfig:Re100_64p8_sl}\includegraphics[width=.45\textwidth]{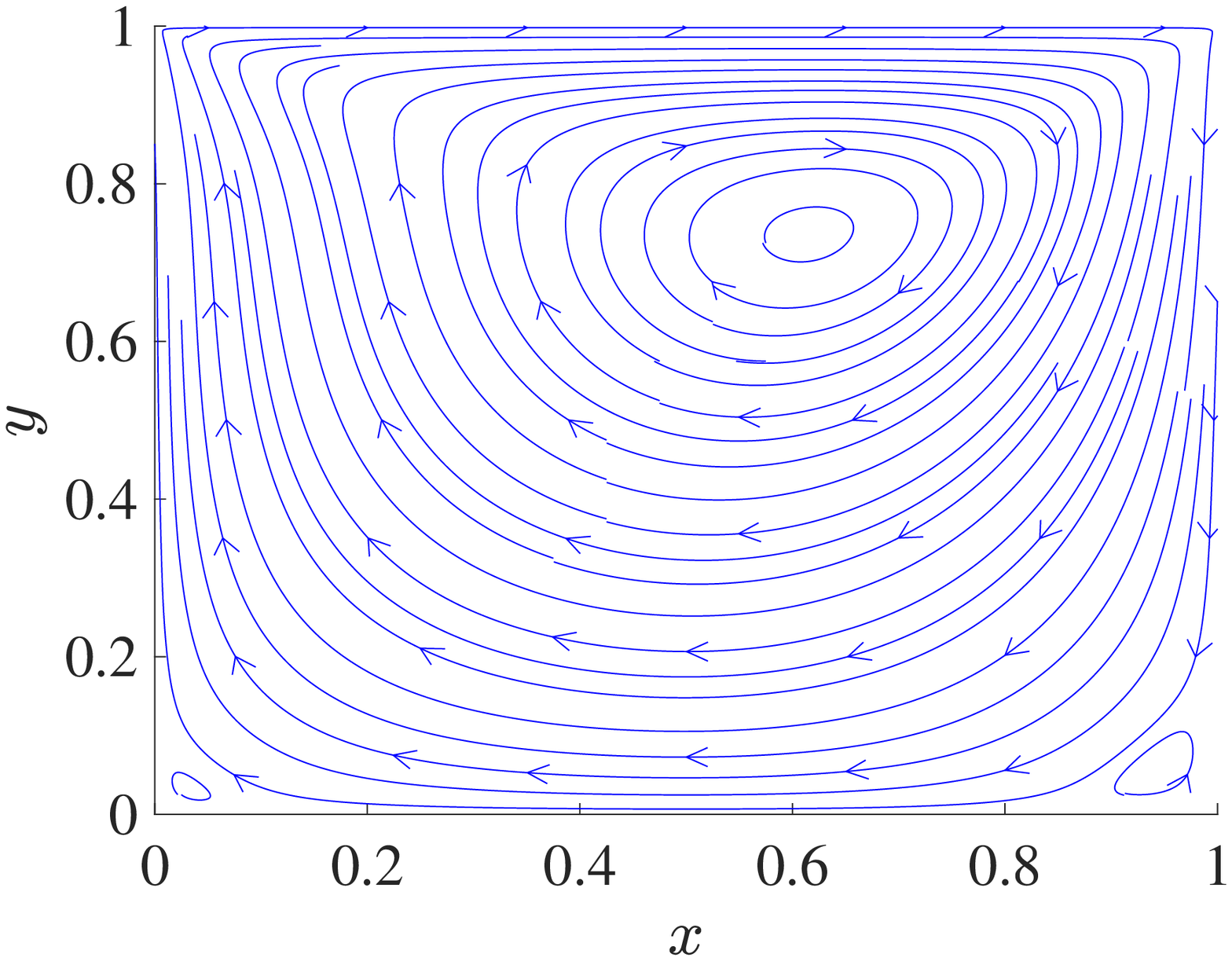}} \\
\subfloat[Re = 400 with 32$^2$ elements]{\label{sfig:Re400_32p8_sl}\includegraphics[width=.45\textwidth]{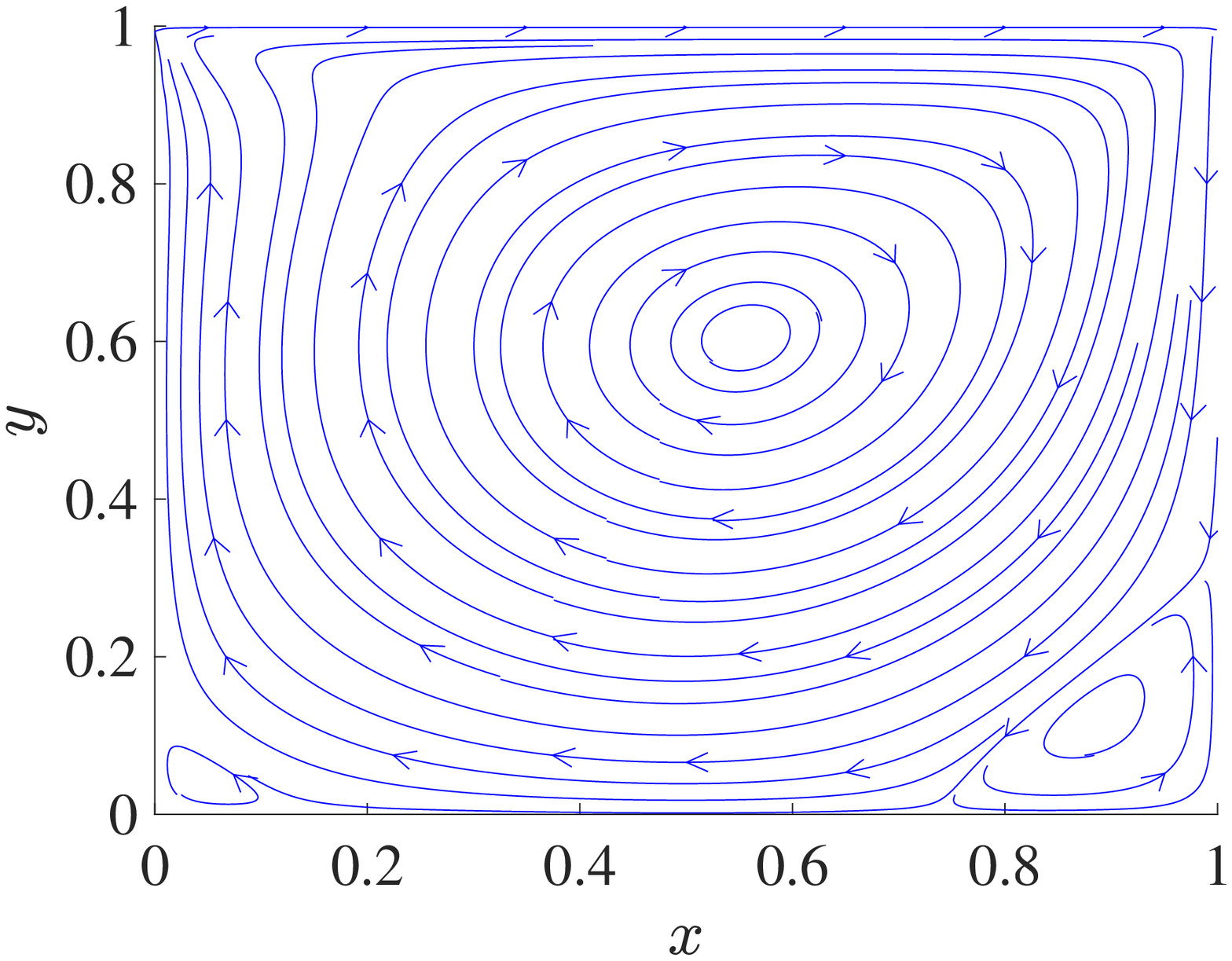}}\hfill
\subfloat[Re = 400 with 64$^2$ elements]{\label{sfig:Re400_64p8_sl}\includegraphics[width=.45\textwidth]{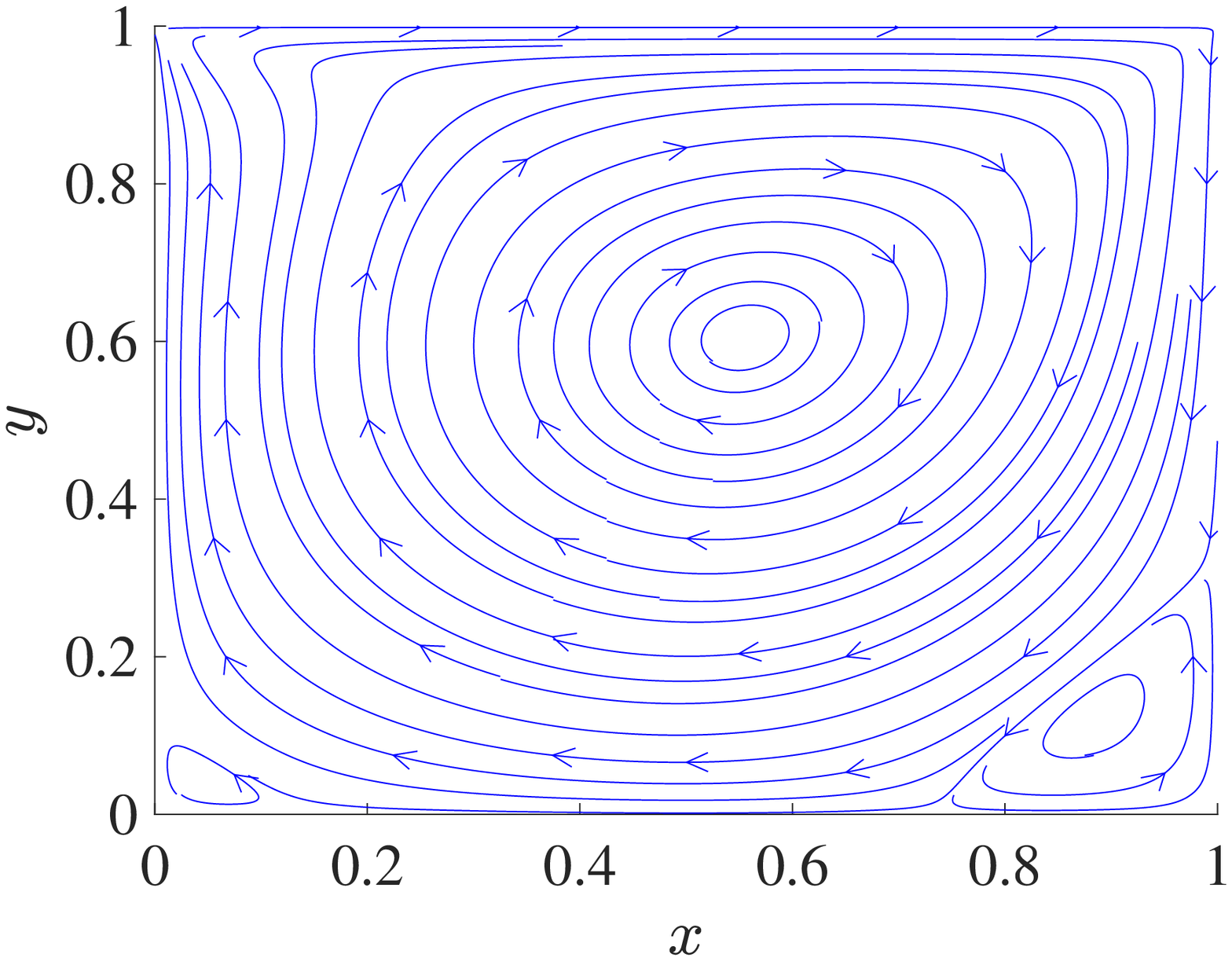}} \\
\subfloat[Re = 1000 with 32$^2$ elements]{\label{sfig:Re1000_32p8_sl}\includegraphics[width=.45\textwidth]{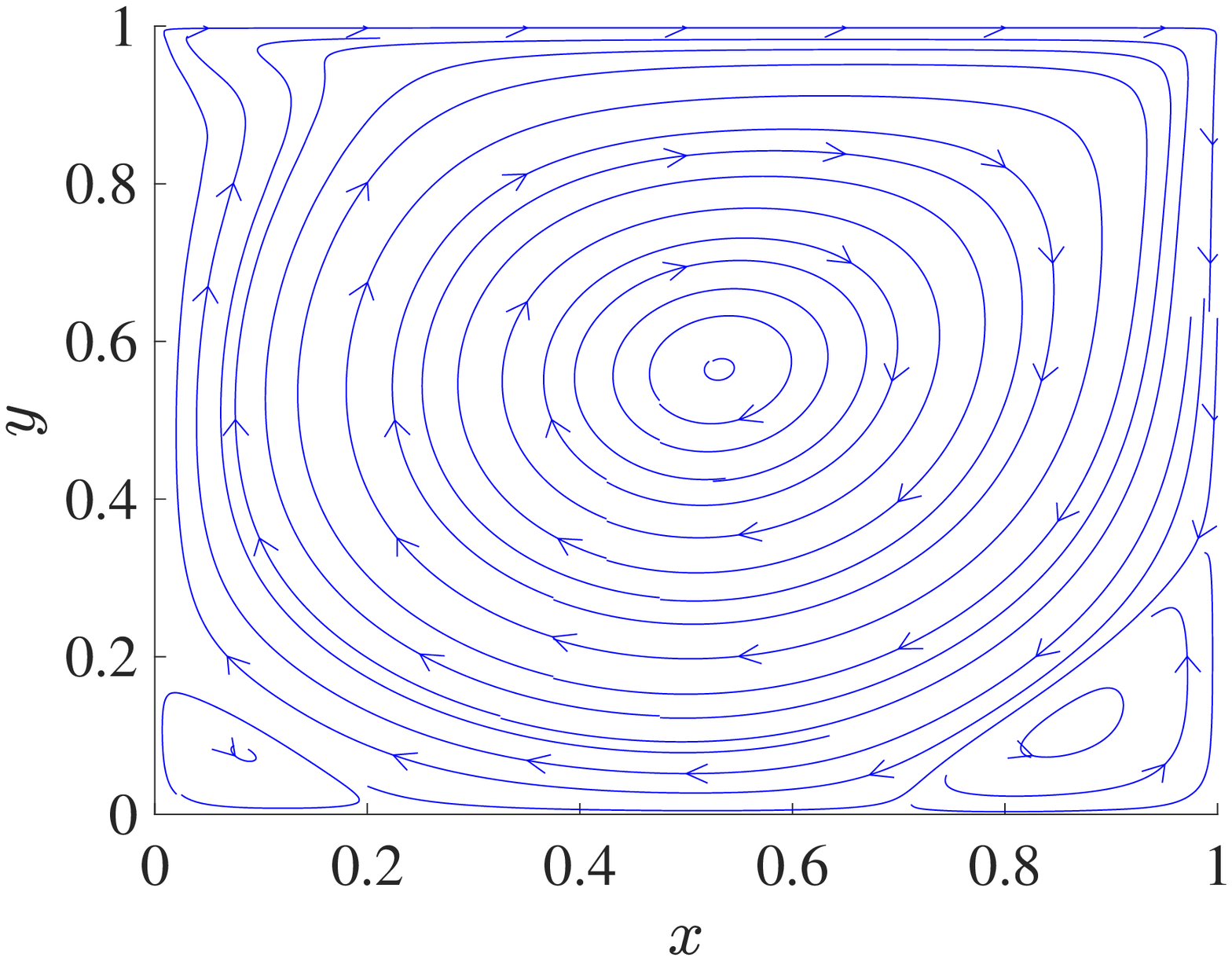}}\hfill
\subfloat[Re = 1000 with 64$^2$ elements]{\label{sfig:Re1000_64p8_sl}\includegraphics[width=.45\textwidth]{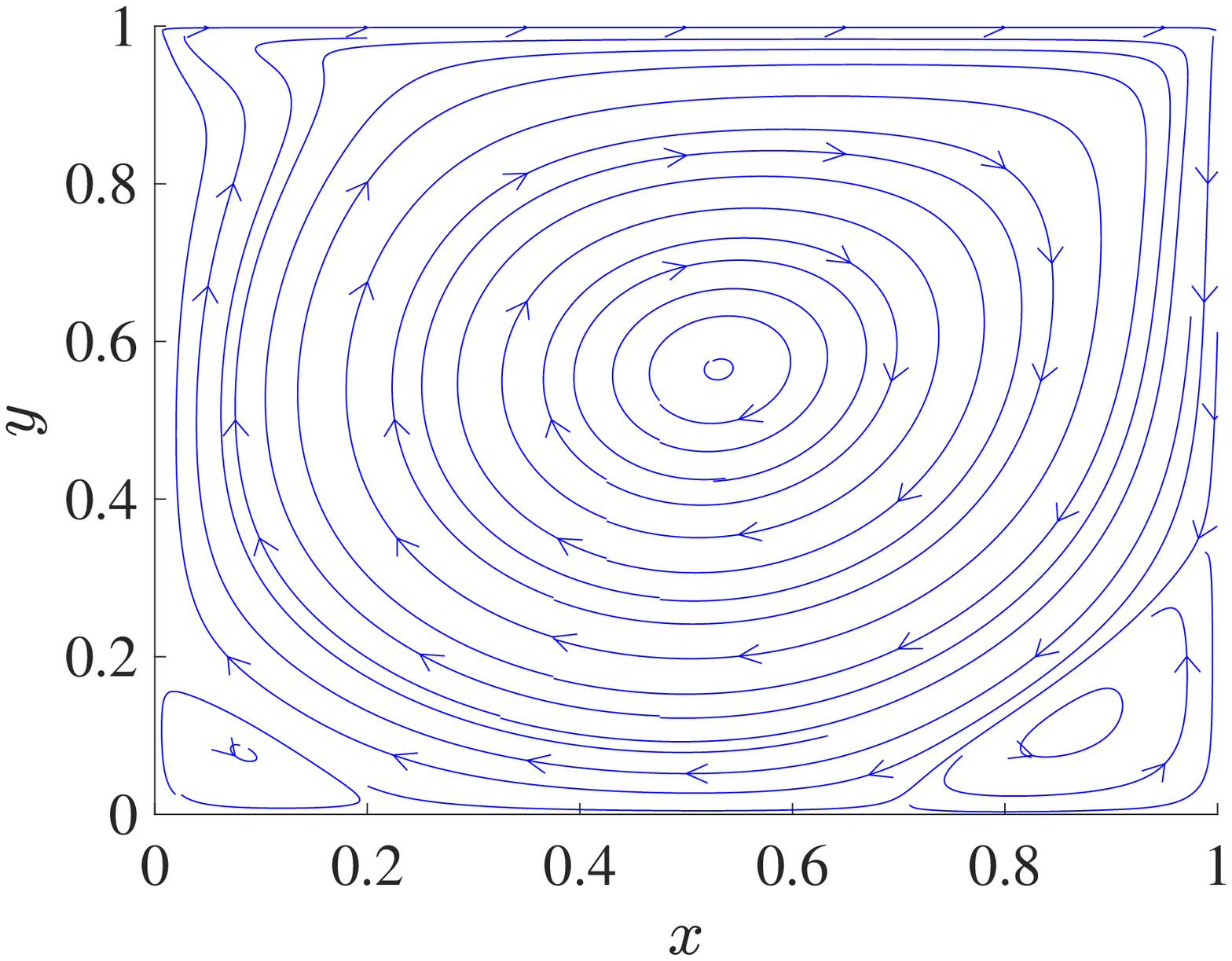}} \\
\caption{Streamlines of velocity for 2D lid-driven cavity flow using velocity-pressure form, $k$ = 8}
\label{fig:ns_2d_cavity_sl_p8}
\end{figure}

\noindent Clearly the scheme recovers the expected behaviors of the flow. As the Reynolds number increases, we see the corner vortices increase in size and the flow become non-symmetric. Note that at a Reynolds number of 1000 the solution on the coarse mesh struggles to accommodate the boundary condition at the top left corner of the domain. This is not unexpected, as the Dirichlet boundary data is extremely sharp in that region. When the mesh is refined this effect disappears as expected. Figures \ref{fig:ns_2d_cavity_sl_p8} shows the same results but with the polynomial degree elevated to $k = 8$. Note that elevating the degree also resolves the flow in the top left corner even on the coarse mesh.

\begin{figure}
\centering
\subfloat[Re = 100 with 32$^2$ elements]{\label{sfig:Re100_32}\includegraphics[width=.45\textwidth]{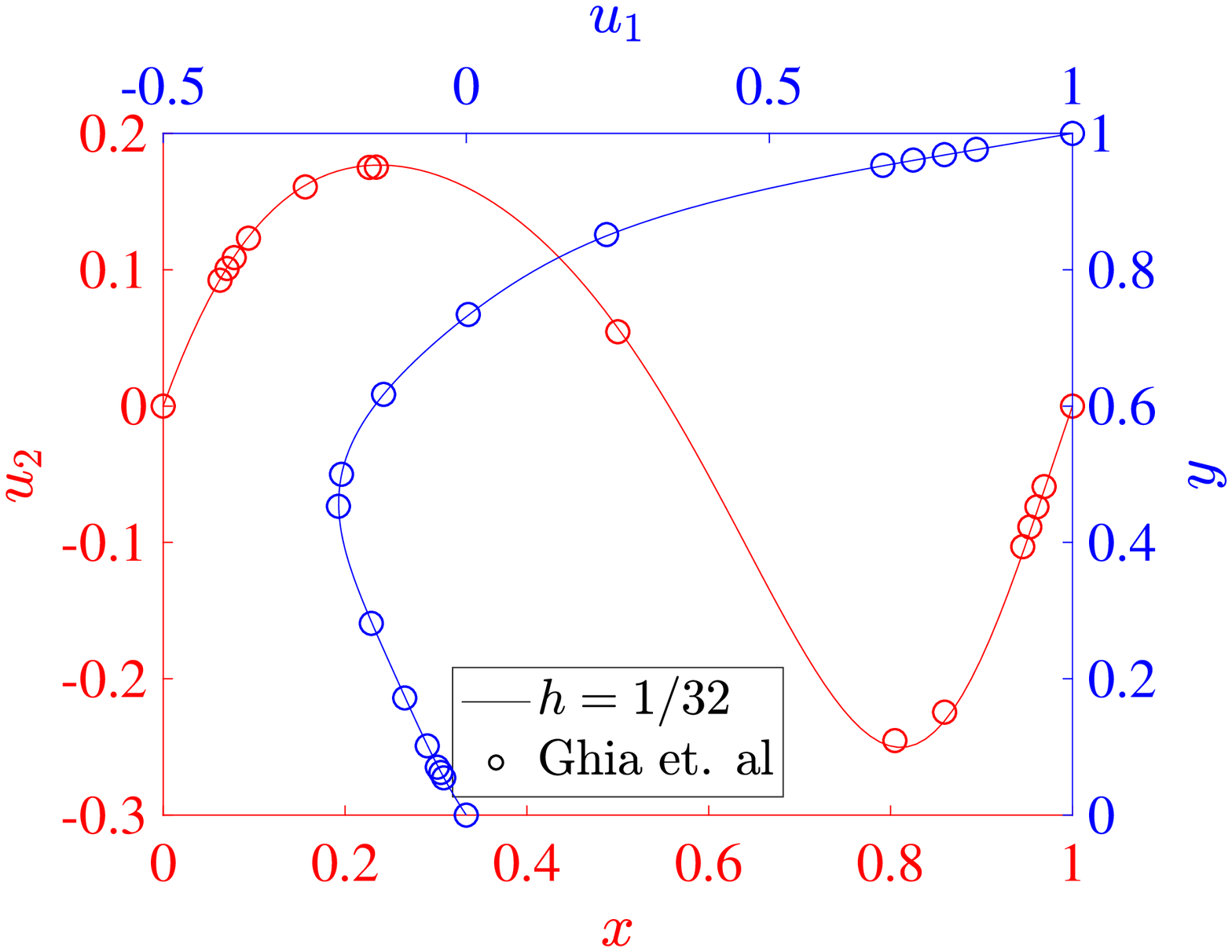}}\hfill
\subfloat[Re = 100 with 64$^2$ elements]{\label{sfig:Re100_64}\includegraphics[width=.45\textwidth]{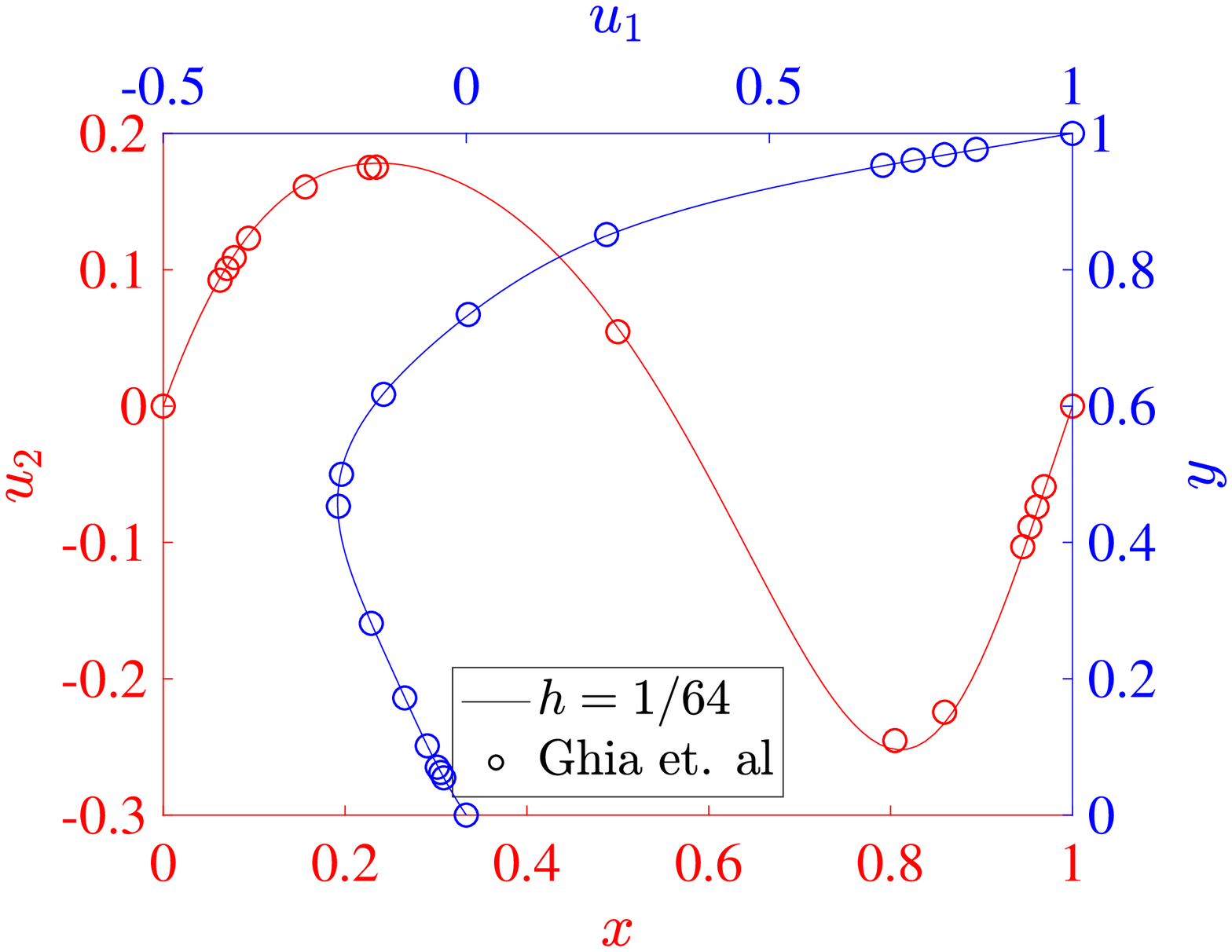}} \\
\subfloat[Re = 400 with 32$^2$ elements]{\label{sfig:Re400_32}\includegraphics[width=.45\textwidth]{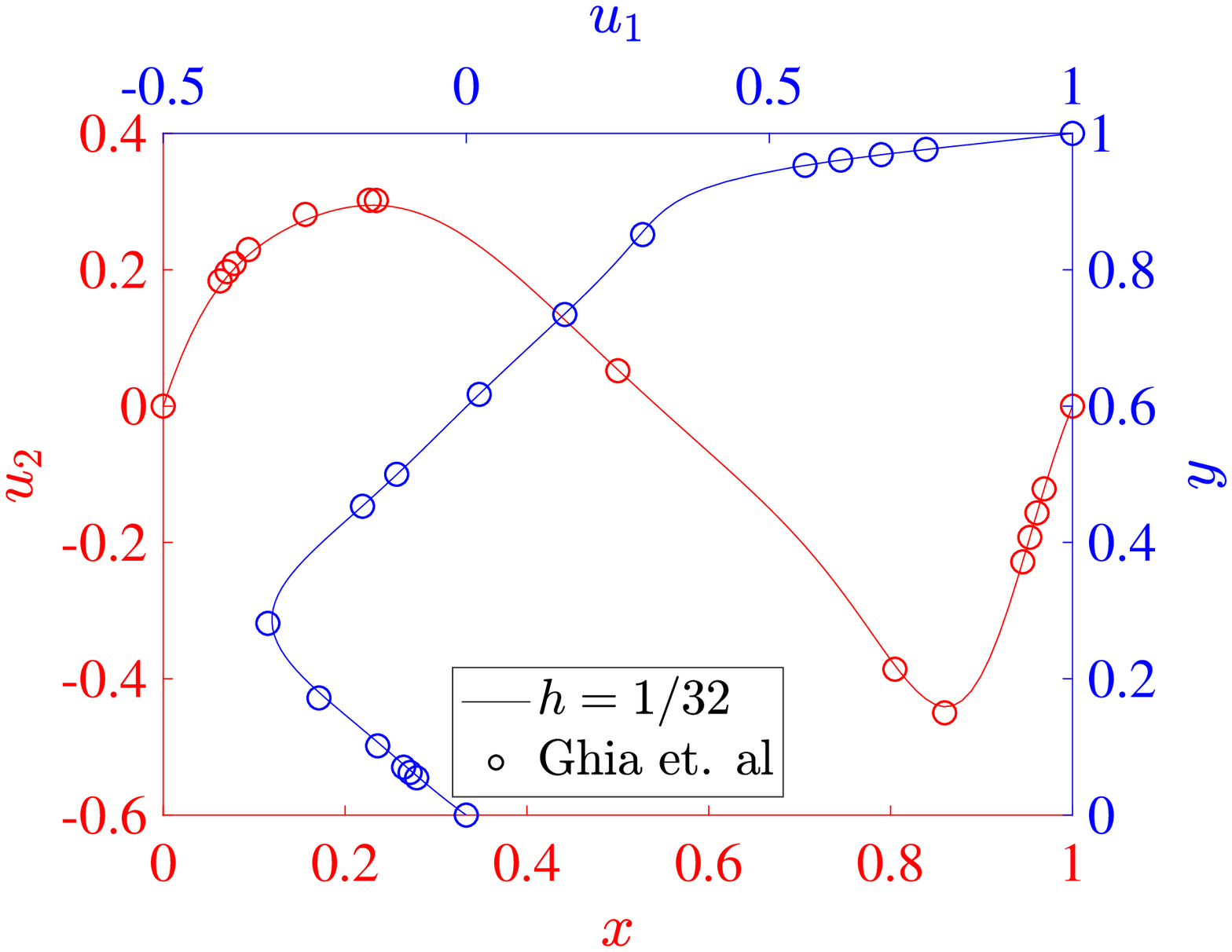}}\hfill
\subfloat[Re = 400 with 64$^2$ elements]{\label{sfig:Re400_64}\includegraphics[width=.45\textwidth]{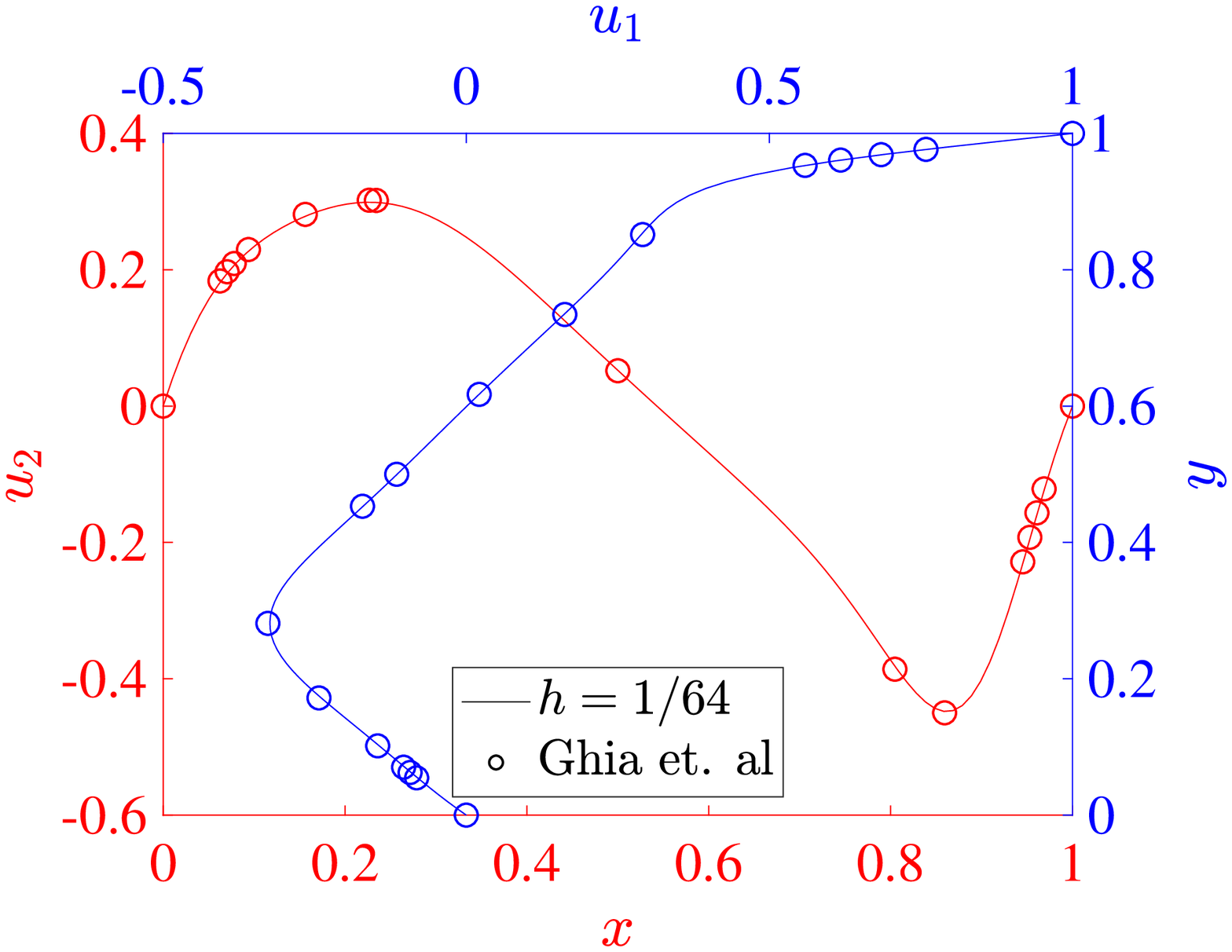}} \\
\subfloat[Re = 1000 with 32$^2$ elements]{\label{sfig:Re1000_32}\includegraphics[width=.45\textwidth]{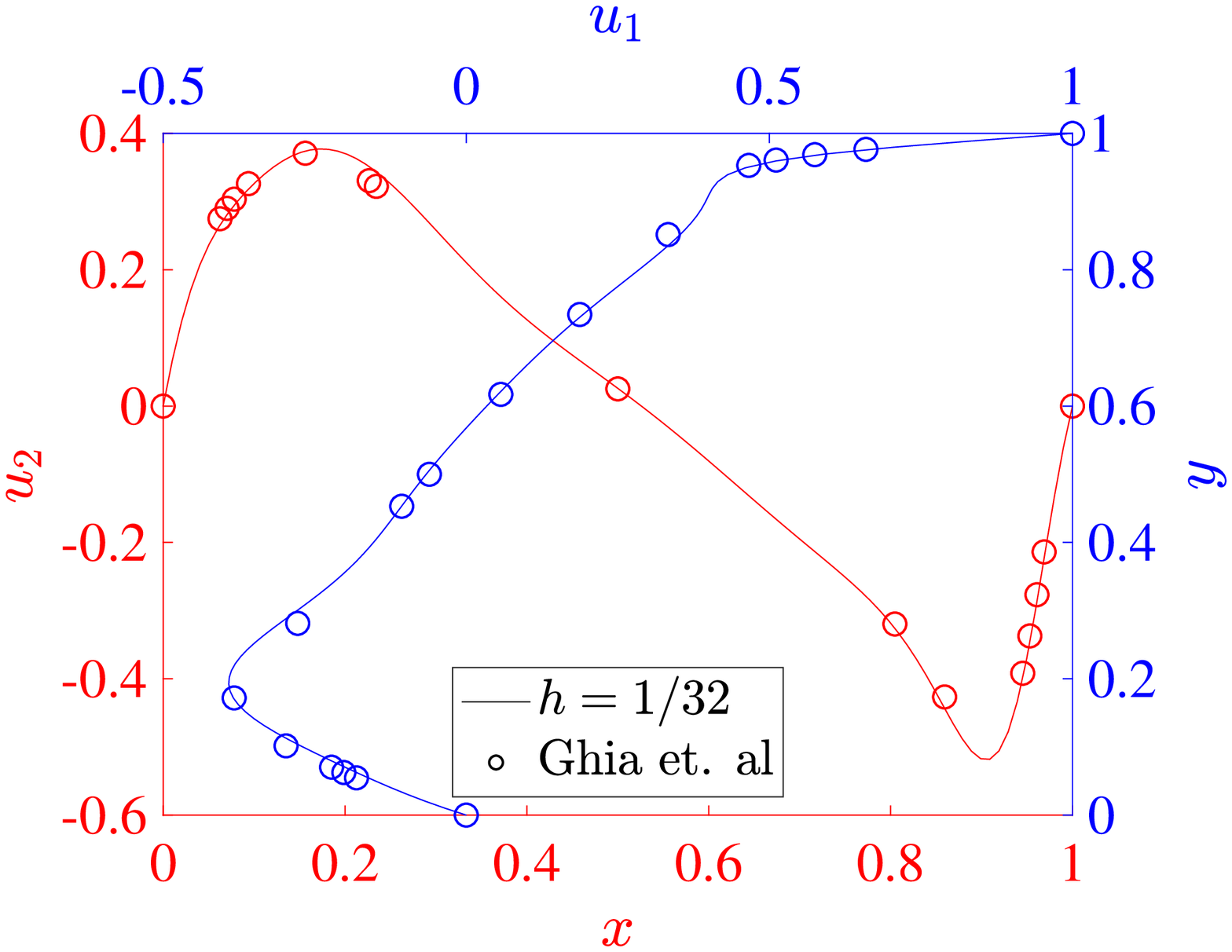}}\hfill
\subfloat[Re = 1000 with 64$^2$ elements]{\label{sfig:Re1000_64}\includegraphics[width=.45\textwidth]{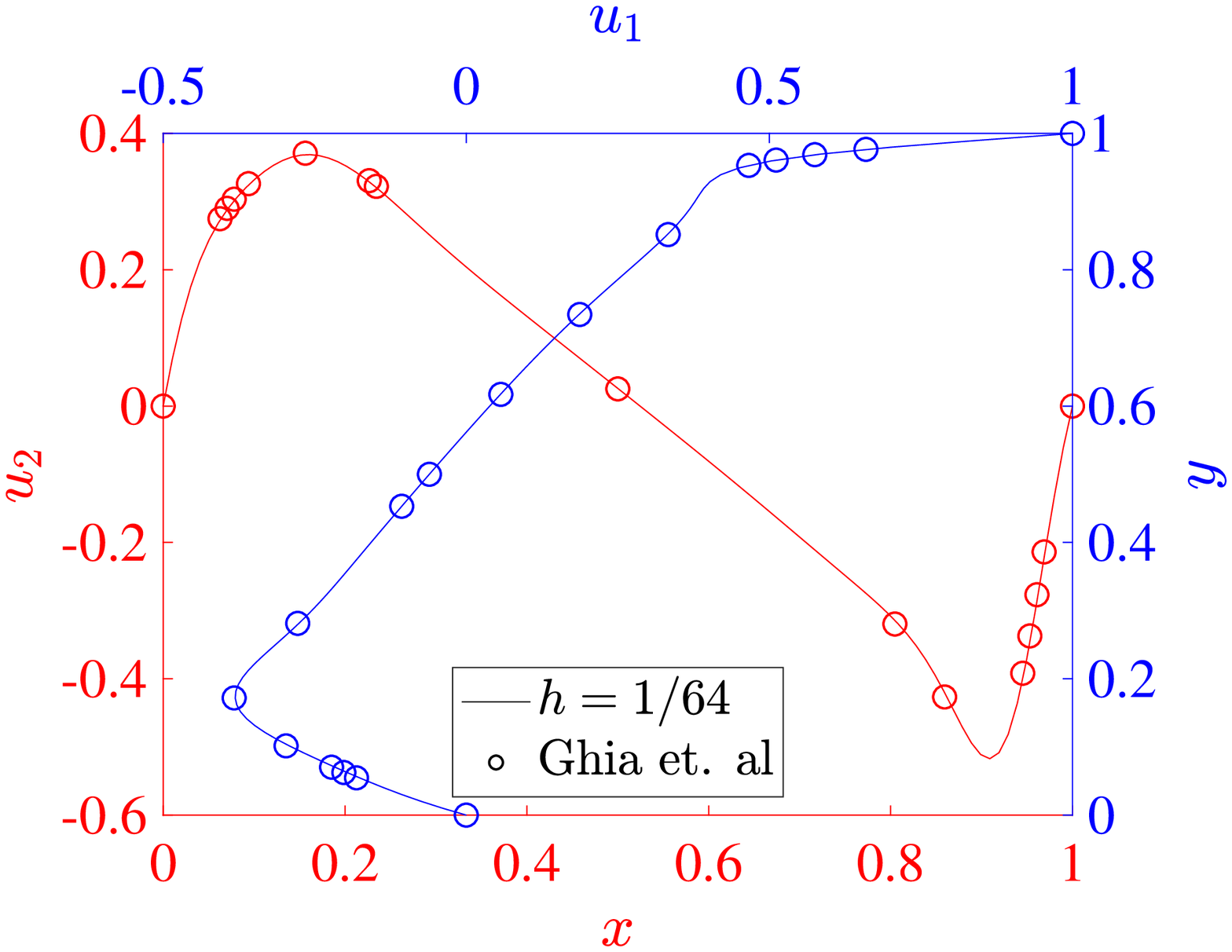}} \\
\caption{Centerline velocity profiles for 2D lid-driven cavity flow using velocity-pressure form, $k$ = 5}
\label{fig:ns_2d_cavity}
\end{figure}

\begin{figure}
\centering
\subfloat[Re = 100 with 32$^2$ elements]{\label{sfig:Re100_32p8}\includegraphics[width=.45\textwidth]{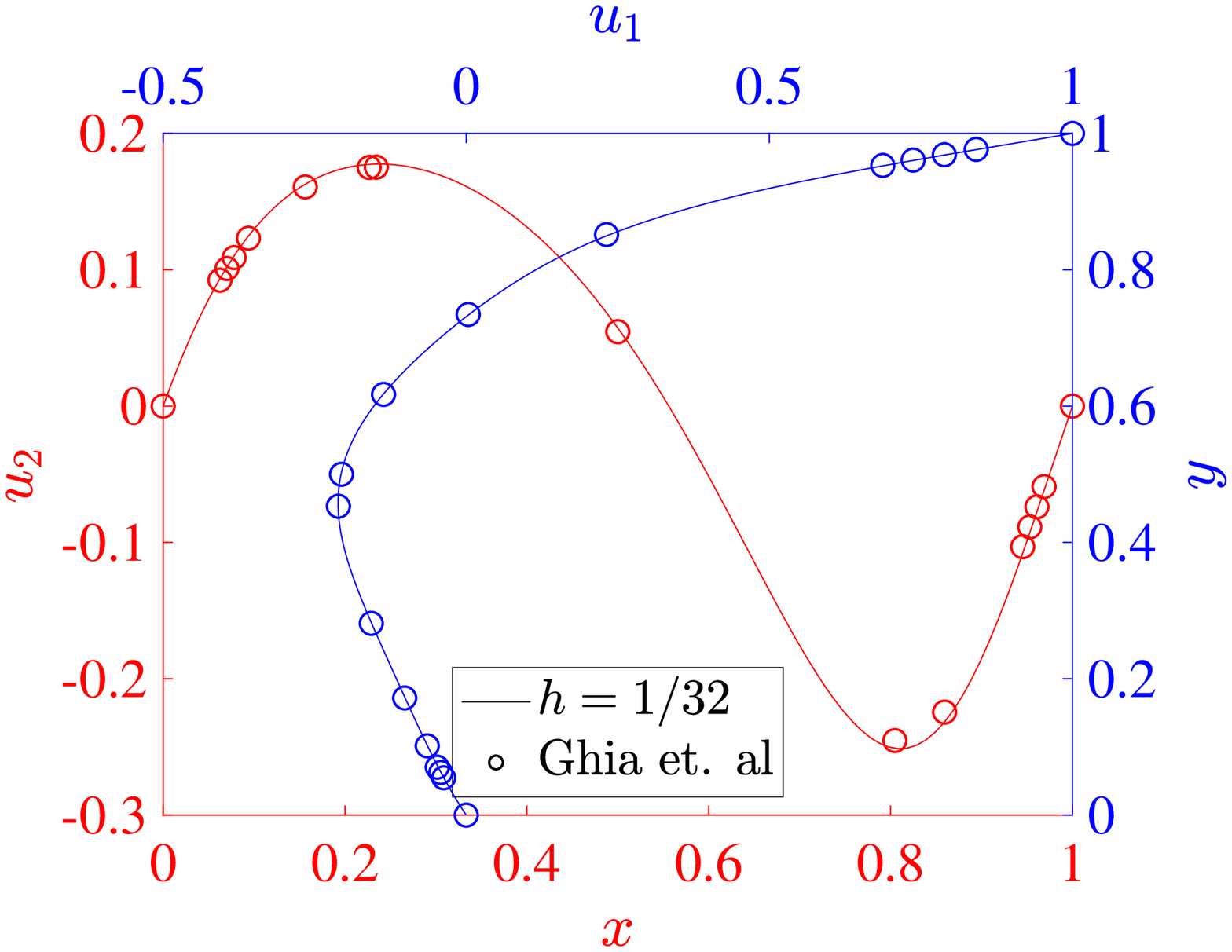}}\hfill
\subfloat[Re = 100 with 64$^2$ elements]{\label{sfig:Re100_64p8}\includegraphics[width=.45\textwidth]{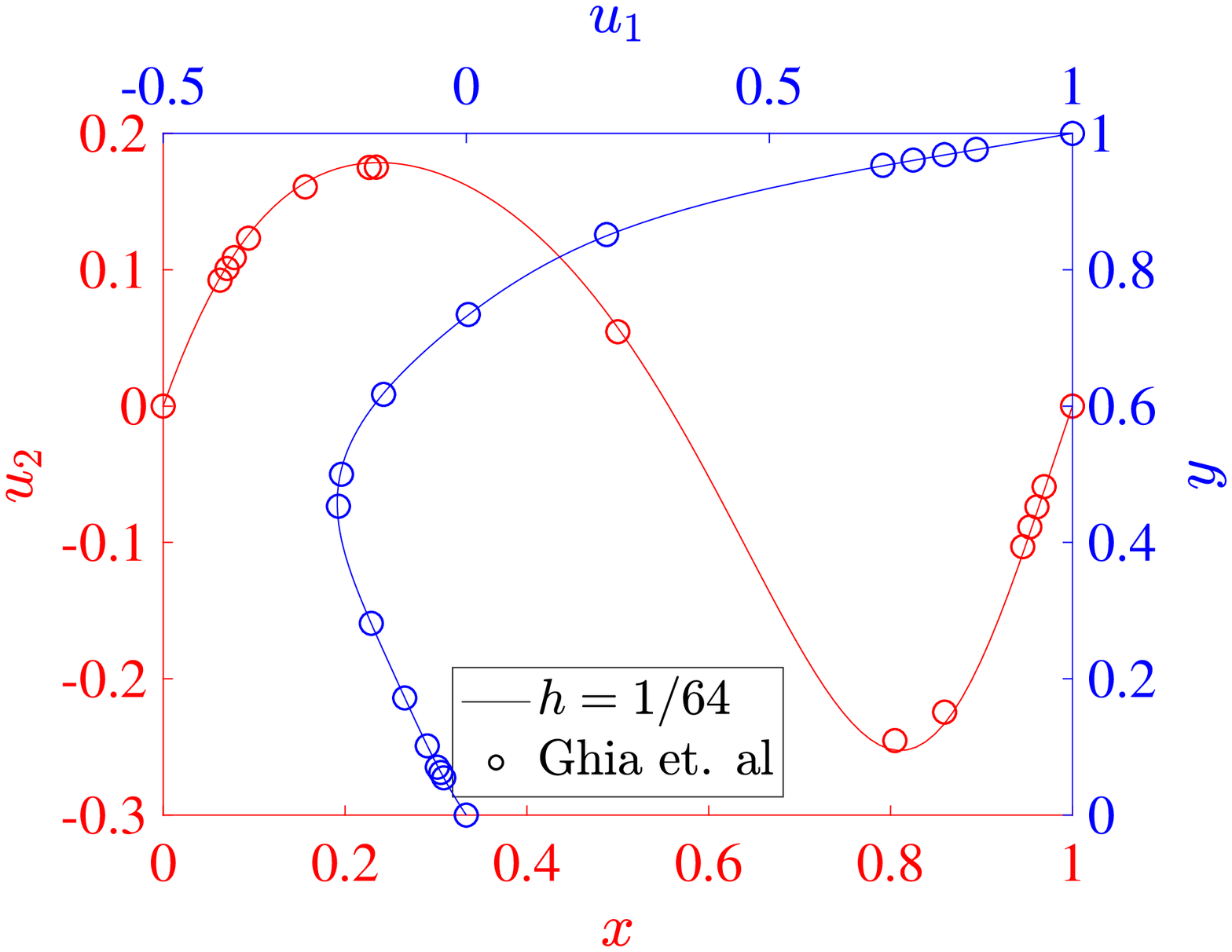}} \\
\subfloat[Re = 400 with 32$^2$ elements]{\label{sfig:Re400_32p8}\includegraphics[width=.45\textwidth]{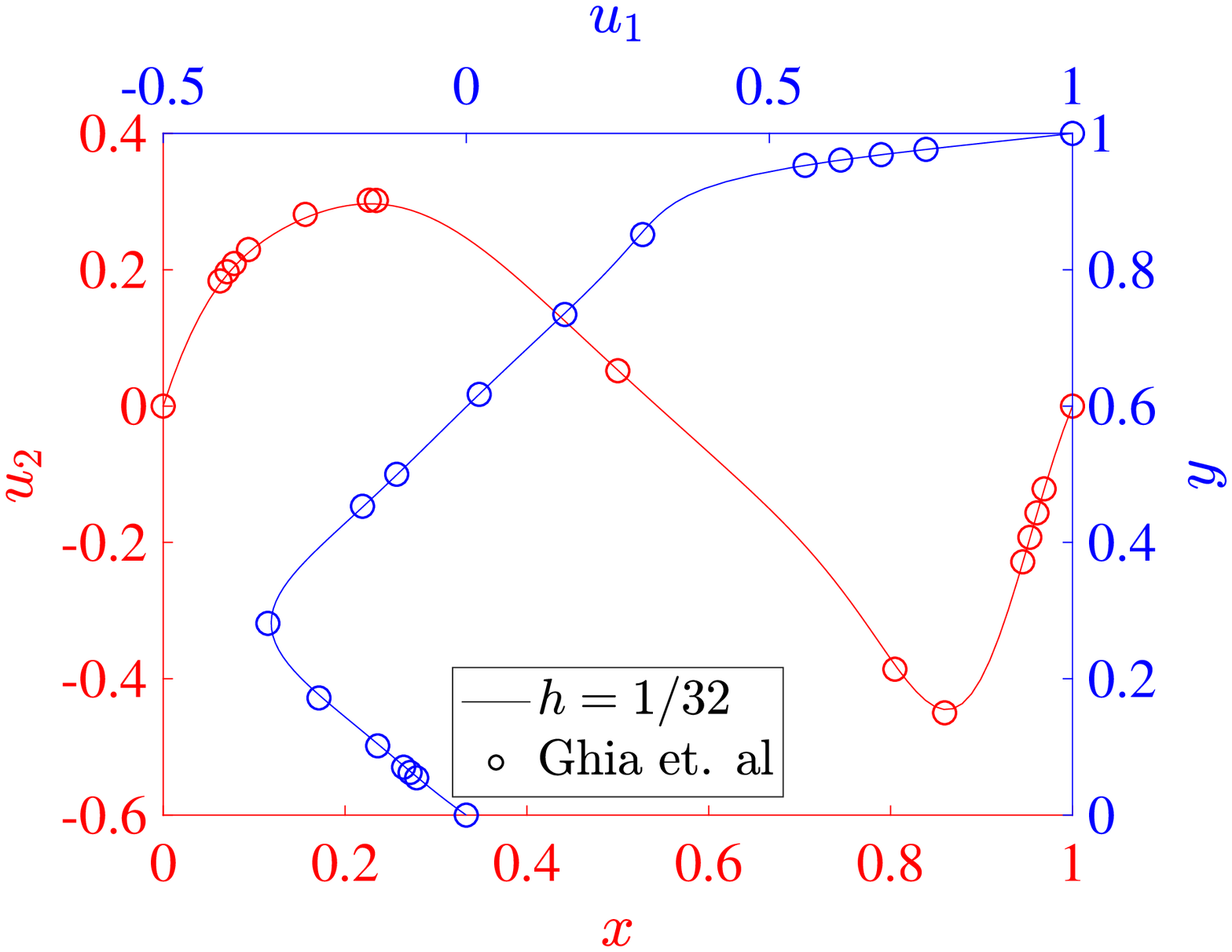}}\hfill
\subfloat[Re = 400 with 64$^2$ elements]{\label{sfig:Re400_64p8}\includegraphics[width=.45\textwidth]{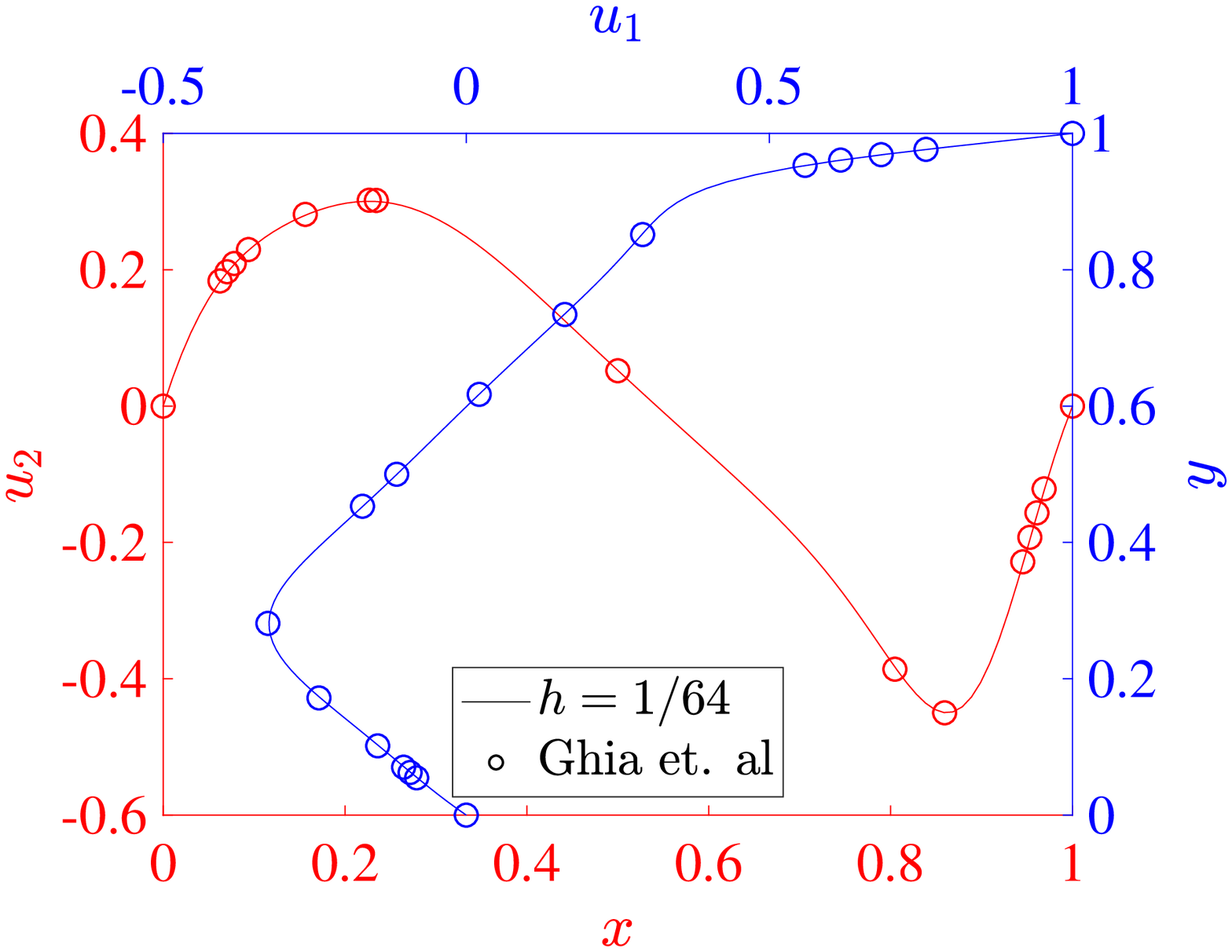}} \\
\subfloat[Re = 1000 with 32$^2$ elements]{\label{sfig:Re1000_32p8}\includegraphics[width=.45\textwidth]{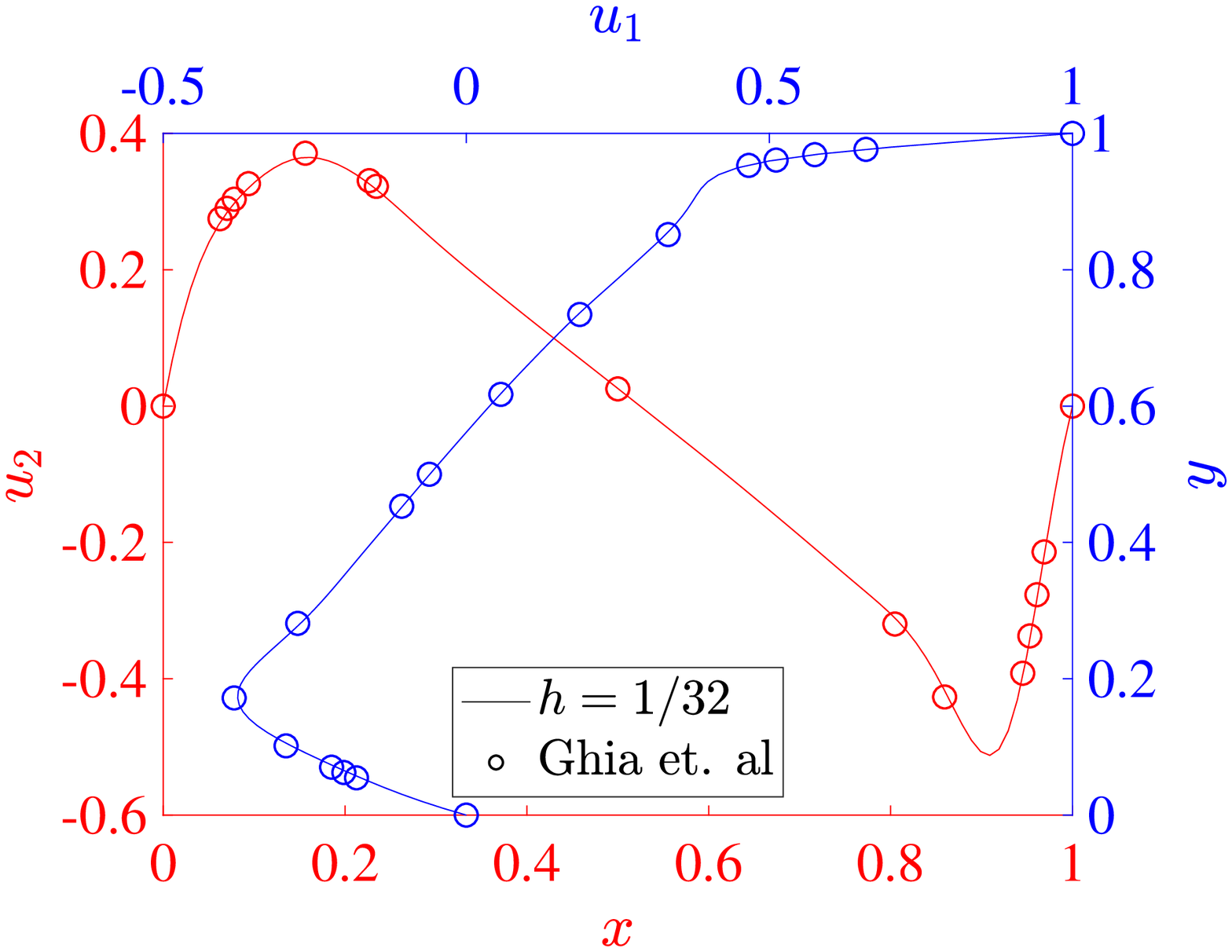}}\hfill
\subfloat[Re = 1000 with 64$^2$ elements]{\label{sfig:Re1000_64p8}\includegraphics[width=.45\textwidth]{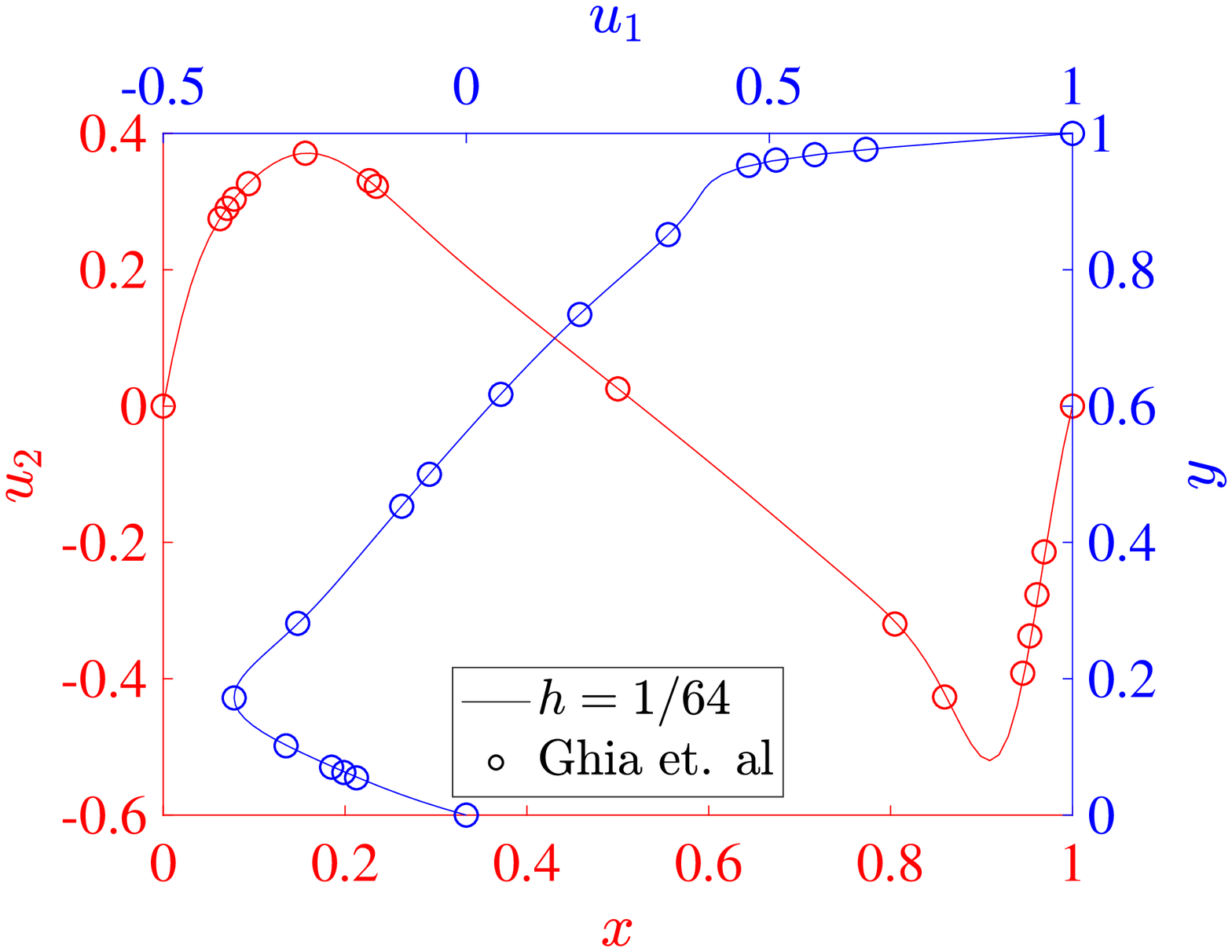}} \\
\caption{Centerline velocity profiles for 2D lid-driven cavity flow using velocity-pressure form, $k$ = 8}
\label{fig:ns_2d_cavity_p8}
\end{figure}

Figure \ref{fig:ns_2d_cavity} shows the horizontal and vertical velocities along the vertical and horizontal centerlines, respectively, for each simulation with $k = 5$. At a Reynolds number of 100 the results match well with the reference results of Ghia at both resolutions. At a Reynolds number of 400 the coarse mesh results produce velocity fields that slightly under-predict the magnitudes of the maximum and minimum velocities in both directions. The results from the refined mesh seem to mitigate these discrepancies. Finally, at the highest Reynolds number the coarse mesh results again possess small discrepancies from the results of Ghia, but mesh refinement aids. The same centerline results obtained with $k = 8$ are shown in Figure \ref{fig:ns_2d_cavity_p8}. Degree elevation seems to produce results in closer agreement to Ghia in all cases.

\begin{figure}
\centering
\subfloat[Re = 100 with 32$^2$ elements]{\label{sfig:Re100_32_sl_rot}\includegraphics[width=.45\textwidth]{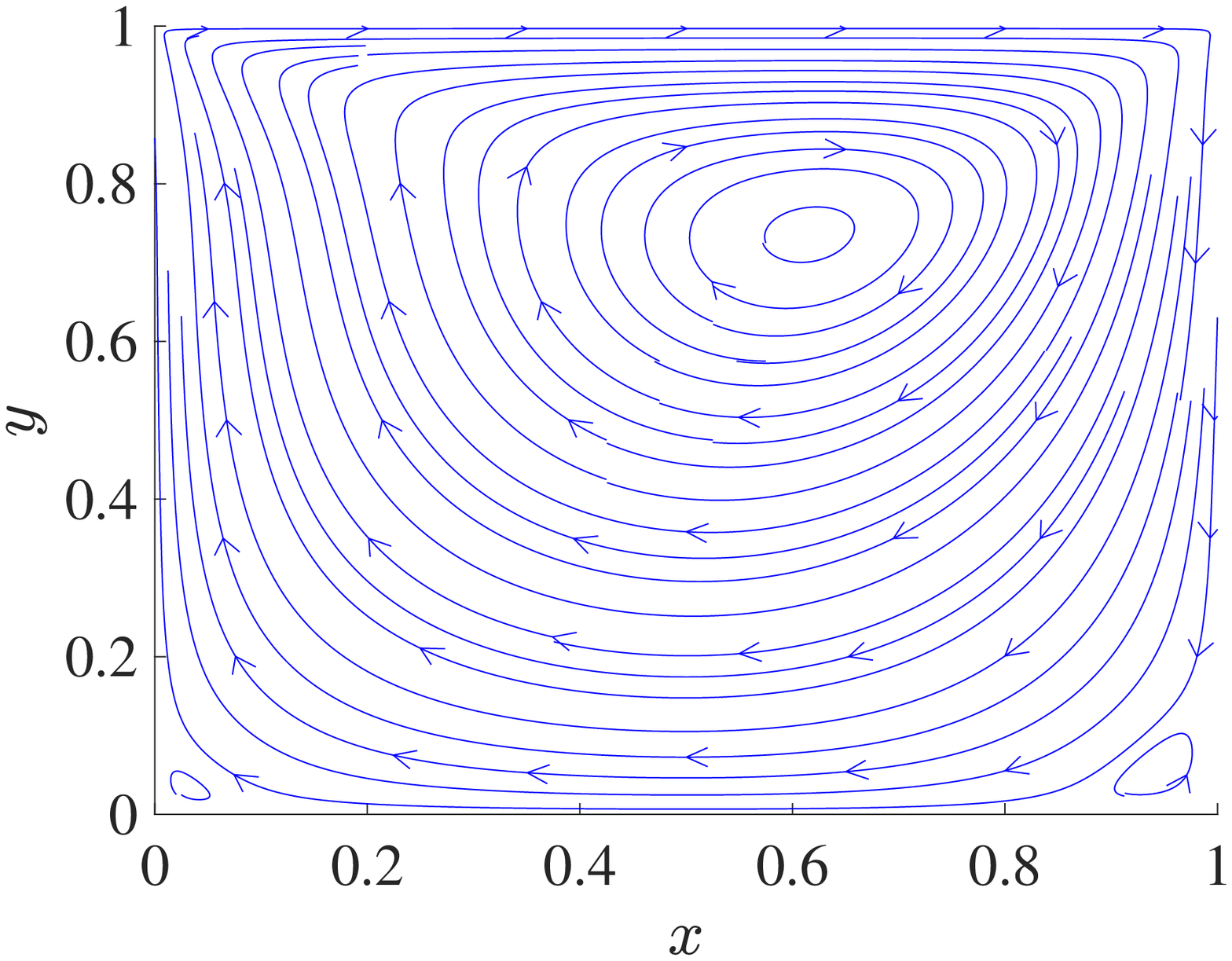}}\hfill
\subfloat[Re = 100 with 64$^2$ elements]{\label{sfig:Re100_64_sl_rot}\includegraphics[width=.45\textwidth]{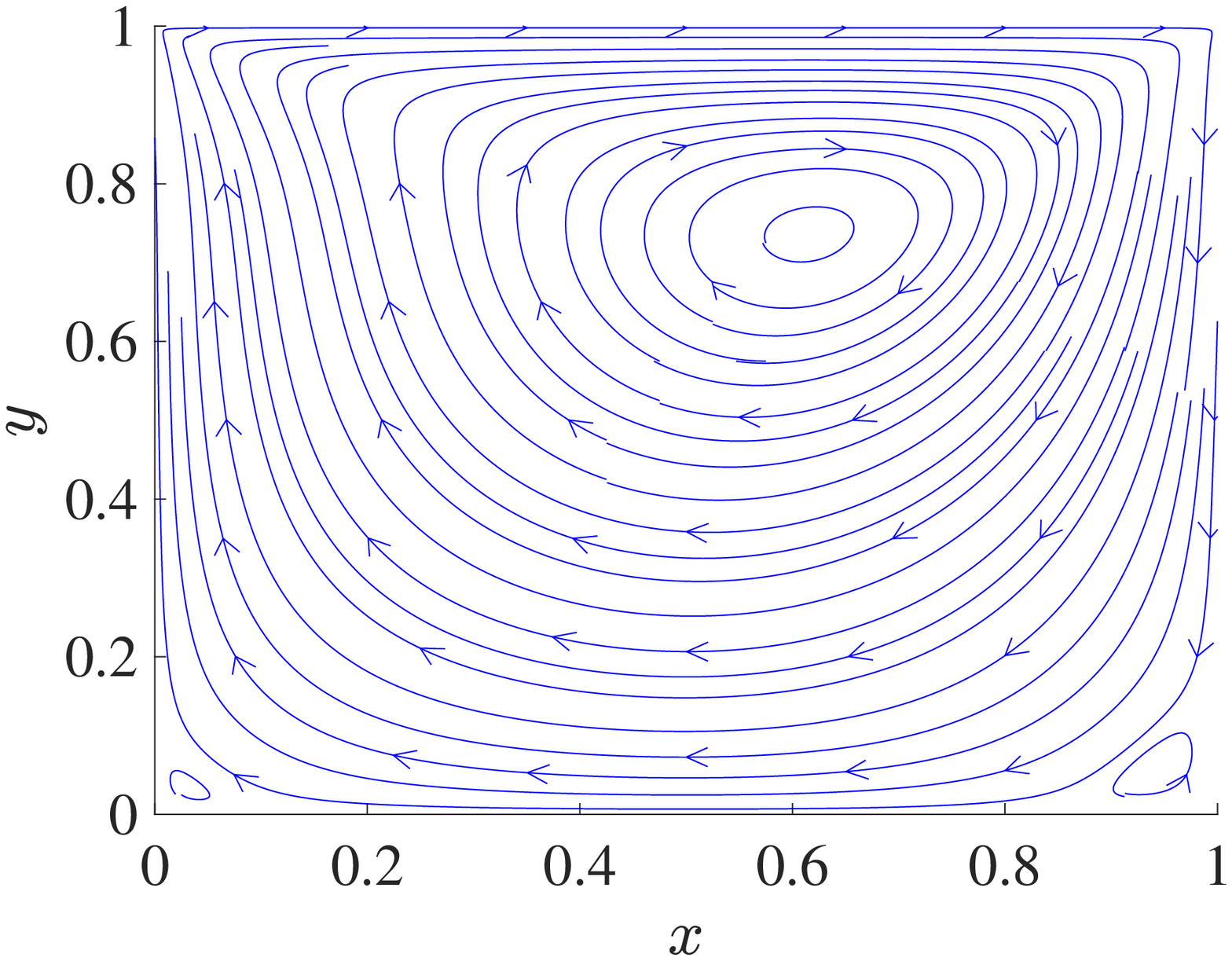}} \\
\subfloat[Re = 400 with 32$^2$ elements]{\label{sfig:Re400_32_sl_rot}\includegraphics[width=.45\textwidth]{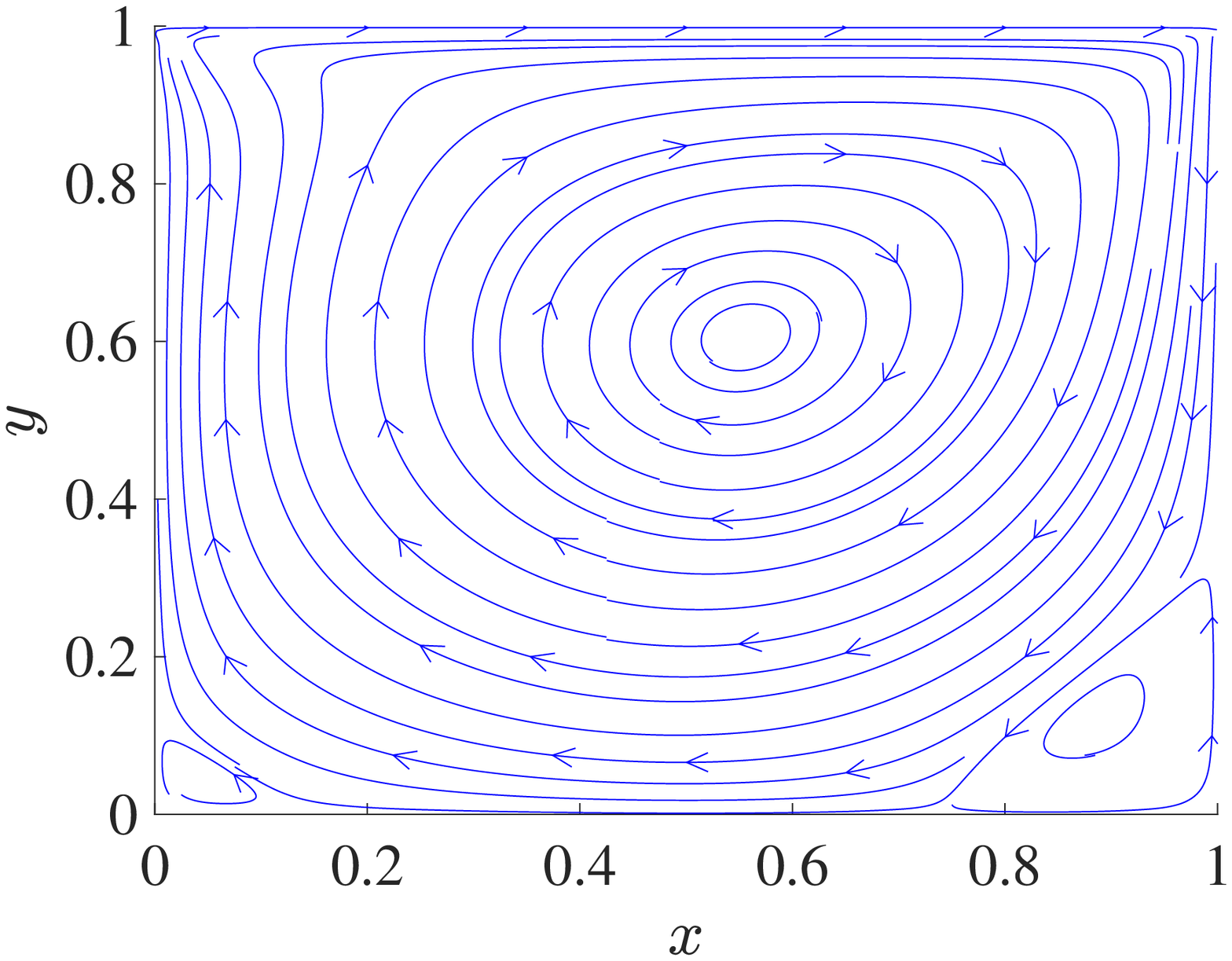}}\hfill
\subfloat[Re = 400 with 64$^2$ elements]{\label{sfig:Re400_64_sl_rot}\includegraphics[width=.45\textwidth]{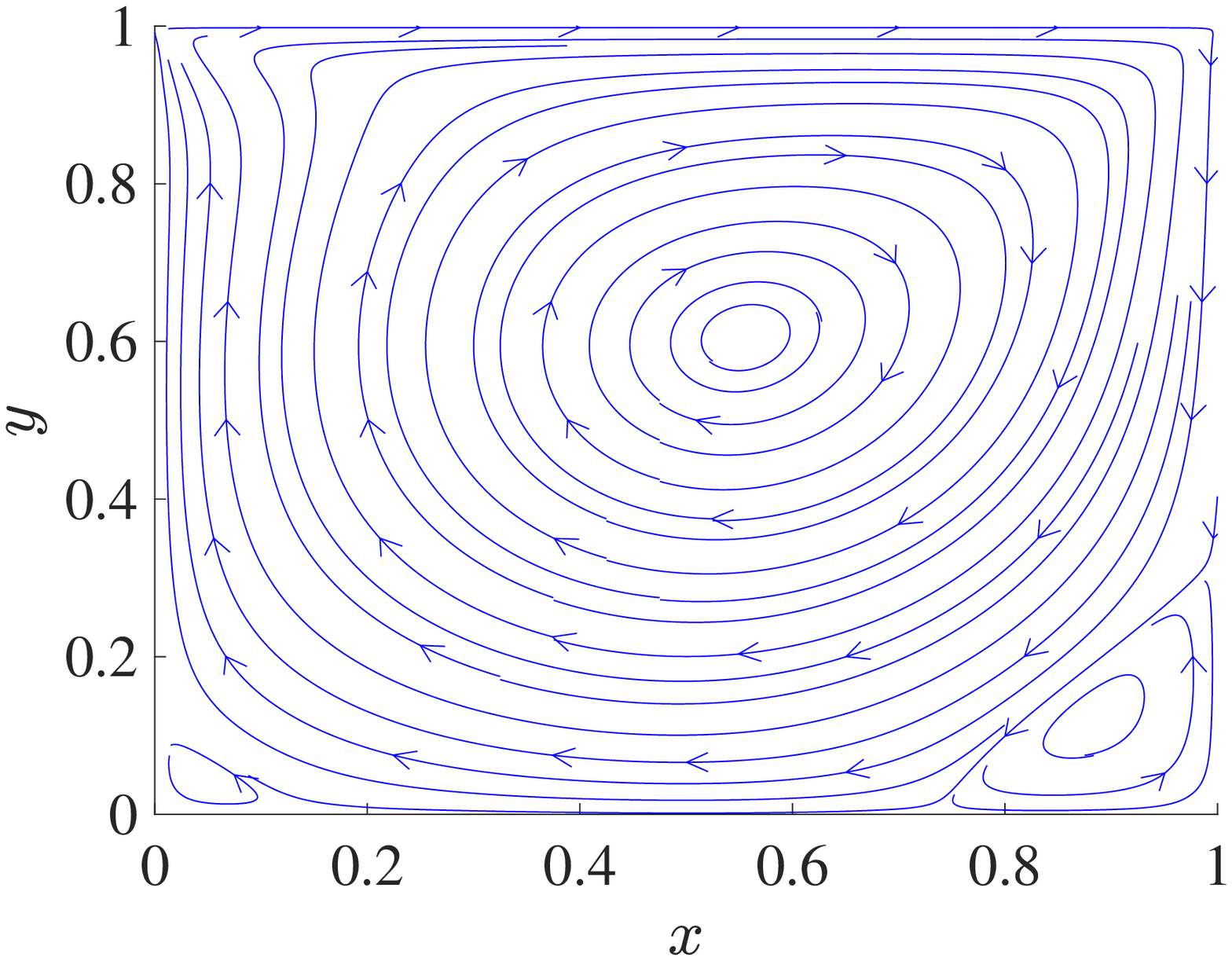}} \\
\subfloat[Re = 1000 with 32$^2$ elements]{\label{sfig:Re1000_32_sl_rot}\includegraphics[width=.45\textwidth]{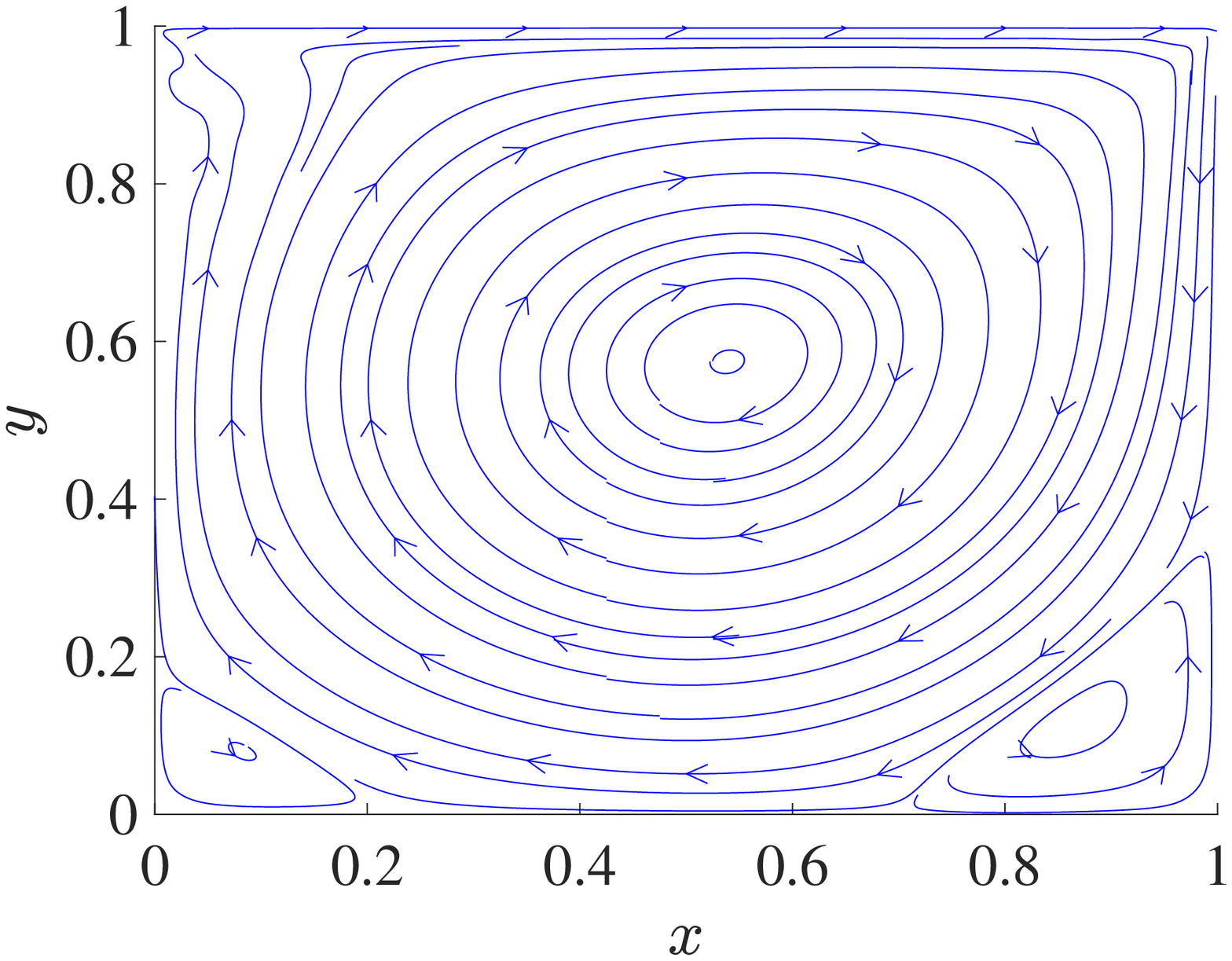}}\hfill
\subfloat[Re = 1000 with 64$^2$ elements]{\label{sfig:Re1000_64_sl_rot}\includegraphics[width=.45\textwidth]{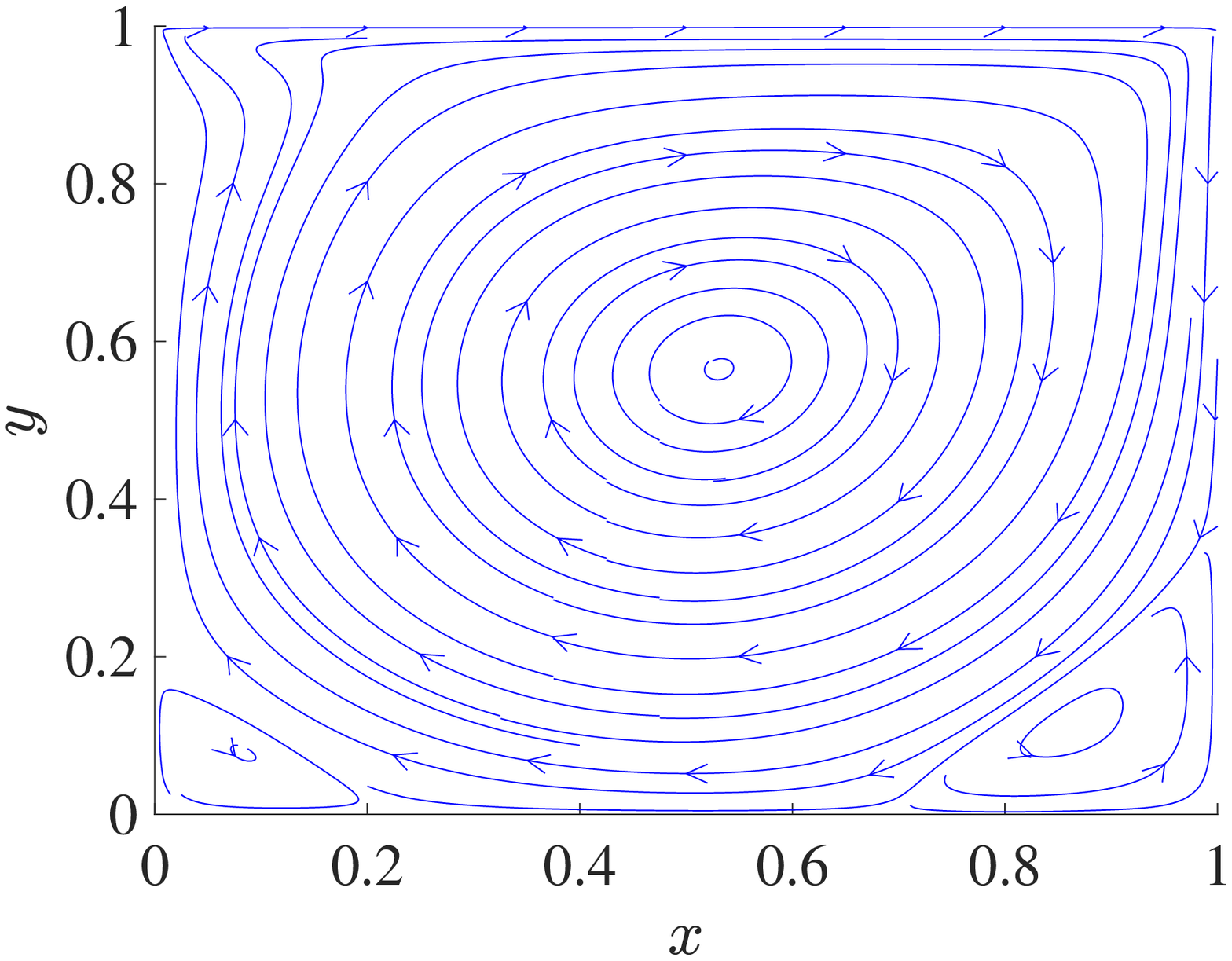}} \\
\caption{Streamlines of velocity for 2D lid-driven cavity flow using rotational form, $k$ = 5}
\label{fig:ns_2d_cavity_sl_rot}
\end{figure}

\begin{figure}
\centering
\subfloat[Re = 100 with 32$^2$ elements]{\label{sfig:Re100_32_rot}\includegraphics[width=.45\textwidth]{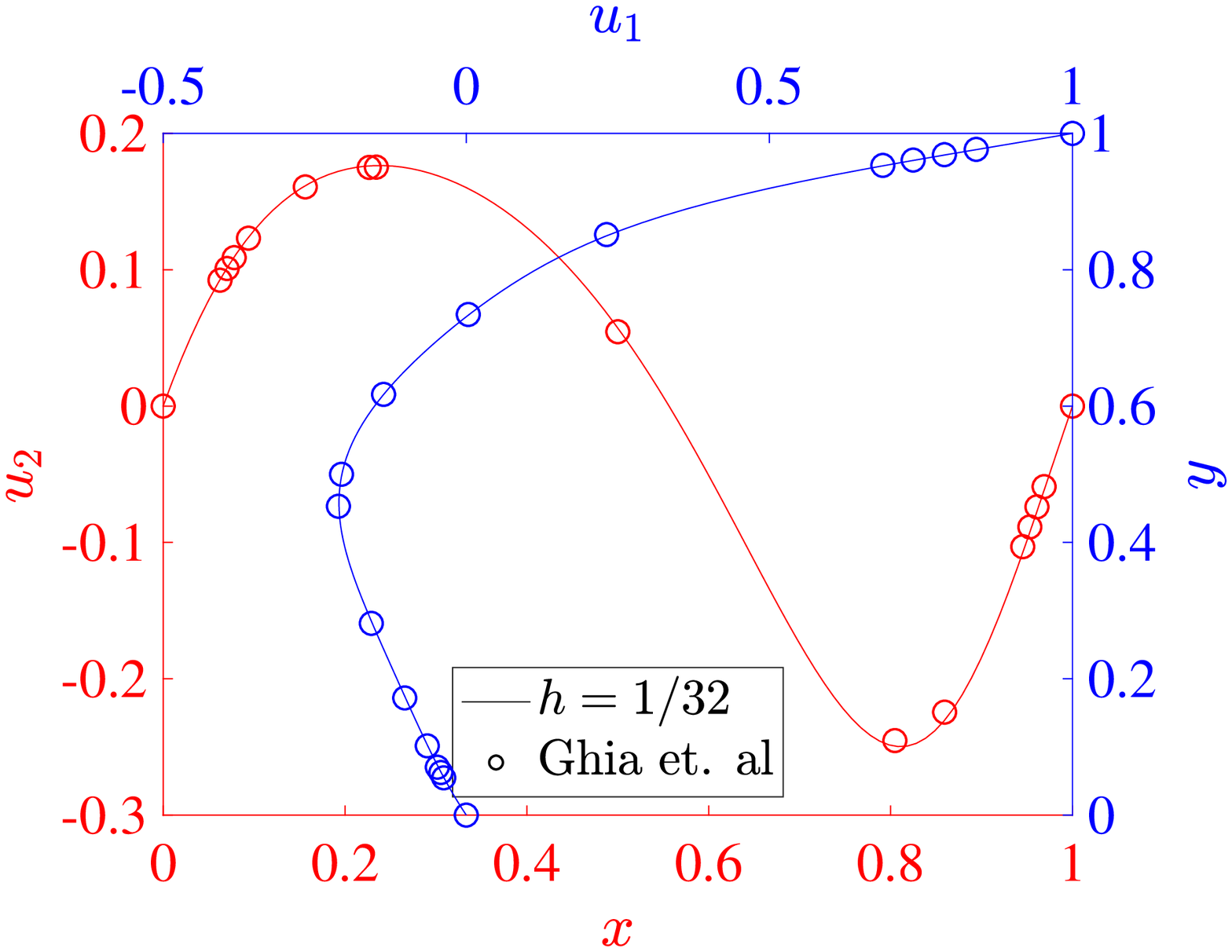}}\hfill
\subfloat[Re = 100 with 64$^2$ elements]{\label{sfig:Re100_64_rot}\includegraphics[width=.45\textwidth]{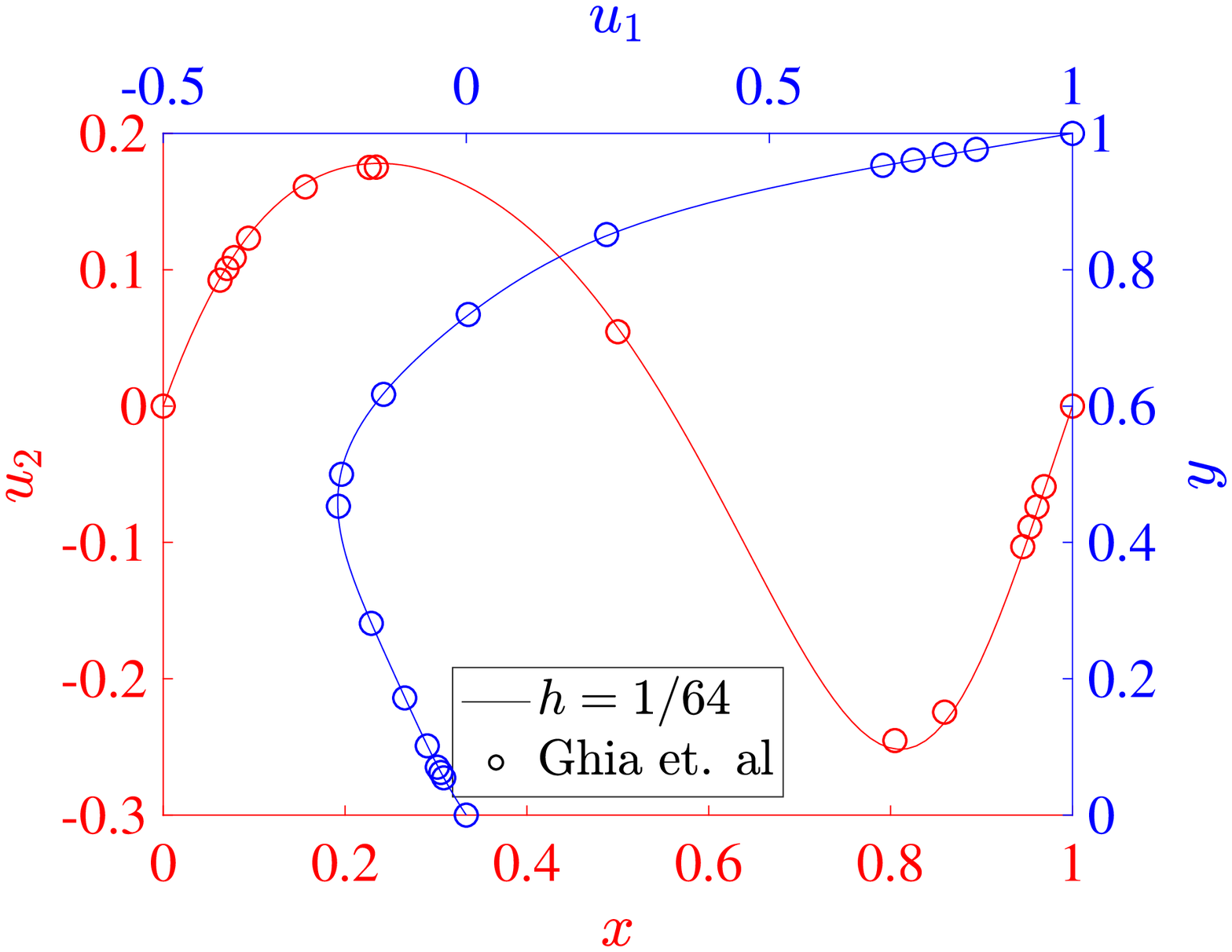}} \\
\subfloat[Re = 400 with 32$^2$ elements]{\label{sfig:Re400_32_rot}\includegraphics[width=.45\textwidth]{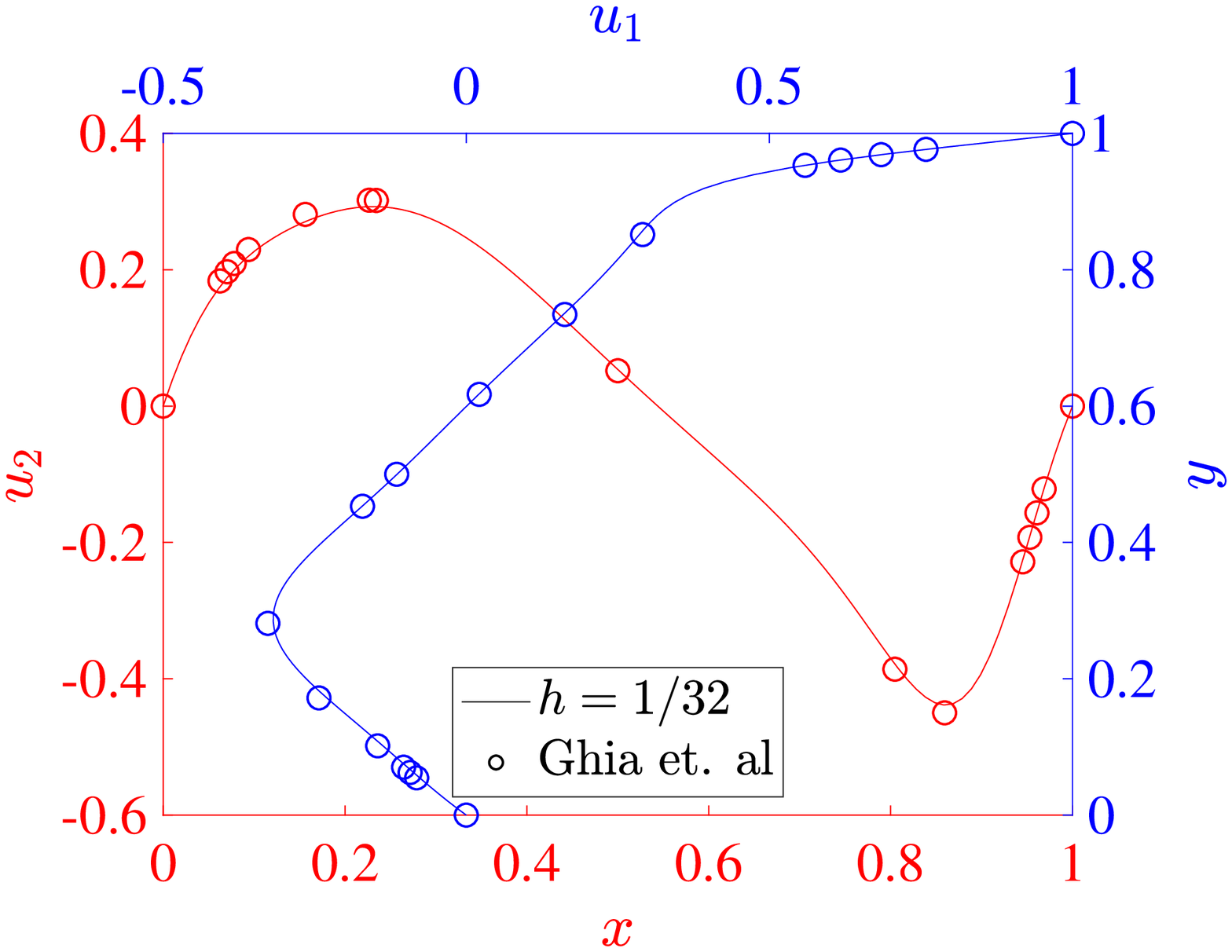}}\hfill
\subfloat[Re = 400 with 64$^2$ elements]{\label{sfig:Re400_64_rot}\includegraphics[width=.45\textwidth]{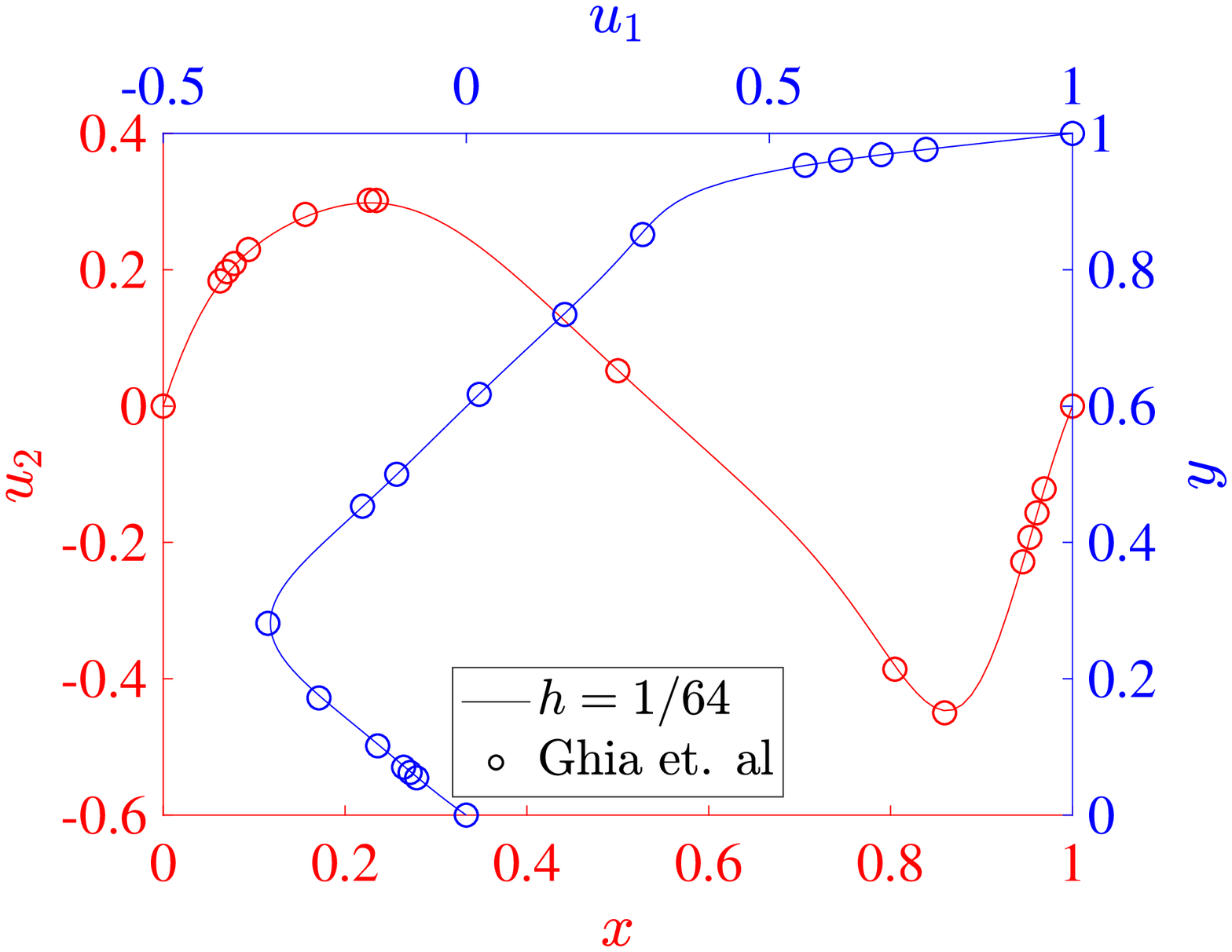}} \\
\subfloat[Re = 1000 with 32$^2$ elements]{\label{sfig:Re1000_32_rot}\includegraphics[width=.45\textwidth]{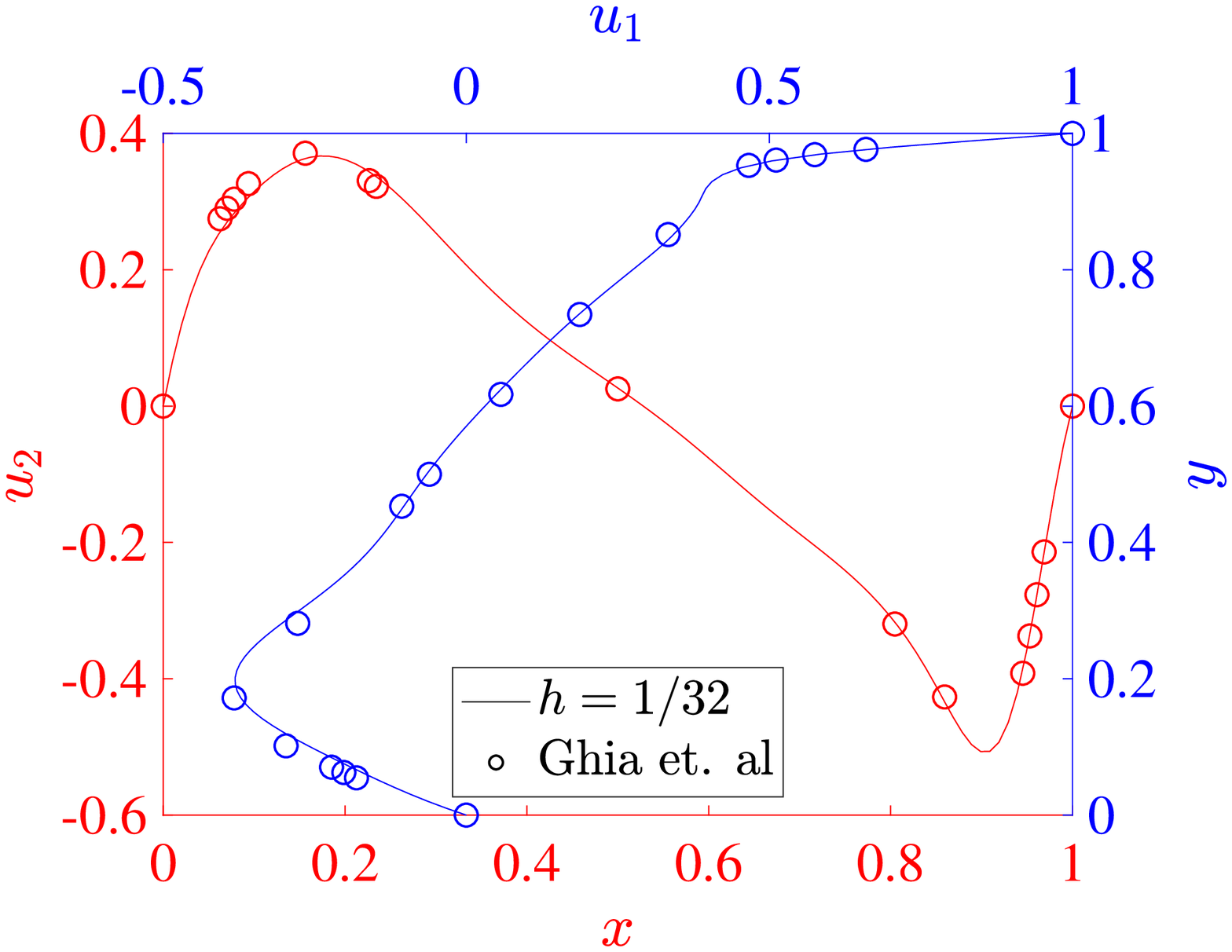}}\hfill
\subfloat[Re = 1000 with 64$^2$ elements]{\label{sfig:Re1000_64_rot}\includegraphics[width=.45\textwidth]{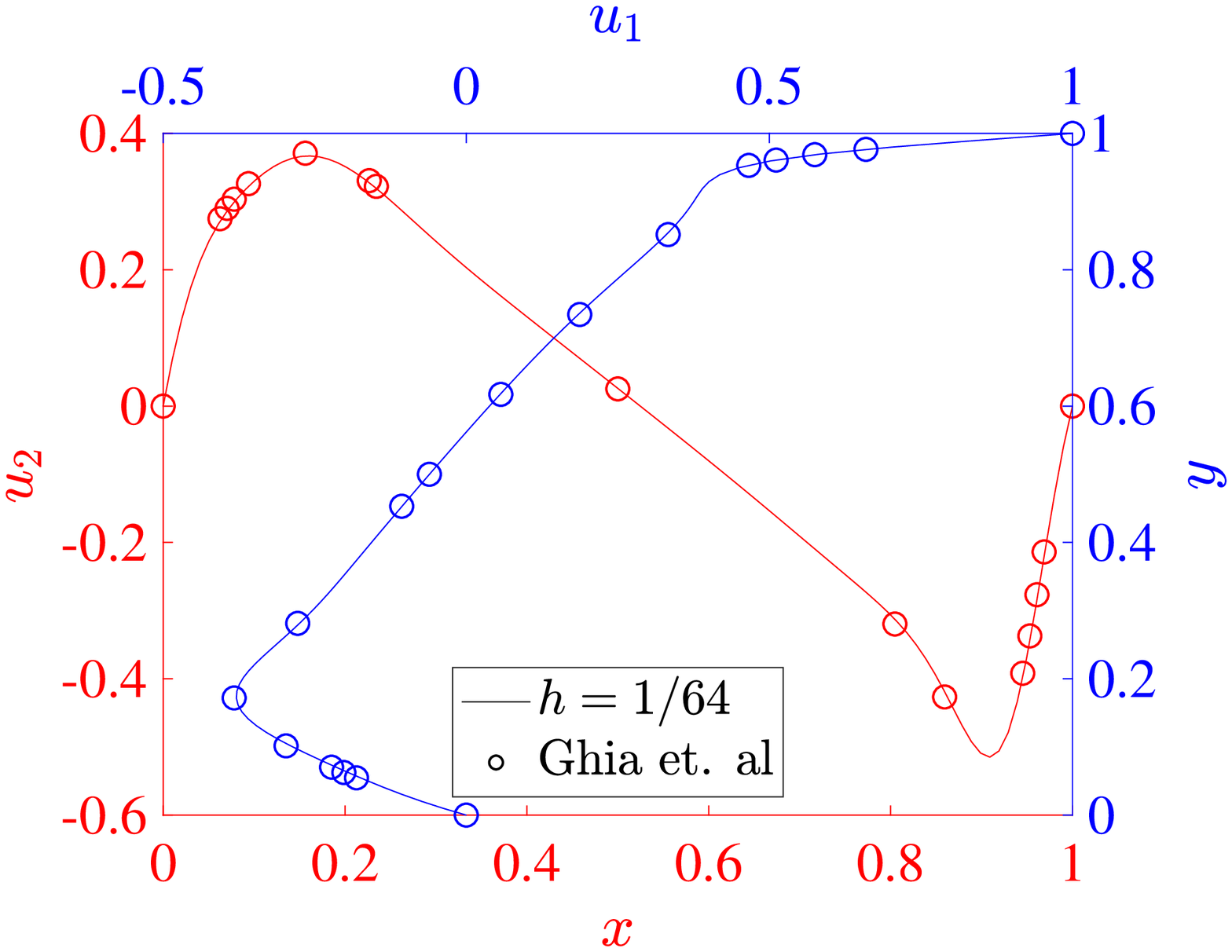}} \\
\caption{Centerline velocity profiles for 2D lid-driven cavity flow using rotational form, $k$ = 5}
\label{fig:ns_2d_cavity_rot}
\end{figure}

\begin{figure}
\centering
\subfloat[Re = 100 with 32$^2$ elements]{\label{sfig:Re100_32p8_sl_rot}\includegraphics[width=.45\textwidth]{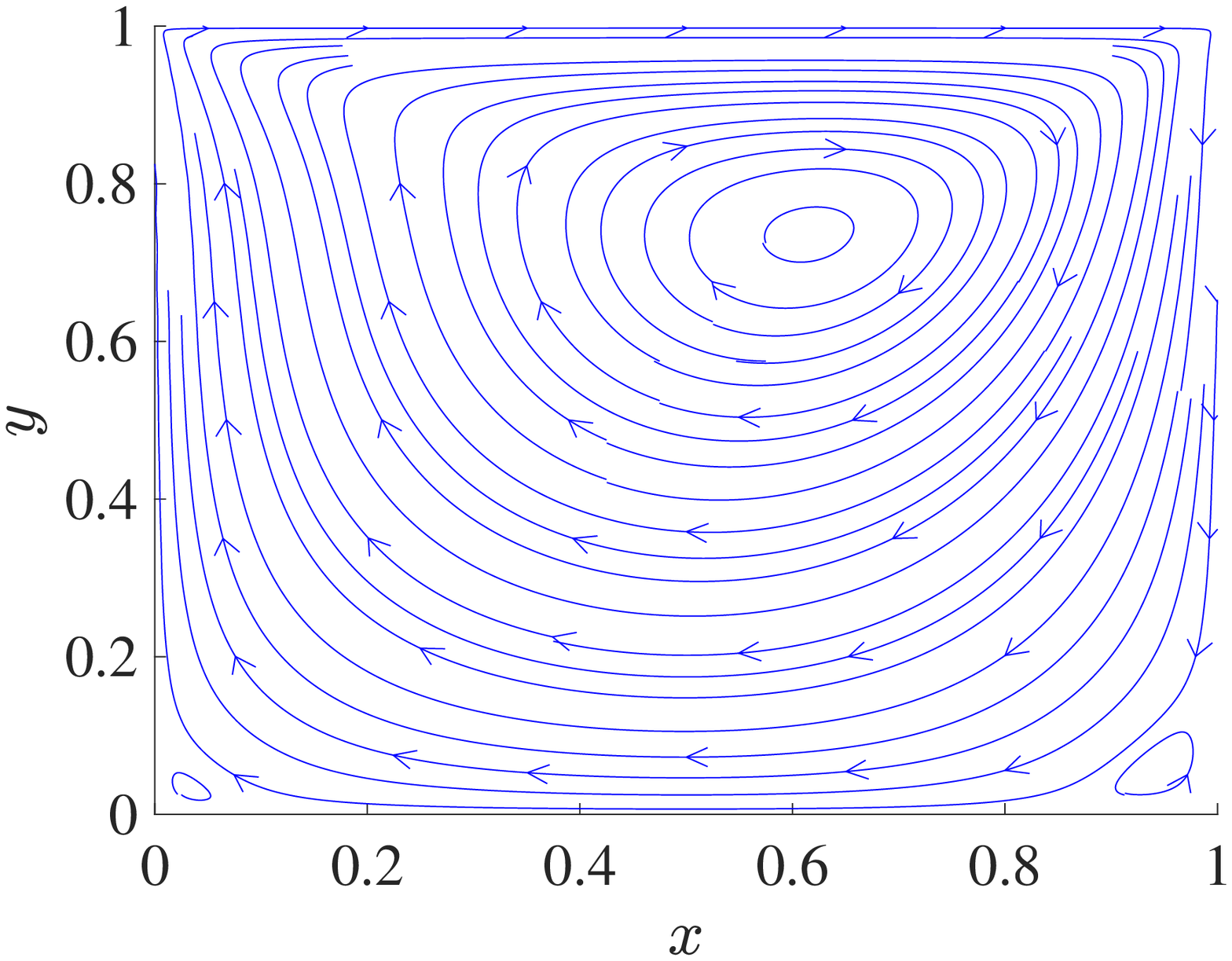}}\hfill
\subfloat[Re = 100 with 64$^2$ elements]{\label{sfig:Re100_64p8_sl_rot}\includegraphics[width=.45\textwidth]{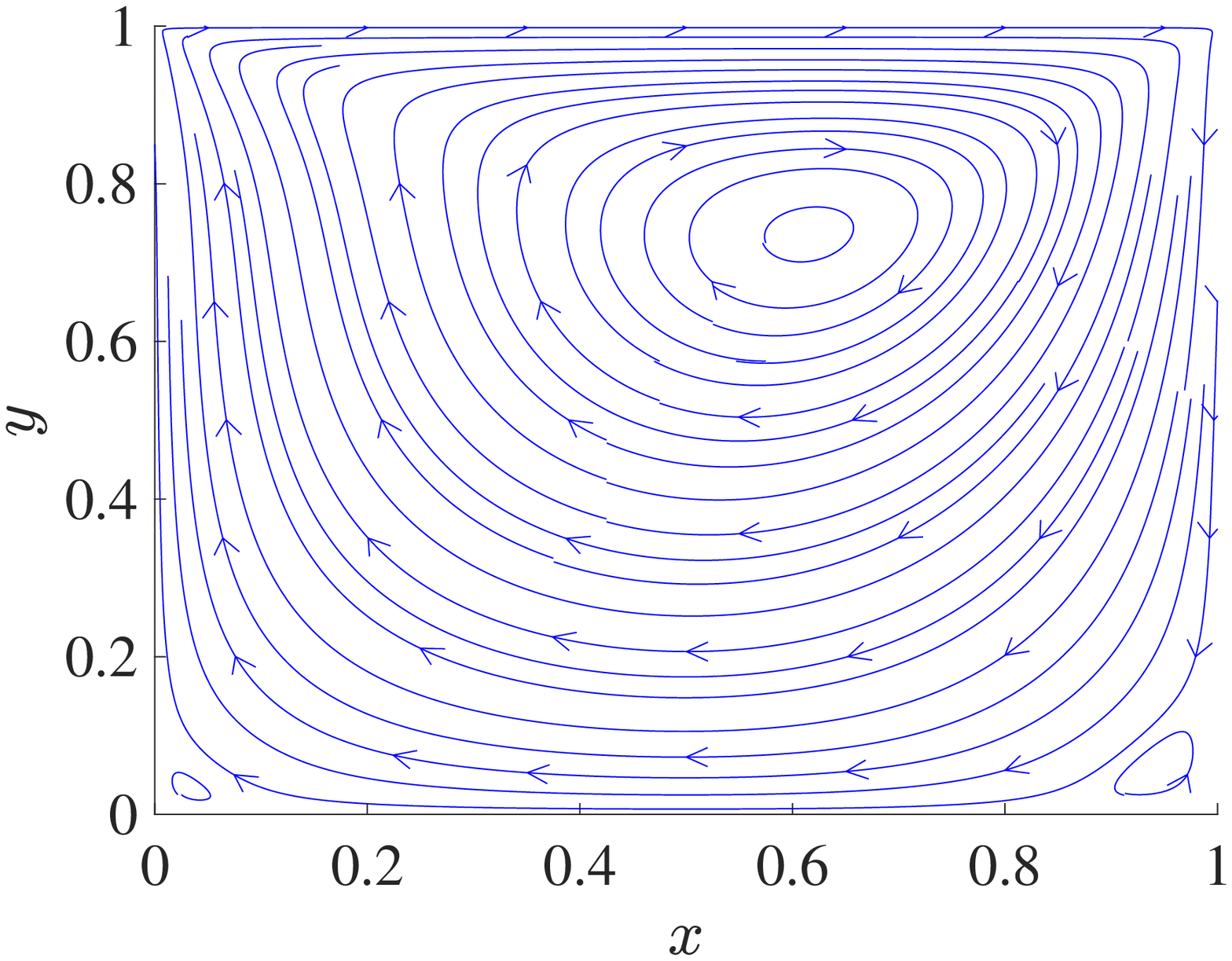}} \\
\subfloat[Re = 400 with 32$^2$ elements]{\label{sfig:Re400_32p8_sl_rot}\includegraphics[width=.45\textwidth]{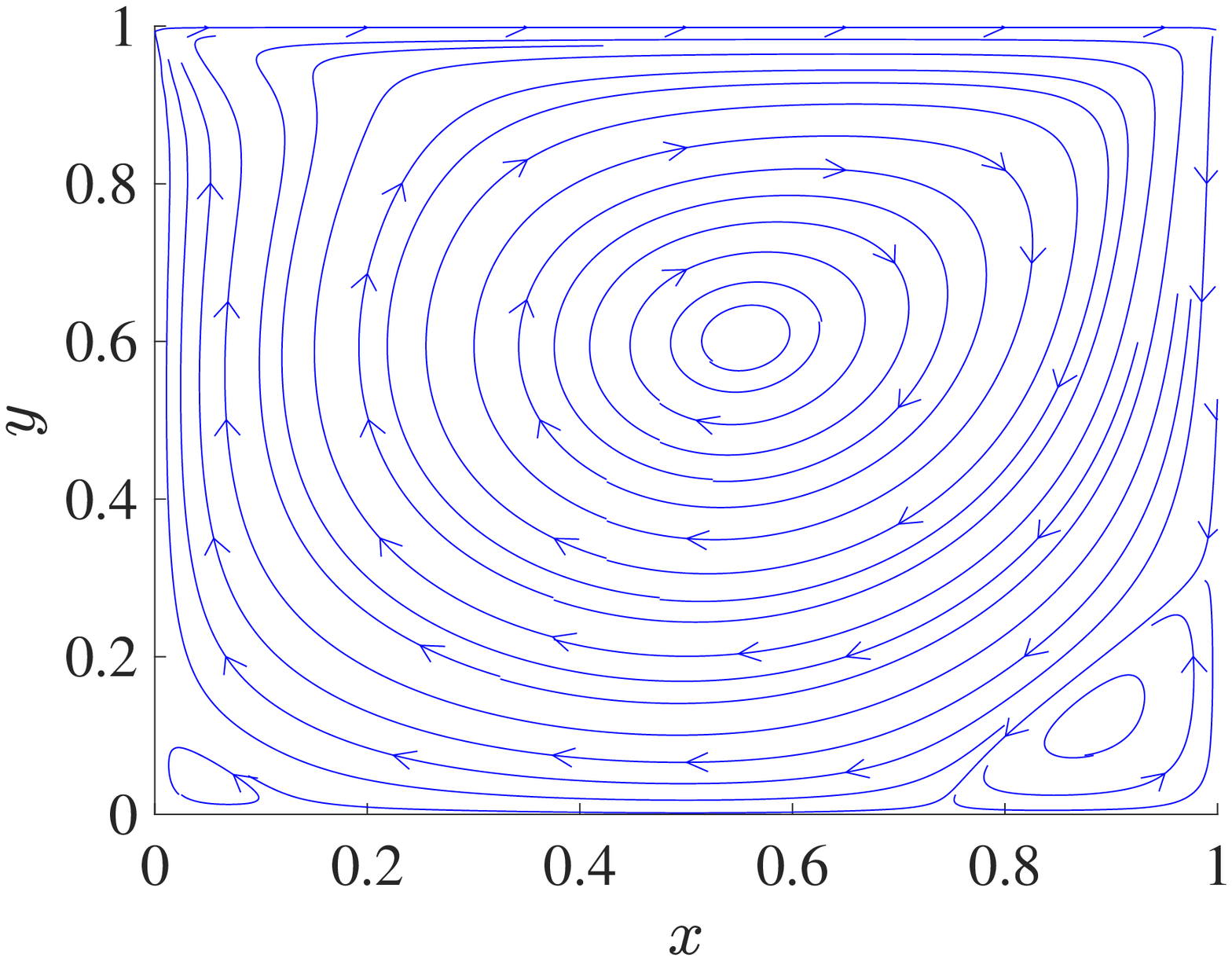}}\hfill
\subfloat[Re = 400 with 64$^2$ elements]{\label{sfig:Re400_64p8_sl_rot}\includegraphics[width=.45\textwidth]{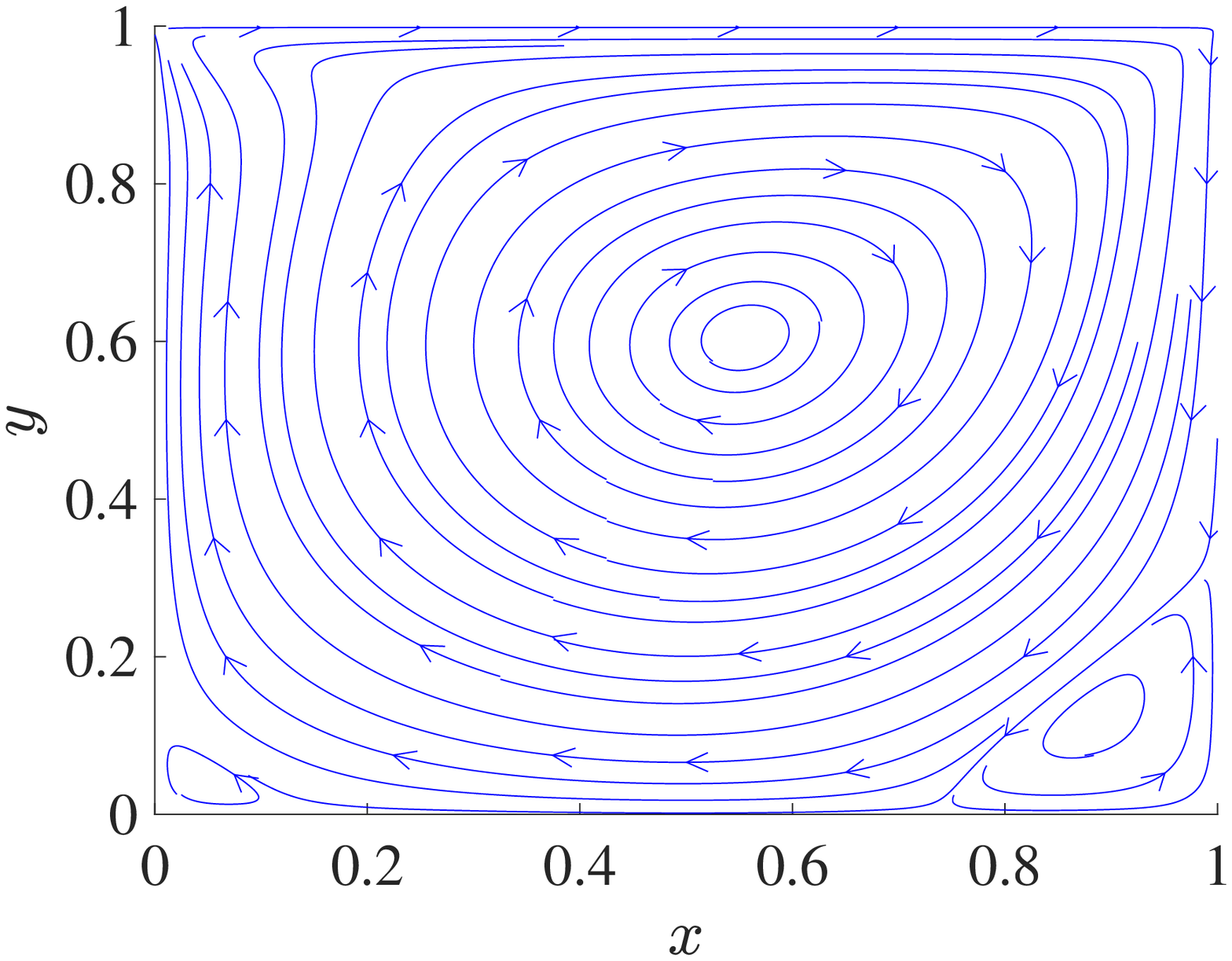}} \\
\subfloat[Re = 1000 with 32$^2$ elements]{\label{sfig:Re1000_32p8_sl_rot}\includegraphics[width=.45\textwidth]{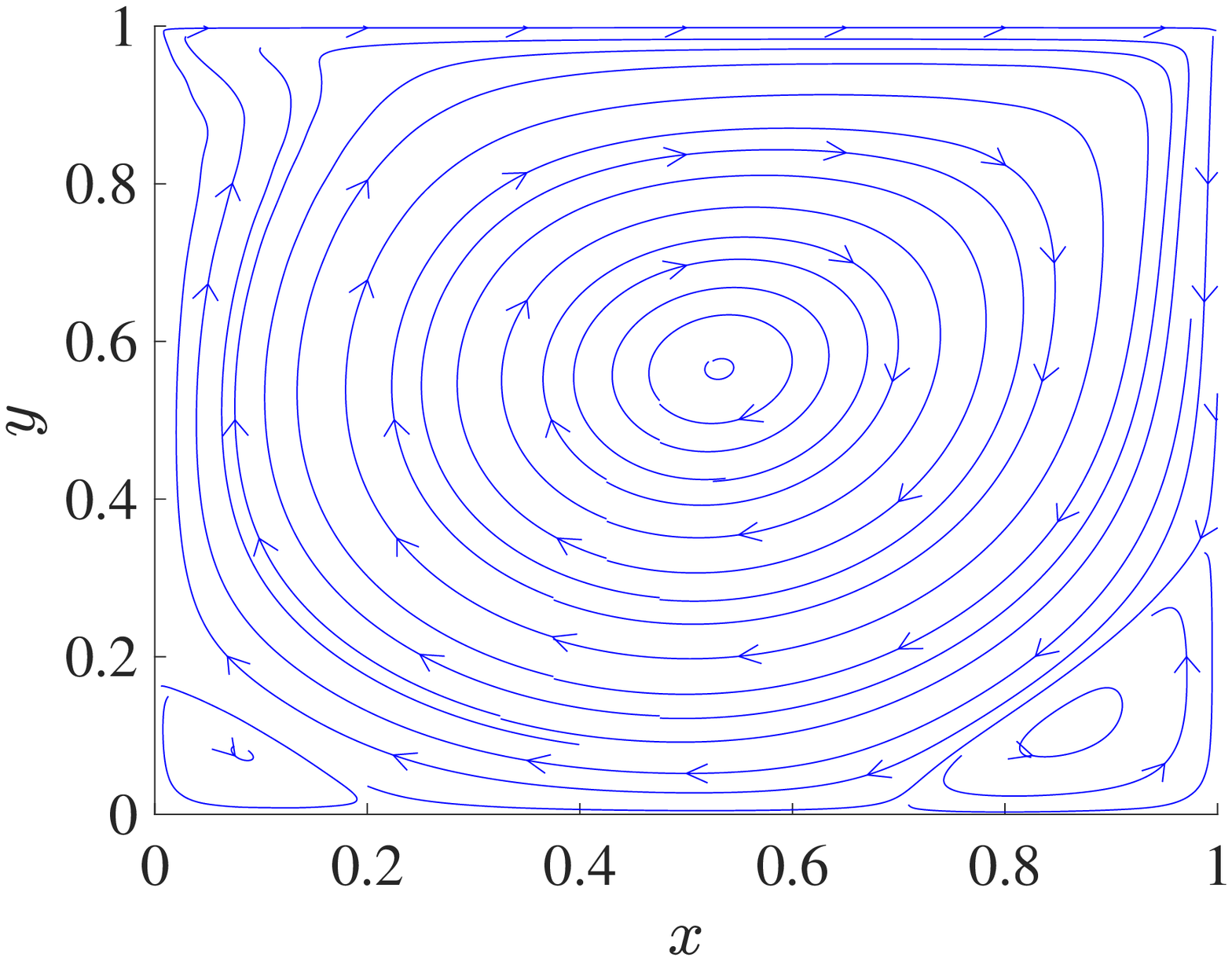}}\hfill
\subfloat[Re = 1000 with 64$^2$ elements]{\label{sfig:Re1000_64p8_sl_rot}\includegraphics[width=.45\textwidth]{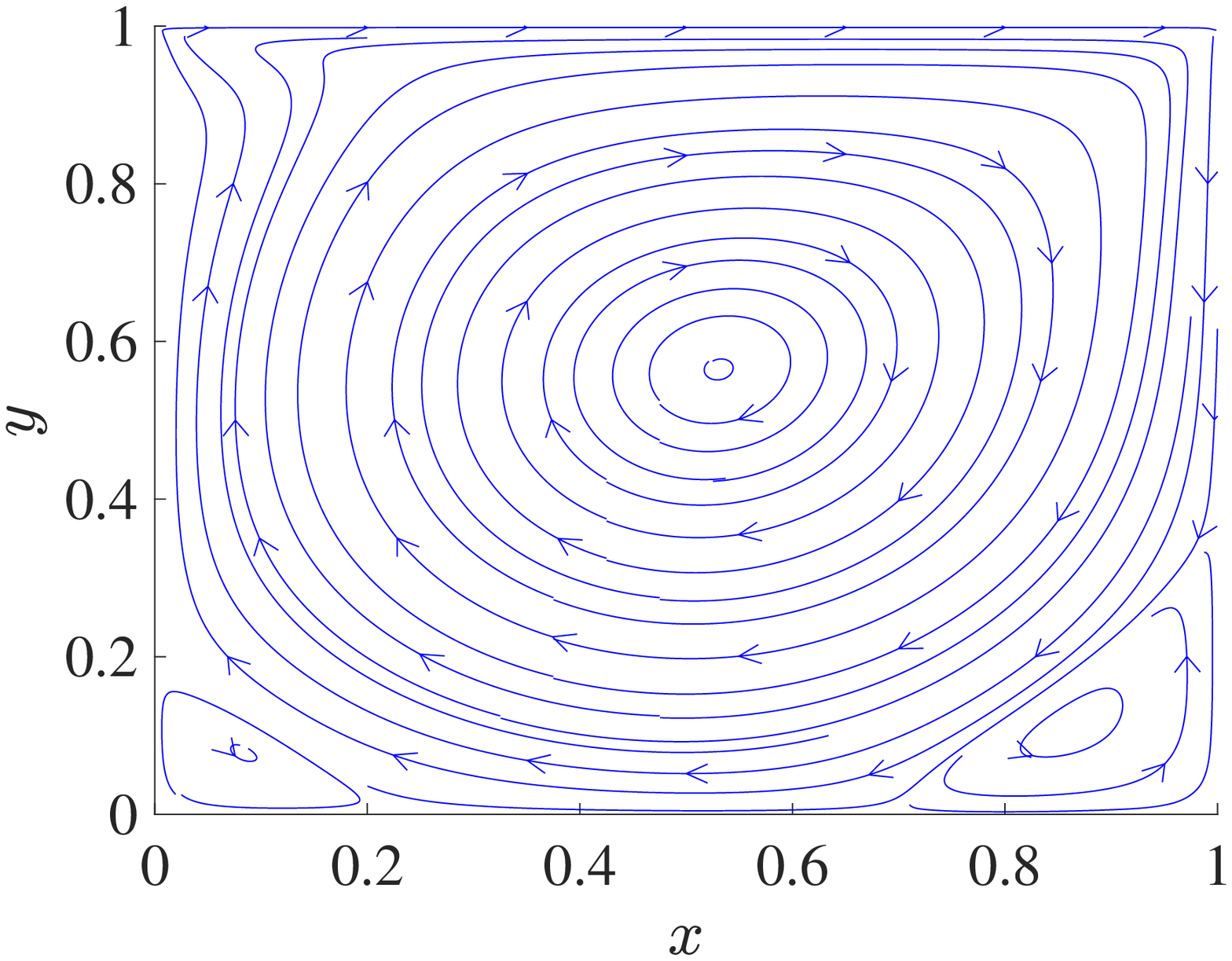}} \\
\caption{Streamlines of velocity for 2D lid-driven cavity flow using rotational form, $k$ = 8}
\label{fig:ns_2d_cavity_sl_rot_p8}
\end{figure}

\begin{figure}
\centering
\subfloat[Re = 100 with 32$^2$ elements]{\label{sfig:Re100_32p8_rot}\includegraphics[width=.45\textwidth]{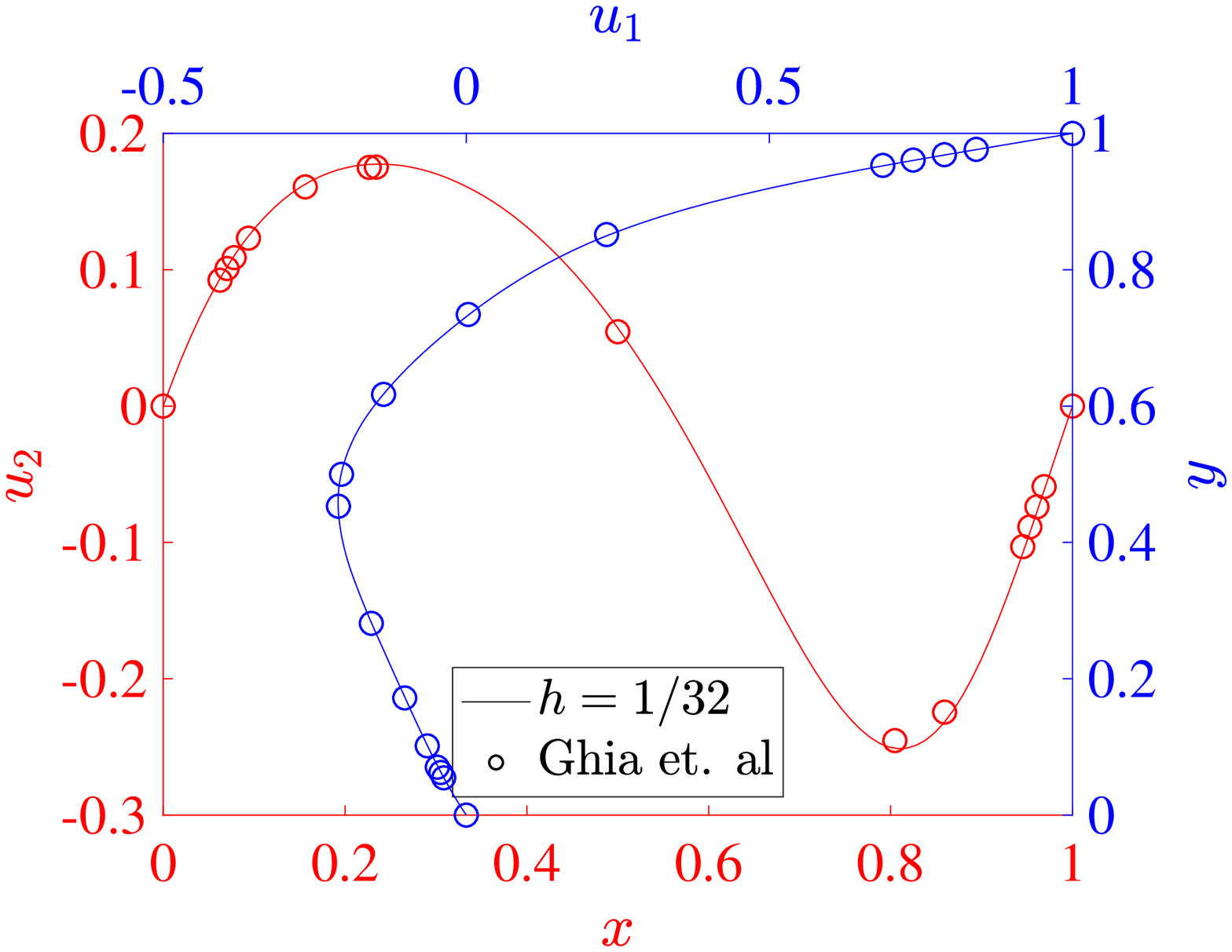}}\hfill
\subfloat[Re = 100 with 64$^2$ elements]{\label{sfig:Re100_64p8_rot}\includegraphics[width=.45\textwidth]{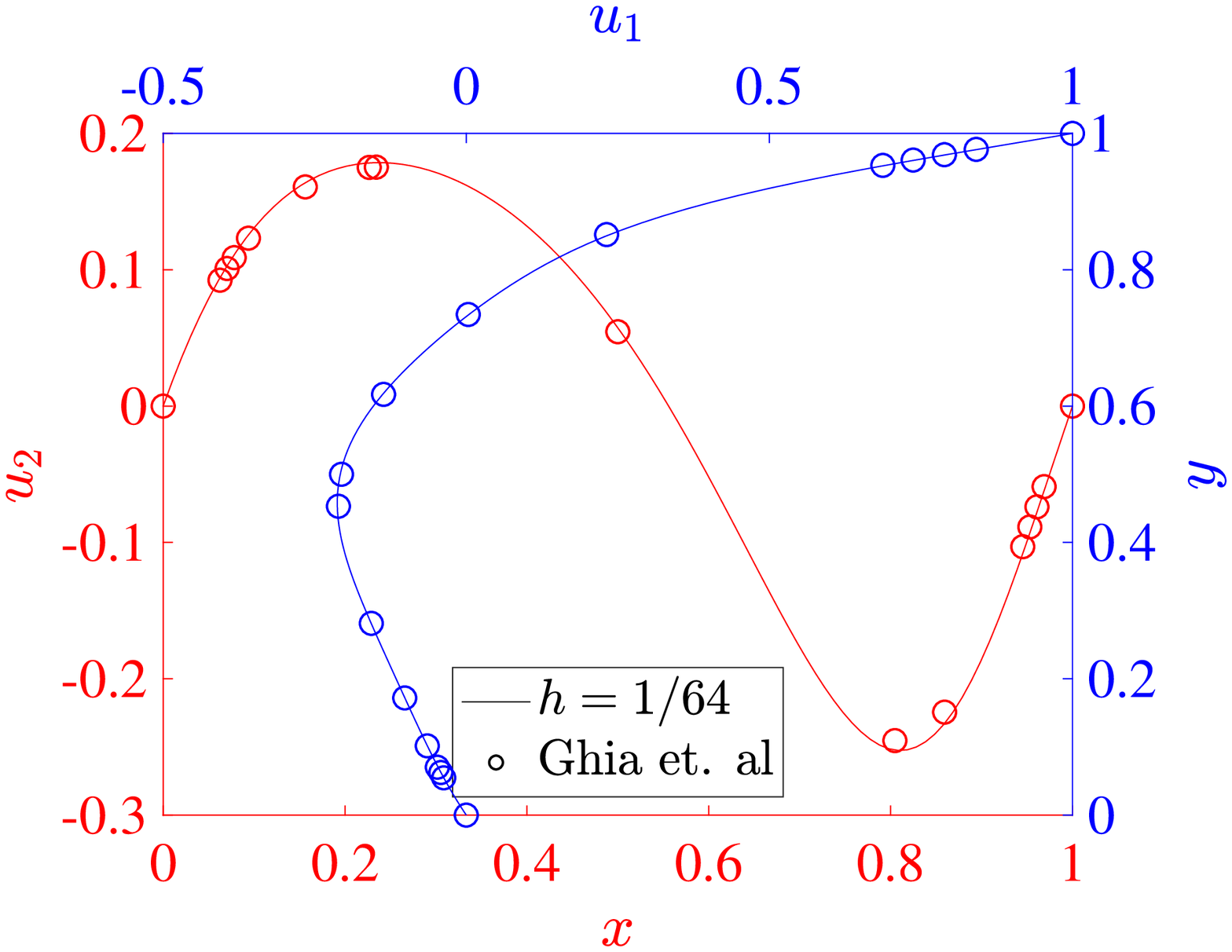}} \\
\subfloat[Re = 400 with 32$^2$ elements]{\label{sfig:Re400_32p8_rot}\includegraphics[width=.45\textwidth]{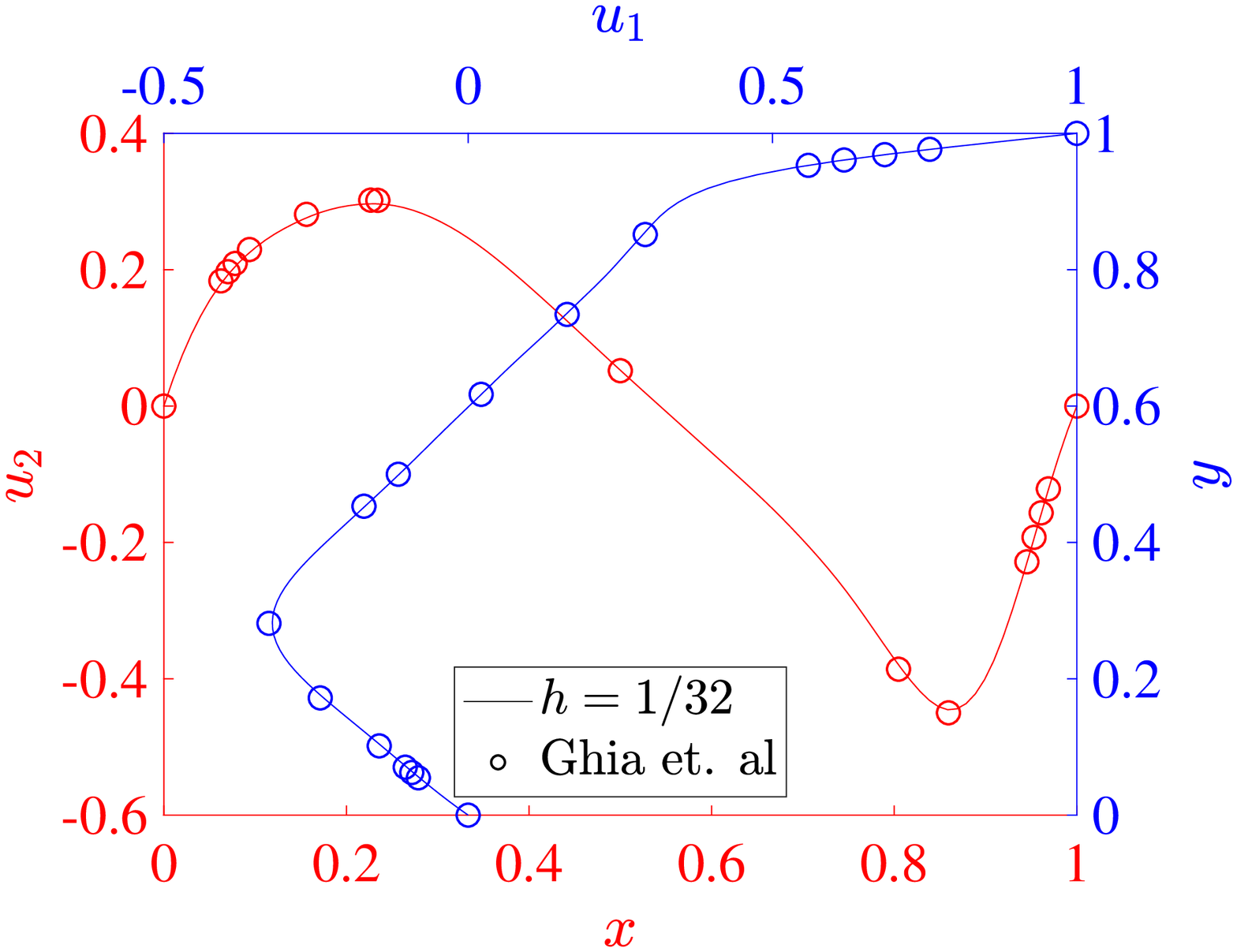}}\hfill
\subfloat[Re = 400 with 64$^2$ elements]{\label{sfig:Re400_64p8_rot}\includegraphics[width=.45\textwidth]{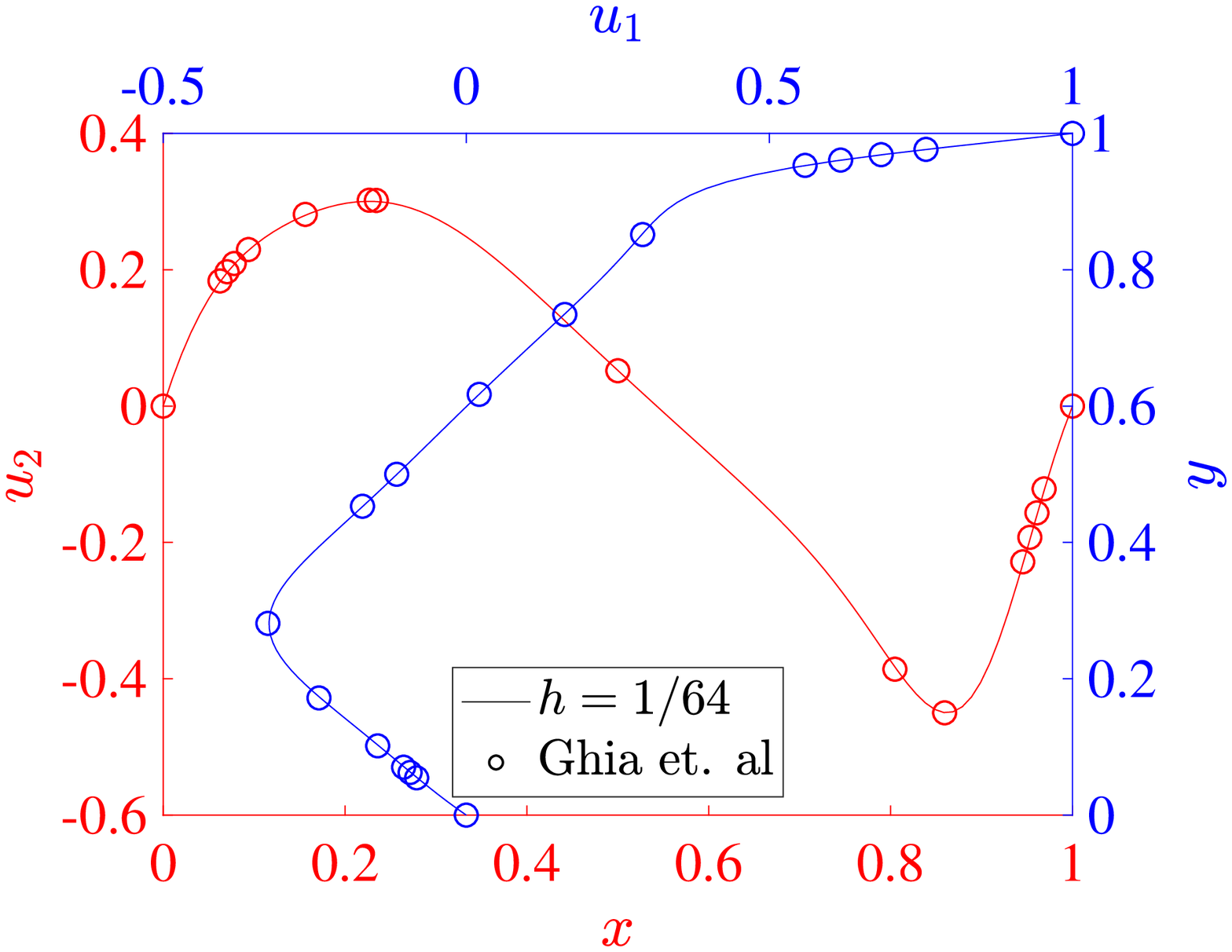}} \\
\subfloat[Re = 1000 with 32$^2$ elements]{\label{sfig:Re1000_32p8_rot}\includegraphics[width=.45\textwidth]{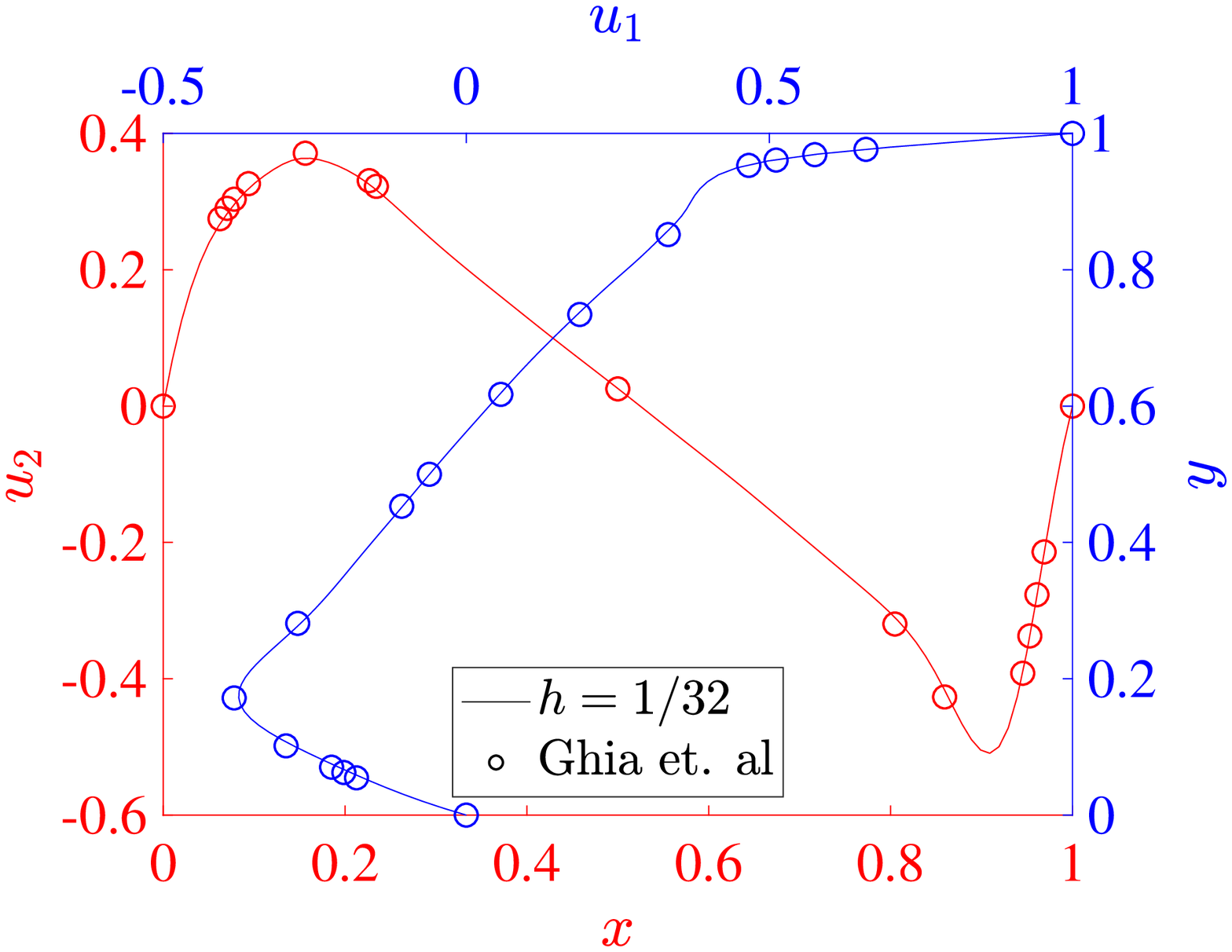}}\hfill
\subfloat[Re = 1000 with 64$^2$ elements]{\label{sfig:Re1000_64p8_rot}\includegraphics[width=.45\textwidth]{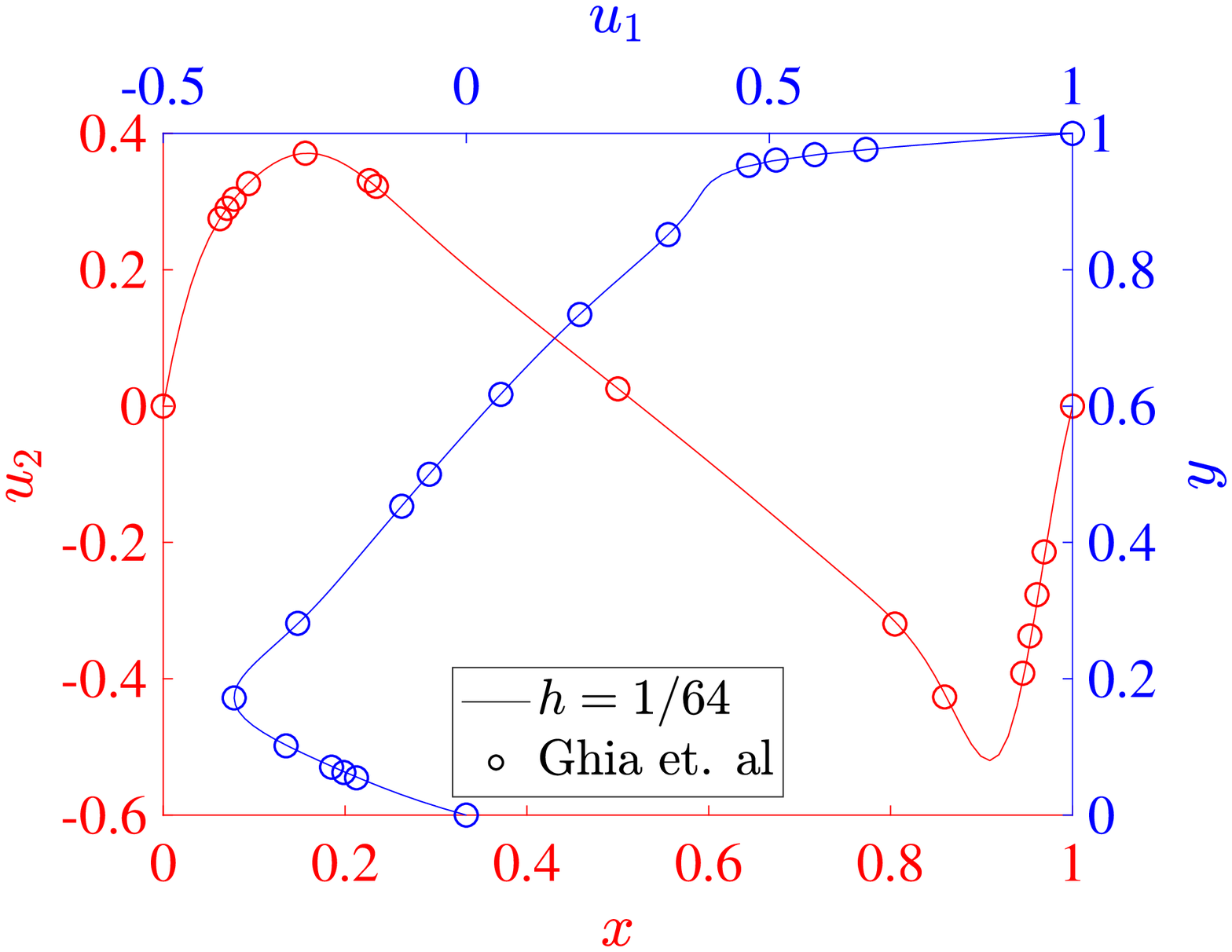}} \\
\caption{Centerline velocity profiles for 2D lid-driven cavity flow using rotational form, $k$ = 8}
\label{fig:ns_2d_cavity_rot_p8}
\end{figure}

We repeat the studies with the rotational formulation collocation scheme. Figure \ref{fig:ns_2d_cavity_sl_rot} details the flow streamlines on the same $32^2$ and $64^2$ element stretched meshes with $k = 5$, which are very similar to those produced with the two-field scheme. Figure \ref{fig:ns_2d_cavity_rot} shows the centerline velocities for these cases. It seems that the rotational form scheme produces very similar results to the velocity-pressure scheme across the Reynolds number spectrum. Finally, Figures \ref{fig:ns_2d_cavity_sl_rot_p8} and \ref{fig:ns_2d_cavity_rot_p8} show these same experiments with the basis degree elevated to $k = 8$, confirming that the results are also converging with polynomial degree for the rotational scheme.

\FloatBarrier

\section{Conclusions}

In this work we have developed stabilized, isogeometric collocation methods for scalar transport and incompressible flow problems. Inspired by residual-based stabilizations for Galerkin discretizations, we extend SUPG and PSPG techniques to the spline collocation setting. Our numerical tests demonstrate that these stabilizations effectively remove spurious oscillations in the solution field (in the case of scalar transport) or the velocity field (in the case of incompressible flow) caused by advection instabilities, as well as spurious oscillations in the pressure field caused by the pressure instabilities associated with equal-order interpolations. 

While we have focused on problems posed on simple domains in this work, stabilized collocation schemes on more complex domains could be generated by utilizing NURBS instead of B-splines, multi-patch discretizations, and possibly immersed collocation techniques \cite{torre2023immersed}. Another interesting future direction is development of a full, unsteady, turbulent, incompressible flow solver based on these techniques. Due to the efficient nature and high resolving power of spline collocation schemes, a solver like this may be much more efficient in performing Direct Numerical Simulations and Large Eddy Simulations of turbulent flows than existing techniques. Moreover, the relative advantages and disadvantages of the rotational form of the equations should be studied in more detail, including more study of the definition of the stabilization constants. Further investigating the connection of the stabilized methods developed here and possible variational multiscale collocation techniques would also be interesting, especially given the connection between variational multiscale techniques and turbulence modeling. 

\section*{Acknowledgements}

This material is based upon work supported by the National Science Foundation Graduate Research Fellowship Program under Grant No. DGE-1656518. Any opinions, findings, and conclusions or recommendations expressed in this material are those of the authors and do not necessarily reflect the views of the National Science Foundation.

\FloatBarrier

\bibliographystyle{elsarticle-num}
\bibliography{ref}

\end{document}